\documentclass[twocolumn]{revtex4-1}
\pdfoutput=1

\usepackage{amsmath}
\usepackage{amssymb}
\usepackage{amsfonts}
\usepackage{graphicx}

\begin{document}

\title{ 
Theorems of Euclidean Geometry through Calculus
}
\author{Martin Buysse}
\email[]{martin.buysse@uclouvain.be}
\affiliation{{Facult\'e d'architecture, d'ing\'enierie architecturale, d'urbanisme -- LOCI}, {UCLouvain}}

\begin{abstract} 
\noindent We re-derive Thales, Pythagoras, Apollonius, Stewart, Heron, al Kashi, de Gua, Terquem, Ptolemy, Brahmagupta and Euler's theorems as well as the inscribed angle theorem, the law of sines, the circumradius, inradius and some angle bisector formulae, by assuming the existence of an unknown relation between the geometric quantities at stake, observing how the relation behaves under small deviations of those quantities, and naturally establishing differential equations that we integrate out. Applying the general solution to some specific situation gives a particular solution corresponding to the expected theorem. We also establish an equivalence between a polynomial equation and a set of partial differential equations. We finally comment on a differential equation which arises after a small scale transformation and should concern all relations between metric quantities.
\end{abstract}

\maketitle


\section{Thales of Miletus}
\label{Thales}

\noindent Imagine that Newton was born before Thales. When considering a triangle with two sides of lengths $x$ and $y$, he could have fantasized about moving the third side parallel to itself and thought: "Well, I am not an ancient Greek geometer but I am rather good in calculus and I feel there might be some connection between the way $x$ and $y$ vary in such circumstances." 
\parbox{43.5mm}{\vspace{-1mm}
\includegraphics[width=43.5mm,keepaspectratio]{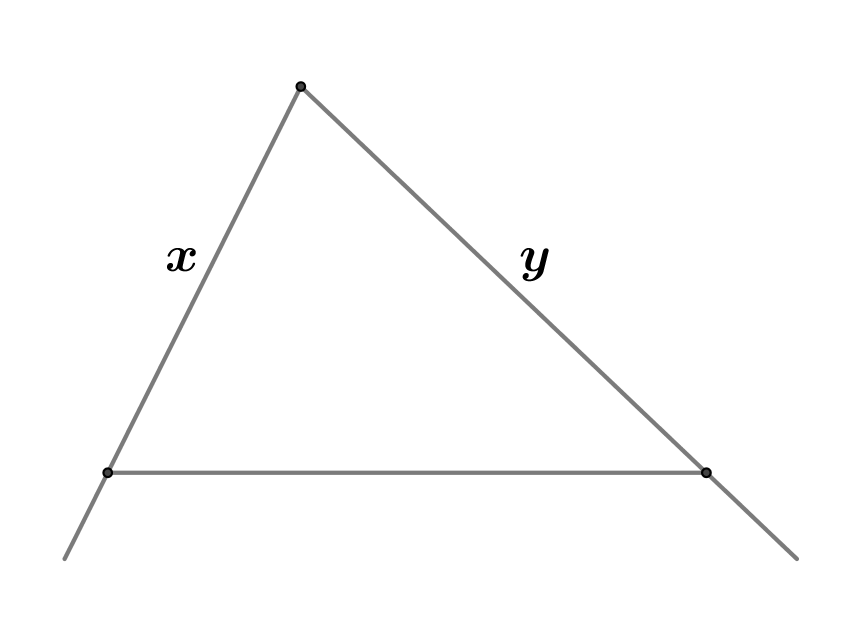}
\vspace*{-8mm}

\includegraphics[width=43.5mm,keepaspectratio]{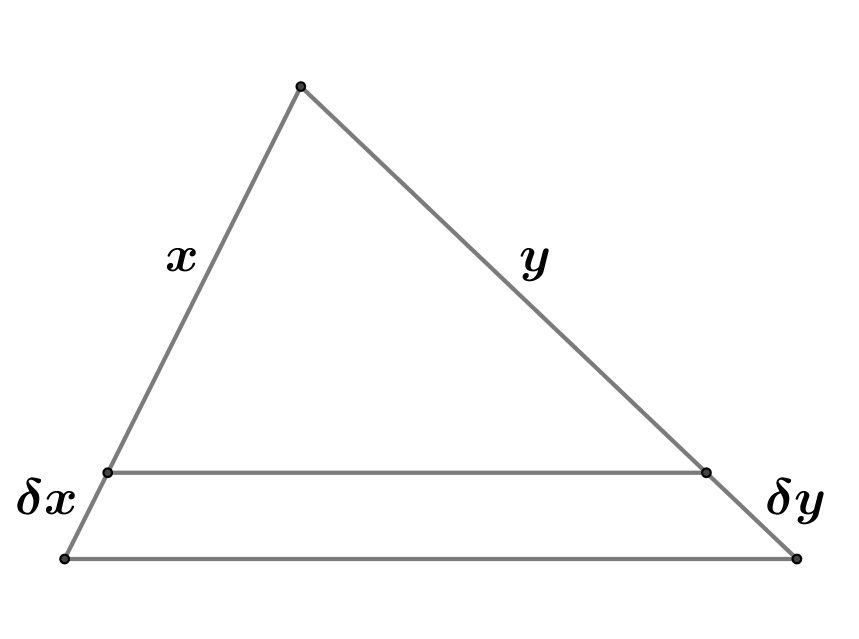}
\vspace*{-5mm}

\includegraphics[width=43.5mm,keepaspectratio]{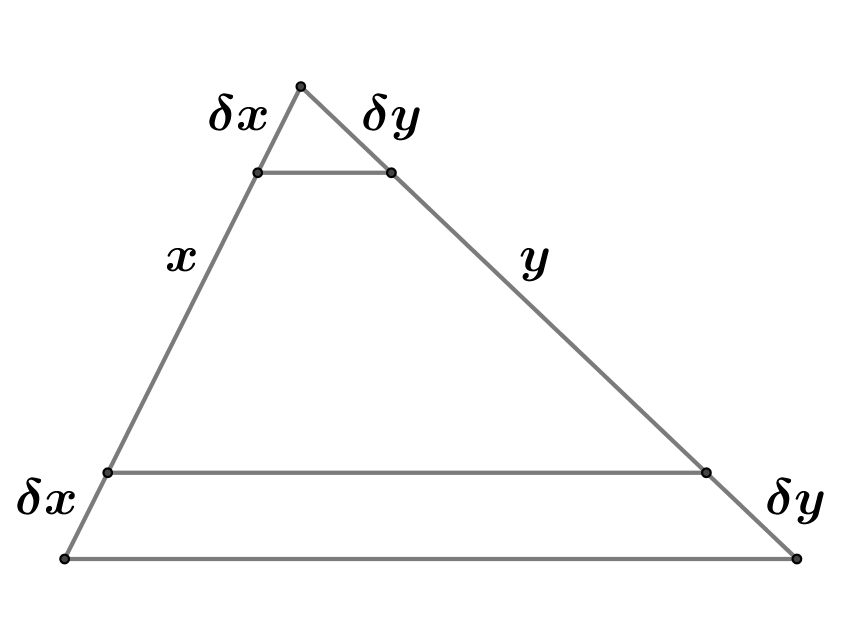}
}%
\parbox{43mm}{\vspace{2mm}
He would have materialized his suspicion in a function 
\begin{equation}
y=y(x)
\label{ThalesFunction_0}
\end{equation}

\noindent connecting $x$ and $y$ whatever the position of the third side, as long as it is moved parallel to itself. In particular, after a slight displacement resulting in small deviations $\delta x$ and $\delta y$, he would have had at first order
\begin{equation}
{\delta y}=y'(x){\delta x}
\label{ThalesDerivative_0}
\end{equation}
But the lengths $\delta x$ and $\delta y$ of the small added segments must themselves obey equation (\ref{ThalesFunction_0}), that is
\begin{equation}
\delta y=y(\delta x)
\label{ThalesFunction_2}
\end{equation}
}
\vspace{-1mm}

\noindent To see it, translate those segments to the $(x,y)$ vertex.
Developing the right-hand-side member of equation (\ref{ThalesFunction_2}) at first order and noticing that $y(0)=0$, we have
\vspace{-1mm}
\begin{equation}
{\delta y}=y'(0){\delta x}
\label{ThalesDerivative_1}
\end{equation}

\vspace{-1mm}
\noindent which, compared to eq. (\ref{ThalesDerivative_0}), implies that $y'(x)$ is constant. Integrating $y'=k$, $k$ being a positive constant since $y(x)$ is an increasing and smooth function, gives the Thales theorem \cite[Book VI, Prop. II]{Euclid:300BC}
\begin{equation}
y(x)=kx
\vspace{5mm}
\label{ThalesTheorem}
\end{equation}

\section{Pythagoras of Samos}
\label{Pythagoras}

\noindent If he was born before Thales, Newton was born before Pythagoras too, so that we do not have to make any further unlikely hypothesis. Imagine that driven by his success in suspecting the existence of a Greek theorem, he moved to consider a right triangle of legs of lengths $x$ and $y$ and of hypotenuse of length $z$. 
\vspace{-3mm}

\noindent\hspace{-1mm}
\parbox{43mm}{\vspace{2mm} He might have been tempted to speculate about the link, if any, between $x$, $y$ and $z$ in every right triangle. And again, as calculus master, he could have postulated that

}%
\parbox{43.5mm}{
\includegraphics[width=43.5mm,keepaspectratio]{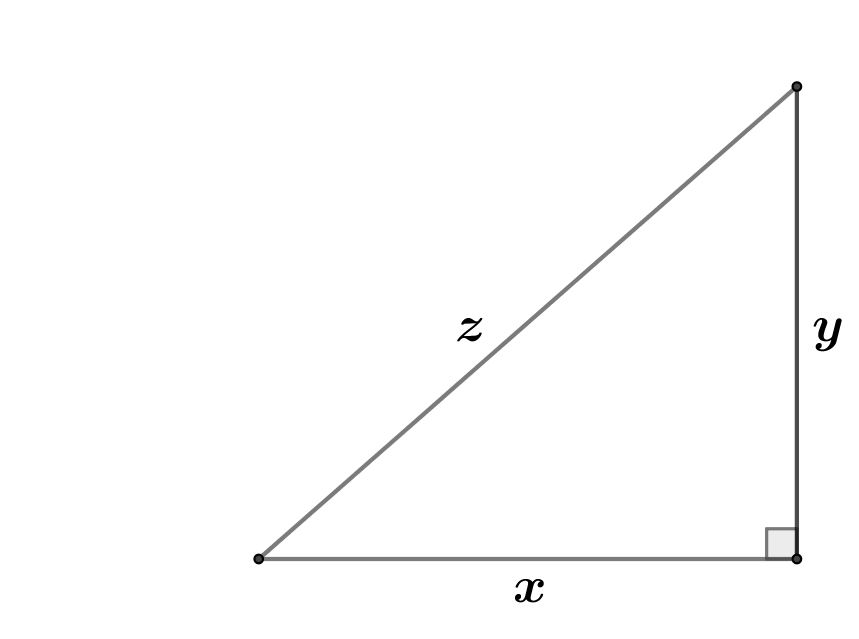}}%
\vspace{-3mm}

\begin{equation}
z=z(x,y)
\label{PythagoreFunction_0}
\end{equation}

\vspace{-1mm}

\noindent a relation that must be true for any $x$, $y$ and $z$ in a right triangle. In particular, after a slight increase in the length of $x$, leading to a small deviation $\delta x$, while $\delta y=0$, he would have found, at first order, that
\begin{equation}
\delta z=\partial_x z \,\delta x 
\label{PartialPyth1}
\end{equation}

\vspace{-4mm}

\noindent\hspace{-.5mm}\parbox{87mm}{
\includegraphics[width=43.5mm,keepaspectratio]{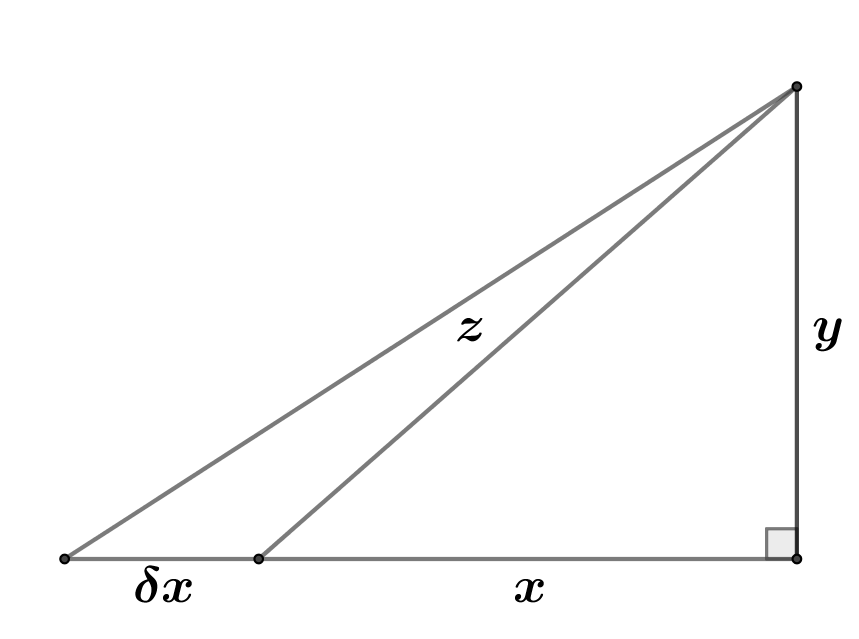}%
\includegraphics[width=43.5mm,keepaspectratio]{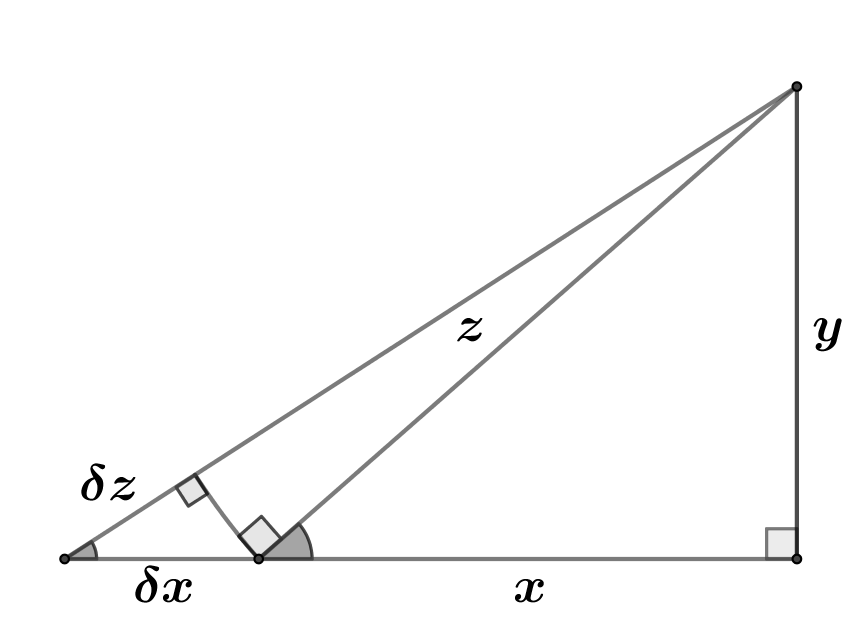}
}

\noindent  Here, $\delta x$ is the length of the hypotenuse of a right triangle with one leg of length $\delta z$. At first order, this small triangle is similar to the initial one. Using the Thales theorem that he has just found out, 
\begin{equation}
\delta z=\frac{x}{z}\delta x 
\end{equation}

\noindent Substituting this result to $\delta z$ in eq. (\ref{PartialPyth1}) would have led him to the partial differential equation
\begin{equation}
z\partial_x z=x 
\label{PartialPyth2}
\end{equation}

\noindent whose general solution is
\begin{equation}
z^2(x,y)=x^2+k(y) 
\end{equation}

\noindent where $k(y)$ is an arbitrary function of $y$. But the function $z$ has to be symmetric in $x$ and $y$ (i.e. he could have made the same reasoning with a non-zero $\delta y$ while $\delta x=0$) and $z(x,0)=x$. Hence
\begin{equation}
z^2(x,y)=x^2+y^2 
\label{PythagoreanTheorem}
\end{equation}

\noindent that is, the Pythagorean theorem, which, like Thales', was probably discovered long before -- and published a few centuries later by Euclid \cite[Book I, Prop. XLVII]{Euclid:300BC}. As will be shown in the last section, the same result can be obtained by considering a small rotation of the hypotenuse around one of its extremities followed by an infintesimal scale transformation. This proof is known and is published in a slightly different form in \cite{Staring:1996, Berndt:1988}.

\section{Apollonius of Perga}
\label{Apollonius}

\noindent\parbox{43.5mm}{\vspace{-6mm}
\includegraphics[width=43.5mm,keepaspectratio]{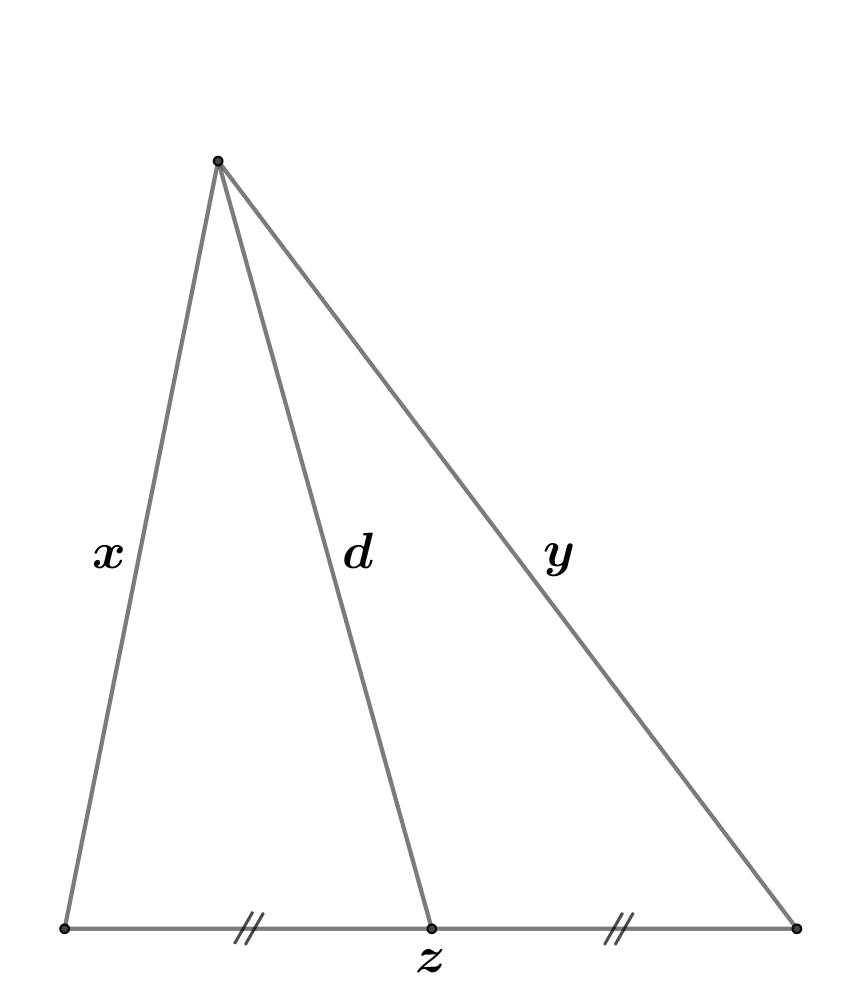}
\vspace*{-4mm}

\includegraphics[width=43.5mm,keepaspectratio]{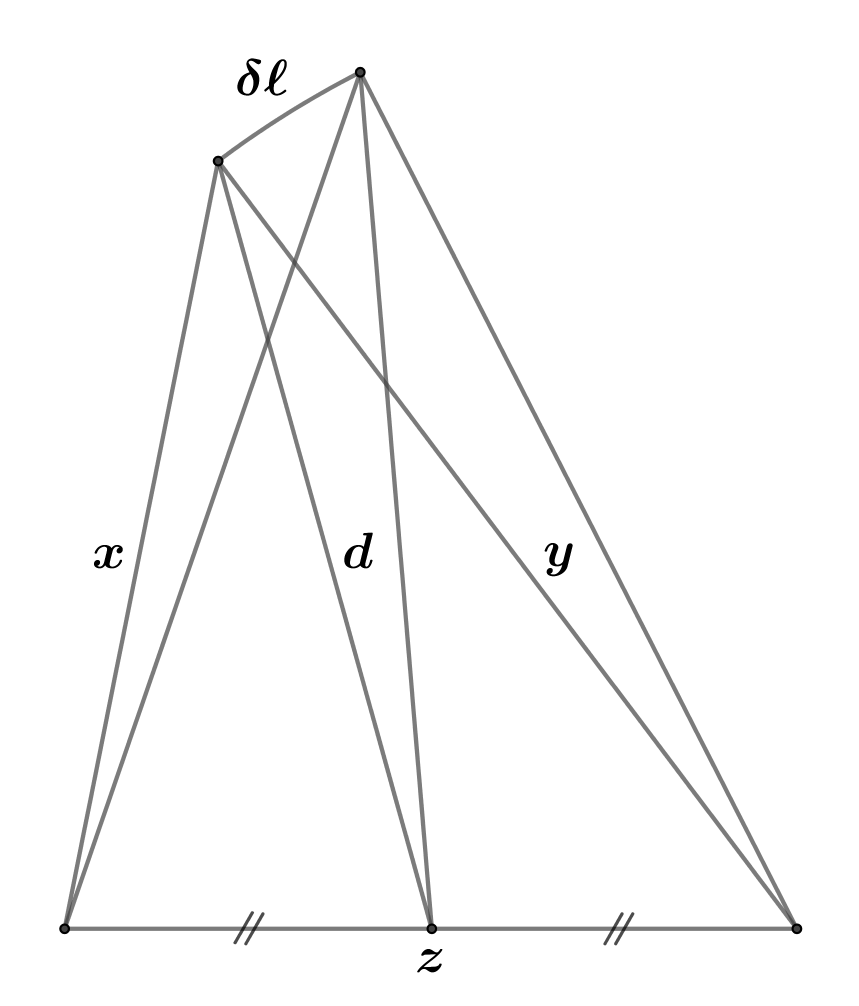}
}%
\parbox{43mm}{\vspace{3mm}

\noindent Thales and Pythagoras theorems are not the only ones that are named before famous Greek geometers. Newton could have gone a step further -- eastwards, a few centuries later -- and assumed that in any triangle of sidelengths $x$, $y$ and $z$, the length $d$ of the median relative to the $z$-length side is a smooth function of $x$, $y$ and $z$, i.e.
\begin{equation}
d=d(x,y,z)
\end{equation}

\noindent After an infinitesimal rotation of the $y$-length side around the $(y,z)$ vertex resulting in a small deviation $\delta x$, with $\delta y=0$ and $\delta z=0$, at first order:  
\begin{equation}
\delta d=\partial_x d \,\delta x 
\label{PartialAppollonius1}
\end{equation}
}

\noindent In order to get an expression for $\delta x$ and $\delta d$ and then a differential equation leading to the would-be theorem, consider the infinitesimal arc travelled by the moved vertex, of length $\delta\ell$.

\noindent\parbox{43.5mm}{\vspace{-3mm}

\includegraphics[width=43.5mm,keepaspectratio]{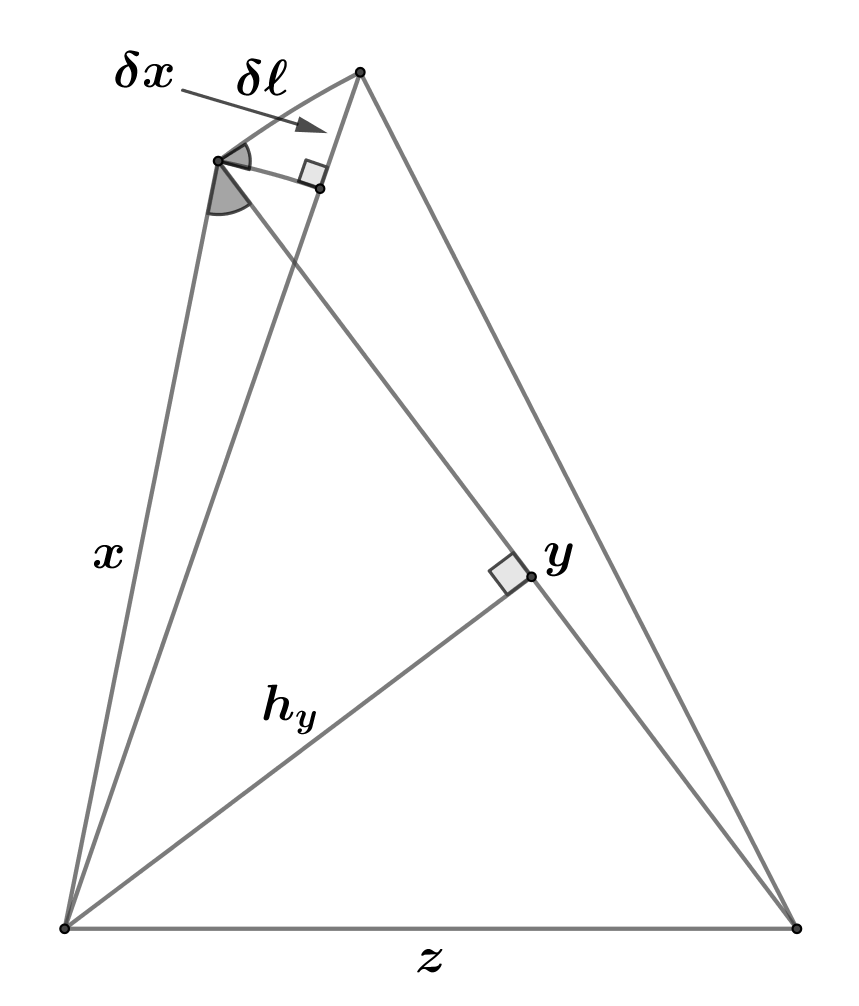}%
}%
\parbox{43mm}{\vspace{2mm}
\noindent At first order, it can be seen as the hypotenuse of a small right triangle with one leg of length $\delta x$, which is similar to a larger right triangle whose corresponding leg is the $h_y$-length height relative to the $y$-length side, and the hypotenuse is the $x$-length side, so that
\begin{equation}
\delta x=\frac{h_y}{x}\delta \ell
\label{DeviationAppollonius1}
\end{equation}
}

\noindent\parbox{43.5mm}{%
\includegraphics[width=43.5mm,keepaspectratio]{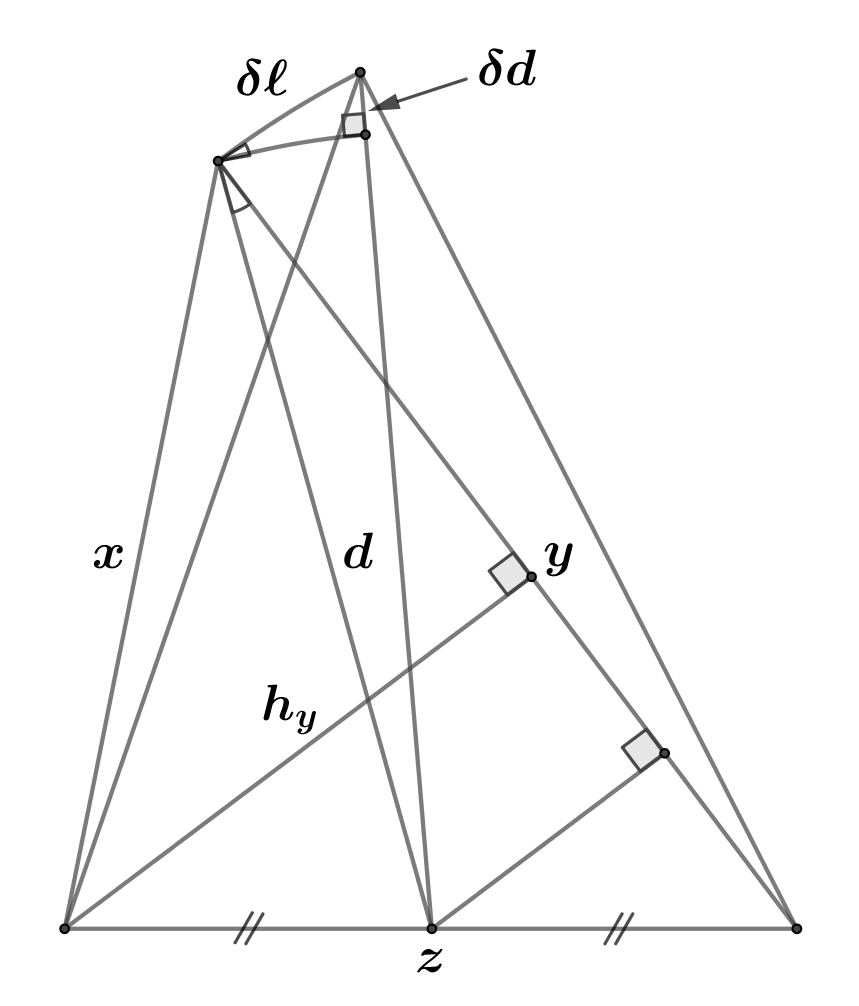}
}%
\parbox{43mm}{%
\noindent The $\delta\ell$-length arc is also the first order hypotenuse of another small right triangle with one leg of length $\delta d$, which is similar to the triangle with a $d$-length hypotenuse and whose corresponding leg is a segment starting from the foot of the median and parallel to -- and thus half of the length of -- the $h_y$-length height, so that
}
\begin{equation}
\delta d=\frac{h_y}{2d}\delta \ell
\label{DeviationAppollonius2}
\end{equation}

\noindent Inserting those deviations in eq. (\ref{PartialAppollonius1}) leads to the partial differential equation
\vspace{-1mm}
\begin{equation}
d\partial_x d=\frac{x}{2} 
\label{PartialAppollonius2}
\end{equation}

\noindent whose general solution is
\begin{equation}
d^2(x,y,z)=\frac{x^2}{2}+k(y,z) 
\end{equation}

\noindent where $k(y,z)$ is a function of $y$ and $z$. But $d(x,y,z)$ has to be symmetric in $x$ and $y$ (i.e. we can make the same reasoning with a non-zero $\delta y$ while $\delta x=0$). Hence
\begin{equation}
d^2(x,y,z)=\frac{x^2+y^2}{2}+c(z) 
\end{equation}

\noindent with $c(z)$ a function of $z$. Furthermore, if $x=0$ (or $y=0$), $y=z$ (or $x=z$) and $d=z/2$. This yields $c(z)=-z^2/4$, which can alternatively be found by invoking the Pythagorean theorem for $y=x$. The particular solution reads 
\begin{equation}
d^2(x,y,z)=\frac{x^2+y^2-2(z/2)^2}{2} 
\label{ApolloniusTheorem}
\end{equation}

\noindent that is, Apollonius's theorem, to be found in a slightly more elaborate form in 
\cite{Apollonius:200BC}.

\section{Matthew Stewart}
\label{Stewart}

\noindent Suppose Newton was born before Stewart, an 18th-century Scottish mathematician (and reverend). Well, he was. Perhaps he was not interested, or did not have the time, otherwise he could have used this tool to generalize Apollonius's theorem to any cevian. 
\noindent\parbox{86mm}{
\includegraphics[width=43mm,keepaspectratio]{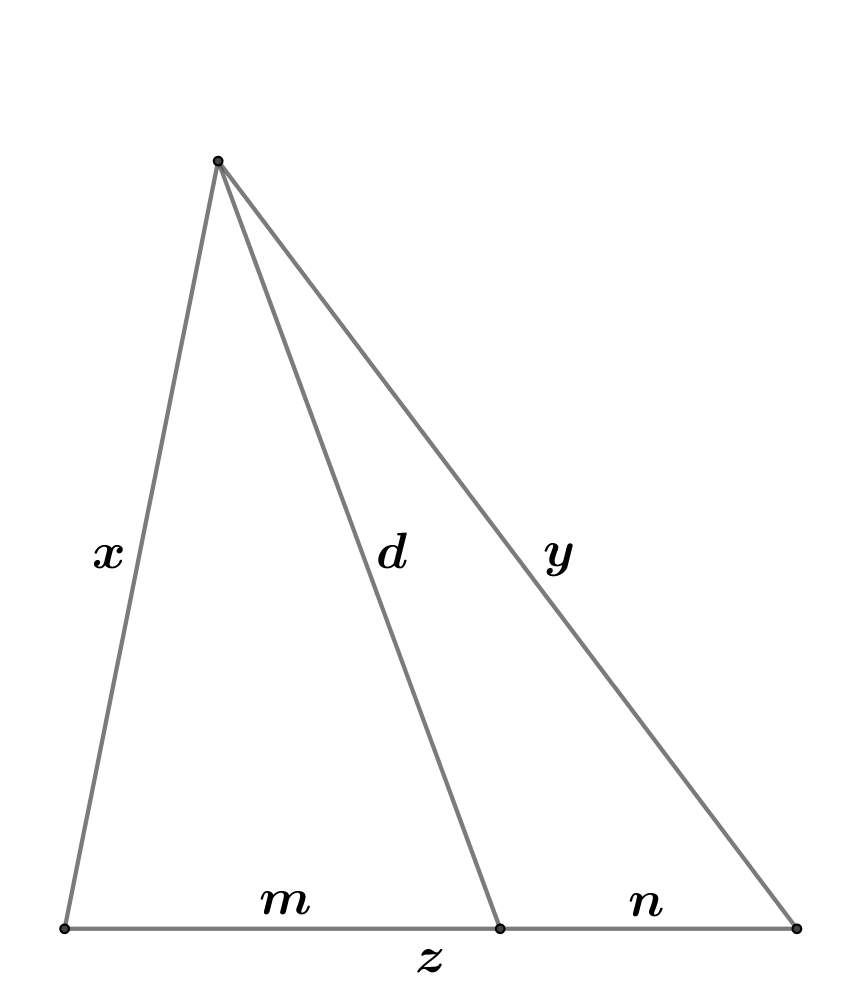}%
\includegraphics[width=43mm,keepaspectratio]{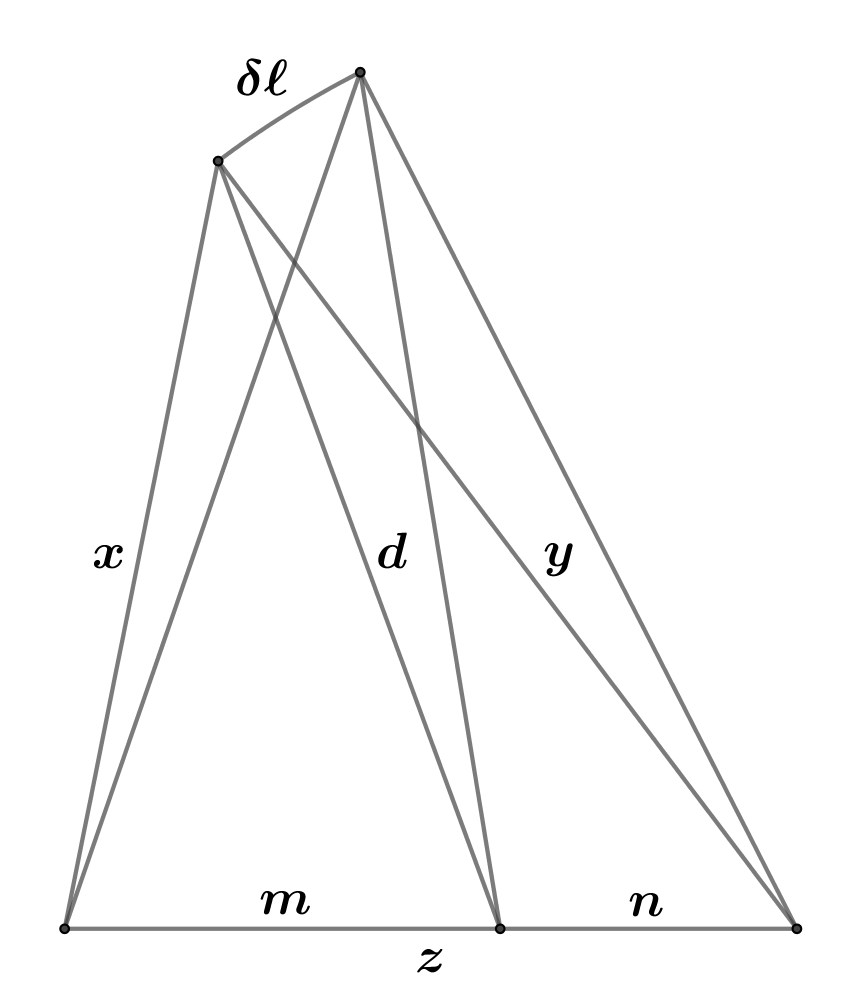}%
}

\noindent In a triangle of sidelengths $x$, $y$ and $z$, assume that the length $d$ of a cevian dividing the side of length $z$ in two segments of lengths $m$ and $n$, is a smooth function of $x$, $y$, $m$ and $n$, that is
\begin{equation}
d=d(x,y,m,n)
\end{equation}

\noindent After an infinitesimal rotation of the $y$-length side around the $(y,z)$ vertex resulting in a small deviation $\delta x$, with $\delta y=0$, $\delta m=0$ and $\delta n=0$, at first order: 
\begin{equation}
\delta d=\partial_x d \,\delta x 
\end{equation}

\noindent\parbox{86mm}{
\includegraphics[width=43mm,keepaspectratio]{GCApollonius_dx.png}%
\includegraphics[width=43mm,keepaspectratio]{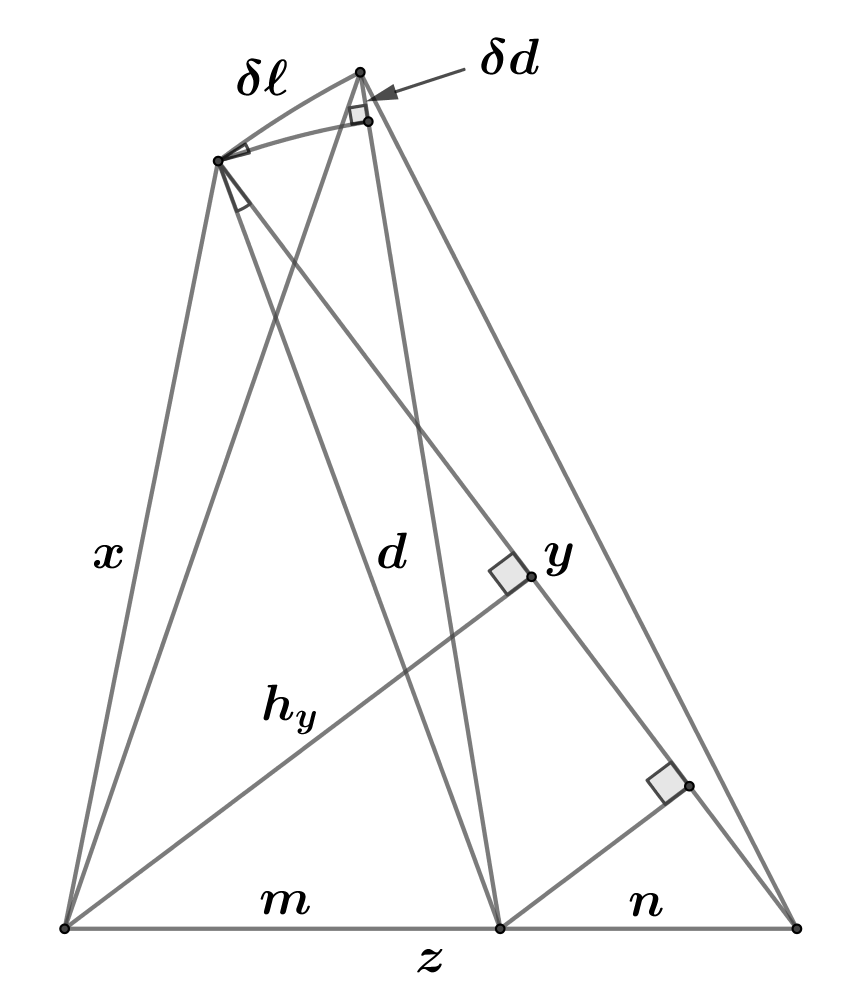}
}%

\noindent Using the same similarities as for the Apollonius's theorem, with the unique difference that the foot of the cevian is not necessarily the middle of the $(m+n)$-length side but falls at a distance $n$ from its right vertex, we find
\begin{equation}
\delta x=\frac{h_y}{x}\delta \ell
\qquad \delta d=\frac{n h_y}{(m+n)d}\delta \ell
\end{equation}

\noindent $h_y$ being the length of the height relative to $y$, we have the partial differential equation
\begin{equation}
d\partial_x d=\frac{n}{m+n} x
\end{equation}

\noindent whose general solution is
\begin{equation}
d^2(x,y,m,n)=\frac{n}{m+n}x^2+k(y,m,n) 
\end{equation}

\noindent where $k(y,m,n)$ is a function of $y$, $m$ and $n$. But $d(x,y,m,n)$ must be symmetric in $(x,m)$ and $(y,n)$ (i.e. we can make the same reasoning with a non-zero $\delta y$ while $\delta x=0$). Hence
\begin{equation}
d^2(x,y,m,n)=\frac{nx^2+my^2}{m+n}+c(m,n) 
\end{equation}

\noindent with $c(m,n)$ a symmetric function of $m$ and $n$. Furthermore, if $x=0$ (or $y=0$), $y=m+n$ (or $x=m+n$) and $d=m$ (or $d=n$). This yields $k(m,n)=-mn$. The particular solution reads 
\begin{equation}
d^2(x,y,m,n)=\frac{n(x^2-m^2)+m(y^2-n^2)}{m+n}
\label{StewartTheorem}
\end{equation}

\noindent that is, Stewart's theorem \cite{Stewart:1746}.

\section{Heron of Alexandria} 
\label{Heron} 

\noindent Intoxicated by his findings, Newton could have switched to a more elaborate, though older, challenge -- as probably did an Ancient Greek Roman Egyptian mathematician... What if, for any triangle, the area $A$ could be a smooth function of the sides lengths $x$, $y$ and $z$? He would have assumed
\begin{equation}
A=A(x,y,z)
\end{equation}

\noindent After an infinitesimal rotation of the $y$-length side around the $(y,z)$ vertex resulting in a small deviation $\delta x$, while $\delta y=0$ and $\delta z=0$, at first order: 
\begin{equation}
\delta A=\partial_x A \,\delta x 
\label{PartialHeron1}
\end{equation}

\noindent\parbox{43.5mm}{\vspace{-4mm}
\includegraphics[width=43.5mm,keepaspectratio]{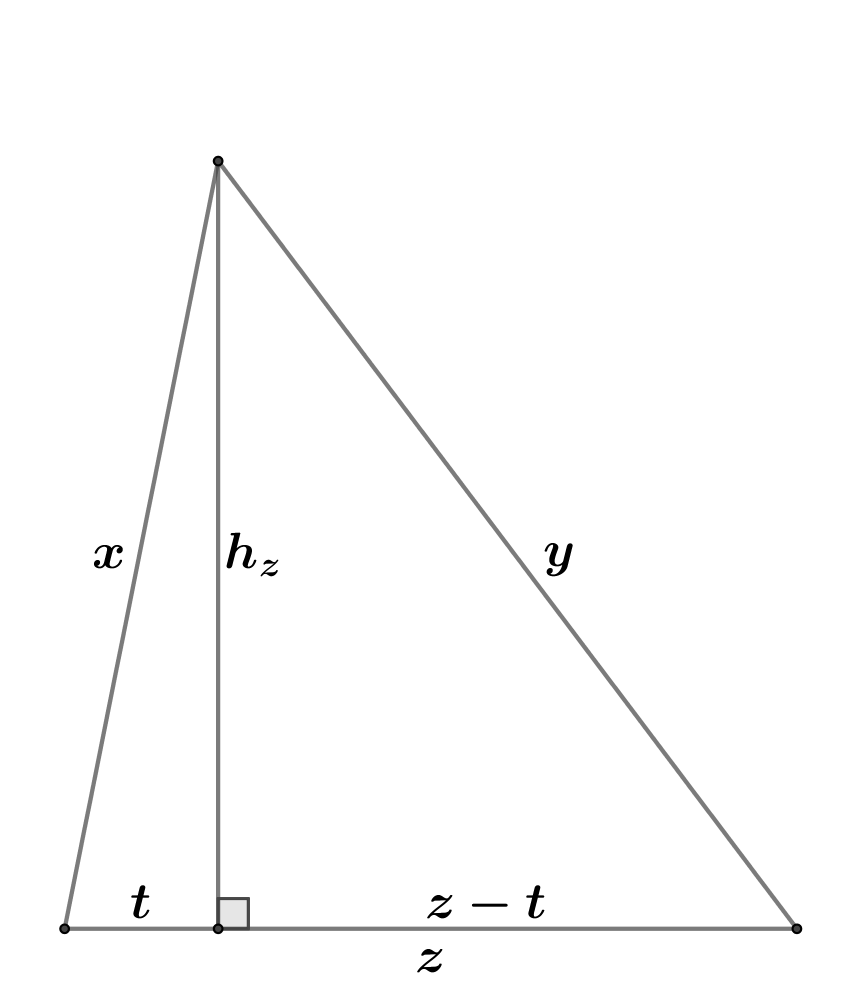}
}%
\parbox{43mm}{%
\noindent First note that the $h_z$-length height relative to the $z$-length side divides the initial triangle in two right triangles of horizontal legs of lengths $t$ and $z-t$ respectively.
One can express $h_z$ as a result of the Pythagorean theorem in both right triangles. }

\noindent Equating those expressions yields $x^2-t^2=y^2-(z-t)^2$ and hence 
\begin{equation}
t=\frac{x^2-y^2+z^2}{2z}\qquad z-t=\frac{y^2-x^2+z^2}{2z}
\label{HauteurDiviseBase}
\end{equation}

\noindent\parbox{86mm}{
\includegraphics[width=43mm,keepaspectratio]{GCApollonius_dx.png}%
\includegraphics[width=43mm,keepaspectratio]{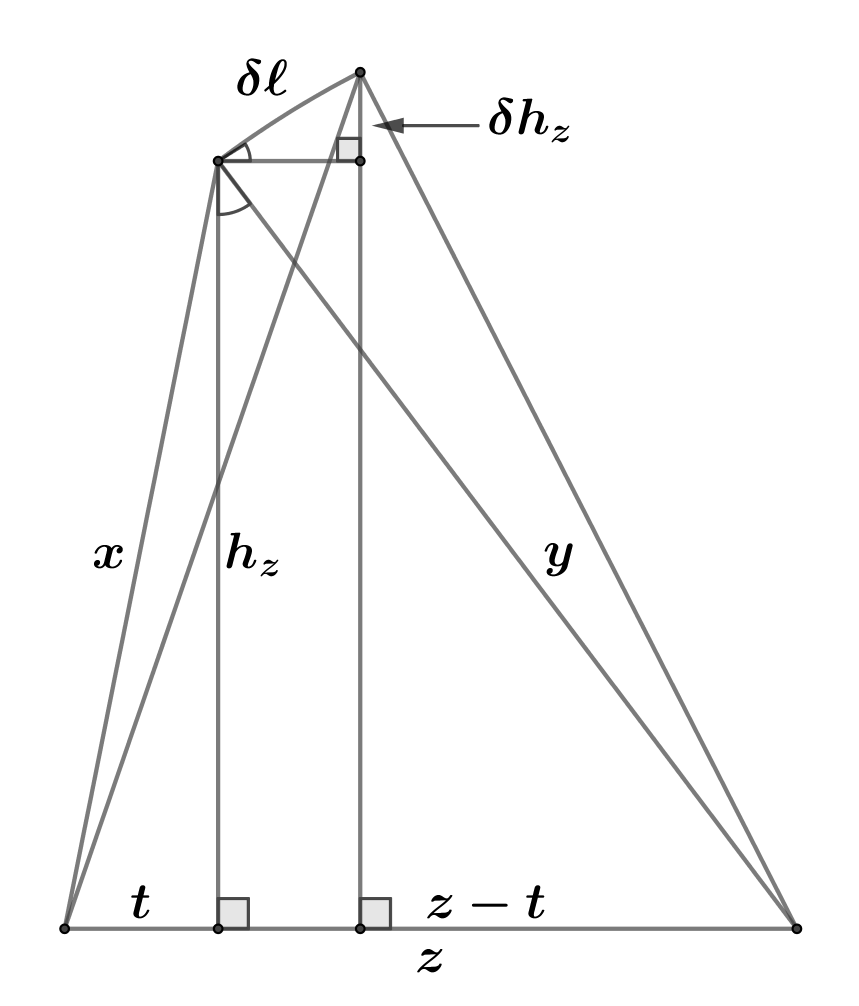}
}%

\noindent Again, $\delta\ell$ is the length of the infinitesimal arc travelled by the moved vertex. Like in the two last sections, it can be considered as the first-order hypotenuse of a small triangle whose similarity with a larger one allows to find $\delta x$. But it is also, at first order, the hypotenuse of another small triangle with one leg of length $\delta h_z$, which is similar to the large right triangle whose corresponding leg is the $(z-t)$-length segment, and the hypotenuse the $y$-length side. Since $h_y=2A/y$ and $\delta h_z=2\delta A/z$, we have
\begin{equation}
\delta x =\frac{2A}{xy}\delta \ell 
\qquad \delta A=\frac{y^2-x^2+z^2}{4y}\delta \ell
\label{DeviationsHeron1}
\end{equation}

\noindent Plugging in results (\ref{DeviationsHeron1}) into equation (\ref{PartialHeron1}), gives the partial differential equation
\begin{equation}
A\partial_x A=\frac{1}{8}[x(y^2+z^2)-x^3]
\label{PartialHeron2}
\end{equation}

\noindent which can be integrated out to give the general solution
\begin{equation}
A^2(x,y,z)=\frac{1}{16}[2x^2(y^2+z^2)-x^4+k(y,z)]
\end{equation}

\noindent where $k(y,z)$ is an homogeneous function of $y$ and $z$. Since $A(x,y,z)$ must be symmetric in $x$, $y$ and $z$ (i.e. we can make the same reasoning with a non-zero $\delta y$ or $\delta z$), $k(y,z)=2y^2z^2-y^4-z^4$. Hence
\begin{equation}
A(x,y,z)=\frac{1}{4}\sqrt{2(x^2y^2+x^2z^2+y^2z^2)-(x^4+y^4+z^4)}
\end{equation}

\noindent which can be factorized into the Heron theorem \cite{Heron:70AD}
\begin{equation}
A(x,y,z)=\sqrt{\frac{x+y+z}{2}\frac{-x+y+z}{2}\frac{x-y+z}{2}\frac{x+y-z}{2}}
\label{HeronTheorem}
\end{equation}
whose discovery could actually be Archimedes' \cite{Heath:1921}.

\section{Jamshid al-Kashi}
\label{alKashi}

\noindent Newton could have chosen to deal with angles -- besides calculus, he knew a bit about trigonometry. Let us send him to Persia, a few centuries before his birth, and wonder wether in any triangle of sidelengths $x$, $y$ and $z$, the angle $\gamma=\widehat{(x,y)}$ could be a smooth function of $x$, $y$ and $z$, that is
\begin{equation}
\gamma=\gamma(x,y,z)
\end{equation}

\noindent After an infinitesimal rotation of the $y$-length side around the $(y,x)$ vertex resulting in a small deviation $\delta z$, with $\delta x=0$ and $\delta y=0$, at first order:  
\begin{equation}
\delta \gamma=\partial_z \gamma \,\delta z 
\label{PartialalKashi1}
\end{equation}

\noindent\parbox{86mm}{
\includegraphics[width=43mm,keepaspectratio]{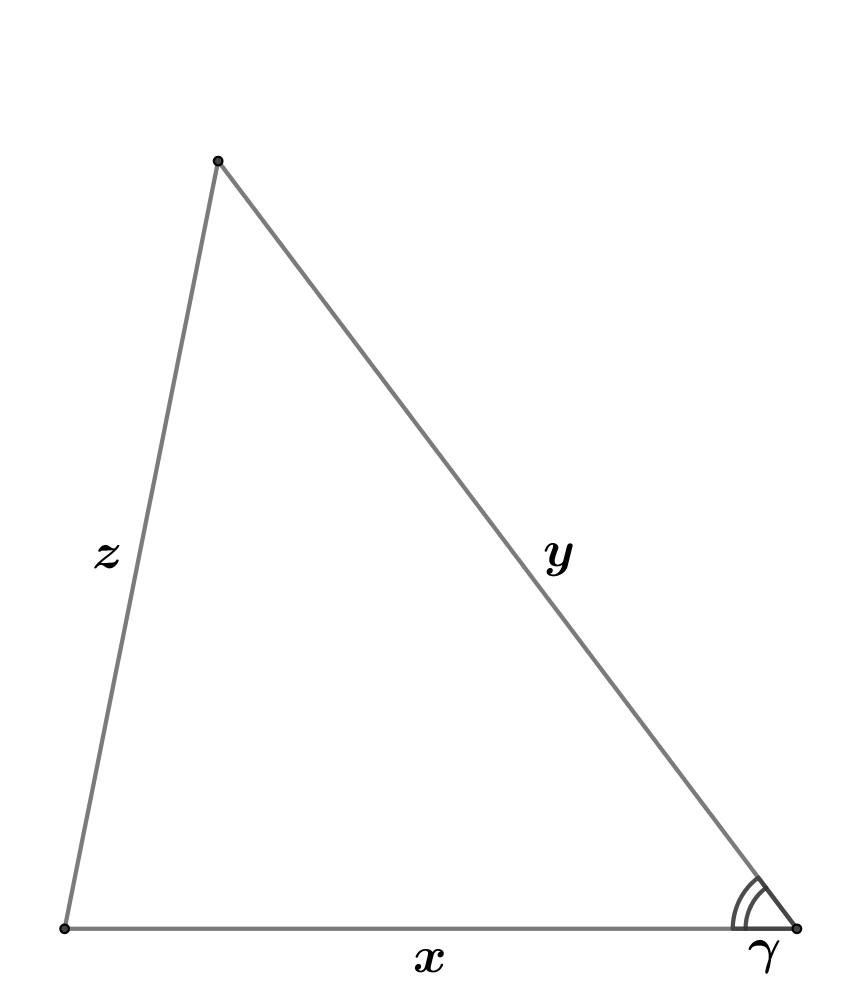}%
\includegraphics[width=43mm,keepaspectratio]{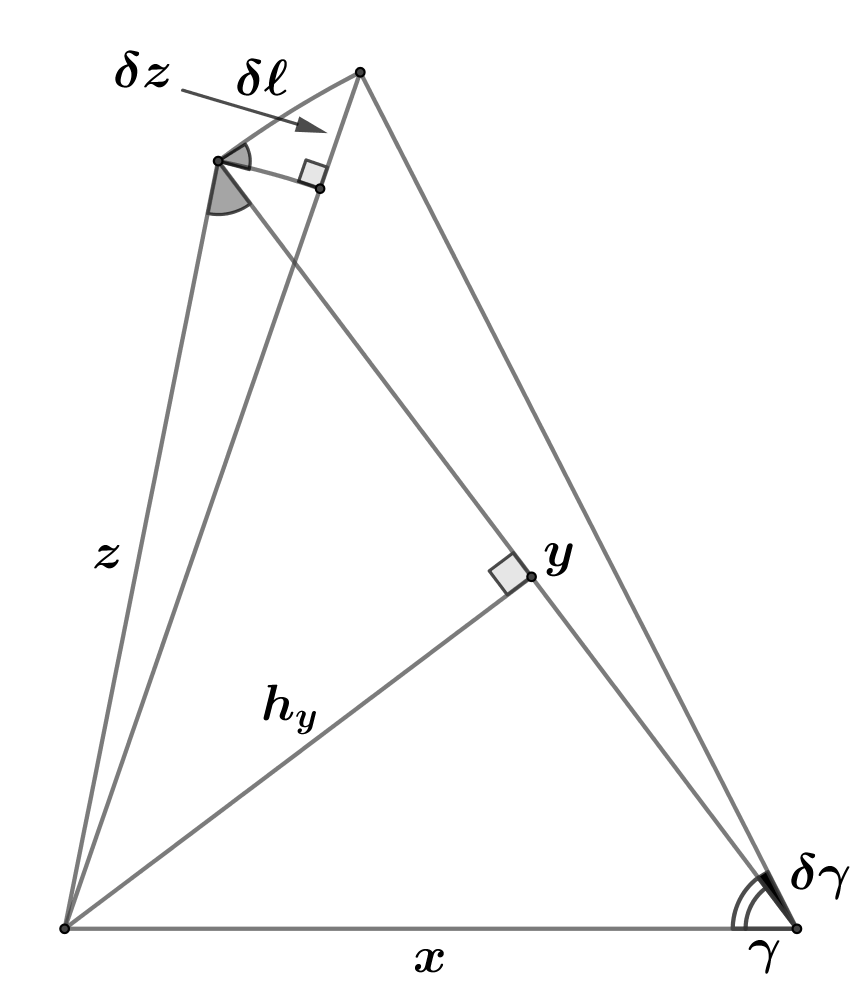}
}

\noindent While $\delta\gamma$ is easy to connect to $\delta\ell$, the length of the arc travelled by the moved vertex (in the illustrative figures, $x$ and $z$ have been swapped for aesthetic reasons), $\delta z$ can be determined thanks to the same similiarity as in the three previous sections. We have thus
\begin{equation}
\delta \gamma=\frac{\delta \ell}{y}
\qquad \delta z=\frac{h_y}{z}\delta \ell
\quad\mbox{with}\quad  h_y=x\sin\gamma
\label{DeviationsalKashi}
\end{equation}

\noindent Inserting those deviations in eq. (\ref{PartialalKashi1}) yields the partial differential equation
\begin{equation}
\sin\gamma\,\partial_z \gamma=\frac{z}{xy} 
\label{PartialalKashi2}
\end{equation}

\noindent whose general solution is
\begin{equation}
\cos[\gamma(x,y,z)]=-\frac{z^2+k(x,y) }{2xy}
\end{equation}

\noindent where $k(x,y)$ is a symmetric, homogeneous function of $x$ and $y$. According to Pythagoras, when $\gamma=\pi/2$, $z^2=x^2+y^2$, i.e. $k(x,y)=-x^2-y^2$. Hence
\begin{equation}
\cos[\gamma(x,y,z)]=\frac{-z^2+x^2+y^2 }{2xy}
\label{alKashiTheorem}
\end{equation}

\noindent that is, al-Kashi's theorem \cite{alKashi:1427} -- also known as the law of cosines or generalized Pythagorean theorem, and already familiar to Euclid \cite[Book II, Prop. XII \& XIII]{Euclid:300BC}.

\section{Olry Terquem}
\label{Terquem}

\noindent\parbox{43.5mm}{\vspace{-6mm}
\includegraphics[width=43.5mm,keepaspectratio]{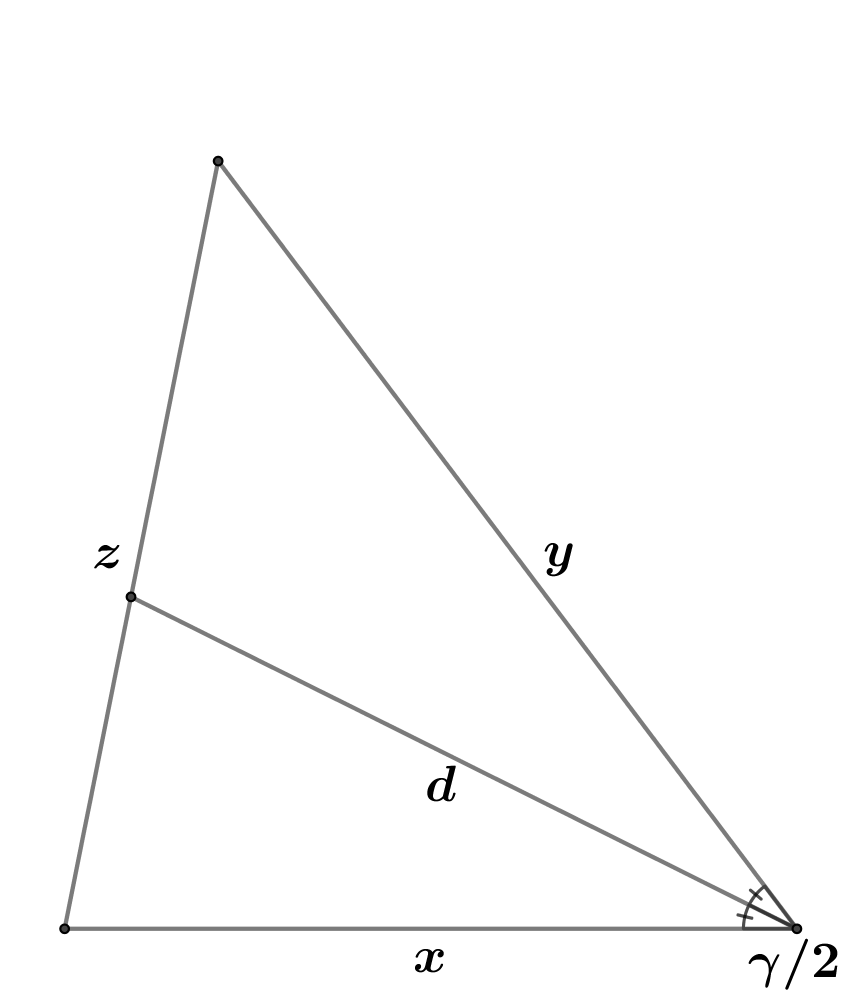}
}%
\parbox{43mm}{\vspace{3mm}
\noindent Completely exhilarated, Newton could have taken on a bigger piece and assumed that in any triangle of sidelengths $x$, $y$ and $z$, the length $d$ of the $\gamma=\widehat{(x,y)}$ angle bisector is a smooth function of $x$, $y$ and $z$, i.e.
\begin{equation}
d=d(x,y,z)
\end{equation}
}
\noindent After an infinitesimal rotation of the $y$-length side around the $(x,y)$ vertex resulting in a small deviation $\delta z$, with $\delta x=0$ and $\delta y=0$, at first order:  
\begin{equation}
\delta d=\partial_z d \,\delta z 
\end{equation}

\noindent Again, thanks to the same similiarity as in the four previous sections, $\delta z$ can easily be linked to $\delta\ell$, the length of the arc travelled by the moved vertex. 

\noindent\parbox{86mm}{
\includegraphics[width=86mm,keepaspectratio]{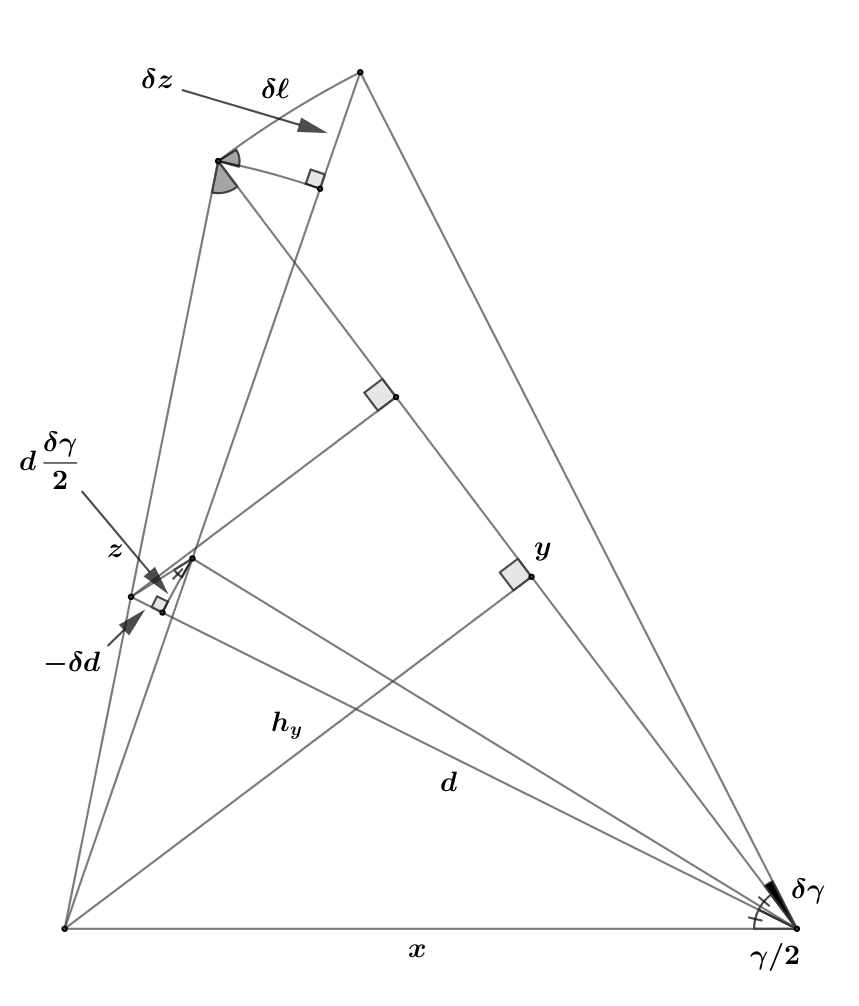}
}

\noindent It is a little more complicated for $\delta d$. First note that in the illustrative figure, $\delta d<0$, so that we will consider the positive length $-\delta d$. Then observe that when the $y$-length side infinitesimally rotates around the $(x,y)$ vertex, the foot of the $\gamma=\widehat{(x,y)}$ angle bisector moves along a perpendicular to the $y$-length side, just like the $(z,y)$ vertex. But the angle between this perpendicular and the angle bisector is the complementary of $\gamma/2$. Thus in the small right triangle of legs of lengths $-\delta d$ and $d\,\delta\gamma/2$, the opposite angle to the $-\delta d$-length leg is, at first order, equal to $\gamma/2$, implying that $\tan(\gamma/2)=-\delta d /(d\,\delta\gamma/2)$. Hence
\begin{equation}
\delta z=\frac{h_y}{z}\delta \ell
\qquad \delta d=-\tan{\frac{\gamma}{2}}\frac{d}{2}\delta \gamma
\quad\mbox{with}\quad  \delta\ell=y\delta\gamma
\end{equation}

\noindent $h_y$ being the length of the height relative to $y$. Using
\begin{equation}
\tan{\frac{\gamma}{2}}=\frac{\sin\gamma}{1+\cos\gamma}
\label{tanhalfangle}
\end{equation}

\noindent with
\begin{equation}
\sin\gamma=\frac{h_y}{x}\quad\mbox{and}\quad\cos\gamma=\frac{-z^2+x^2+y^2 }{2xy}
\end{equation}

\noindent we have the partial differential equation
\begin{equation}
\frac{\partial_z d}{d}=\frac{-z}{-z^2+(x+y)^2}
\end{equation}

\noindent whose general solution is
\begin{equation}
d(x,y,z)=k(x,y)\sqrt{(x+y)^2-z^2}
\end{equation}

\noindent where $k(x,y)$ is a symmetric function of $x$ and $y$. To determine it, note that in the particular case of a right triangle with hypotenuse of length $z$, the angle bisector is the diagonal of the inscribed square of sidelength $xy/(x+y)$ -- as can be deduced from similarities between the right triangles generated by the square in the initial triangle. We find $k(x,y)=\sqrt{xy}/(x+y)$. Hence
\begin{equation}
d(x,y,z)=\sqrt{xy\left(1-\frac{z^2}{(x+y)^2}\right)}
\label{TerquemTheorem}
\end{equation}

\noindent that is, the length of the angle bisector, as Terquem computed in the 19th century \cite{Terquem:1842}.

\section{Jean-Paul de Gua de Malves} 
\label{deGua} 

\noindent Armed with this powerful theorem-finding tool, Newton could have moved on to even bolder challenges, like leaving the plane for the real space, and imagining, say, a generalization of the Pythagorean theorem in three dimensions! Let him consider a trirectangular tetrahedron, that is a tetrahedron with a right angle corner, like the corner of a cube: what if, for any of them, the area of the face opposite to the right angle was a function of the areas of the other faces?

\noindent\hspace{-2.5mm}\parbox{43.5mm}{%
\includegraphics[width=43.5mm,keepaspectratio]{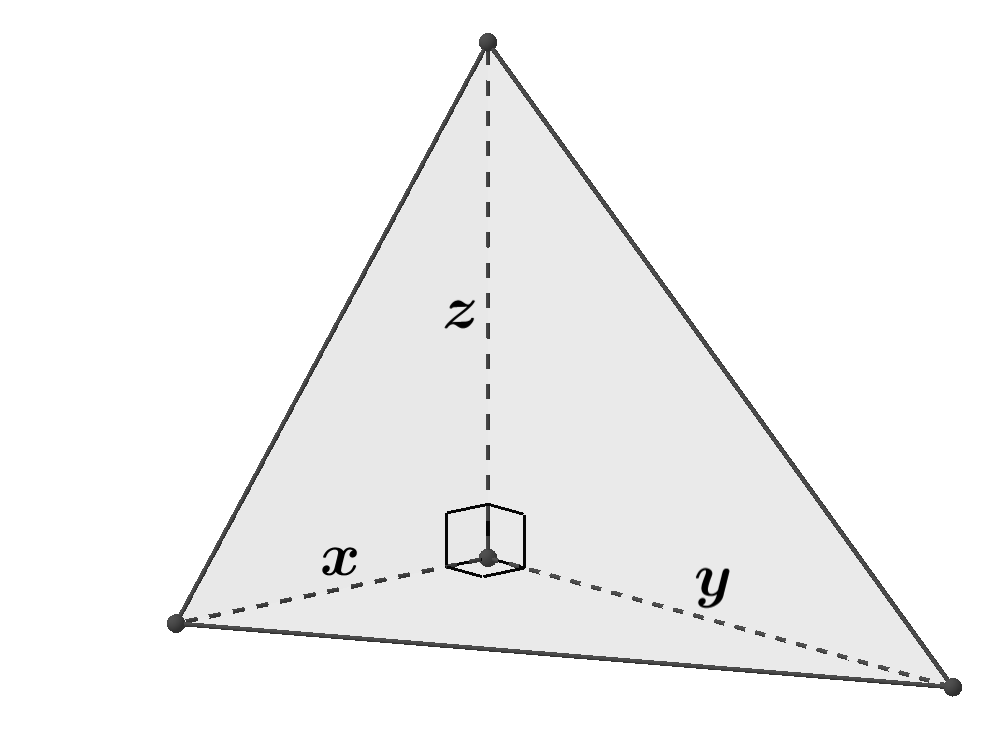}

\includegraphics[width=43.5mm,keepaspectratio]{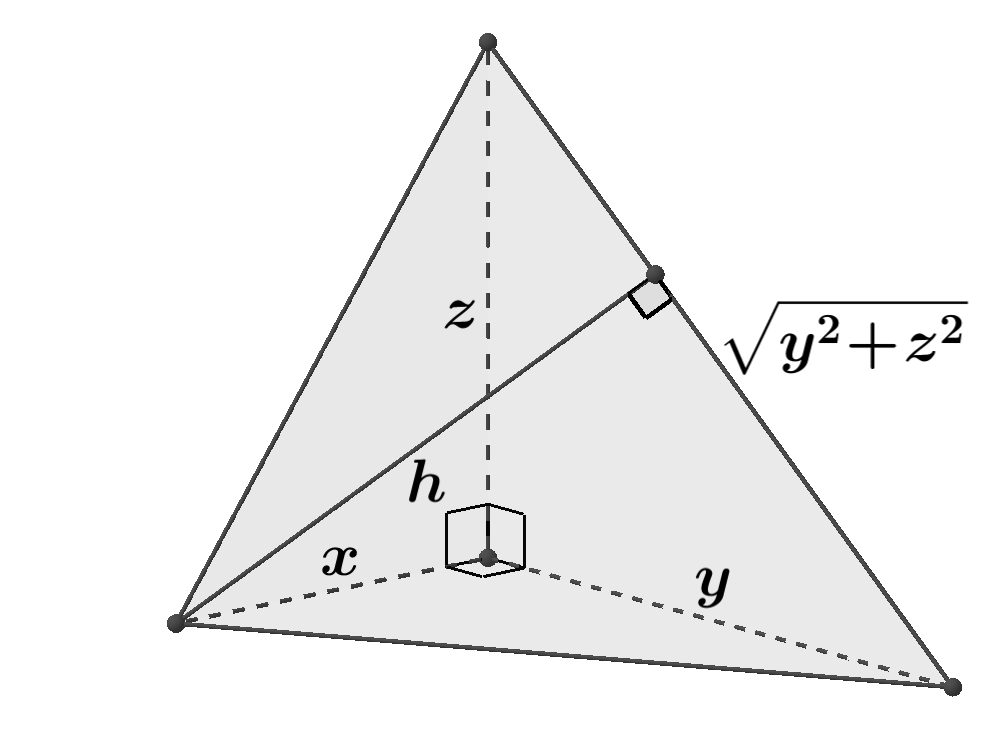}
}\hspace{1.9mm}%
\parbox{43mm}{%

A convenient way to parametrize the problem is to give arbitrary lengths to the three edges from the right angle vertex, say $x$, $y$ and $z$. The areas of the three right triangle faces are $xy/2$, $xz/2$ and $yz/2$. For the area of the last face, opposite to the right angle, say $A$, we can have an expression by choosing a base, say the edge of length $\sqrt{y^2+z^2}$ (thanks Pythagoras) and the relative height of length $h$. We have
\begin{equation}
A=\frac{1}{2}\sqrt{y^2+z^2}\,h
\label{ExpressionA}
\end{equation}
}

\noindent Let us go back to Newton and his obsession. He could have stated that $A$ is a smooth function of $x$ and $y$:
\begin{equation}
A=A(x,y,z)
\end{equation}

\noindent Choosing to slightly increase $x$, while leaving $y$ and $z$ invariants, that is, an infinitesimal deviation $\delta x$, with $\delta y=0$ and $\delta z=0$, we find
\begin{equation}
\delta A=\partial_x A \,\delta x 
\label{PartialGua}
\end{equation}

\noindent Eq. (\ref{ExpressionA}) implies that
\begin{equation}
\delta A=\frac{1}{2}\sqrt{y^2+z^2}\,\delta h
\label{DeviationAGua}
\end{equation}

\noindent\parbox{86mm}{
\includegraphics[width=43mm,keepaspectratio]{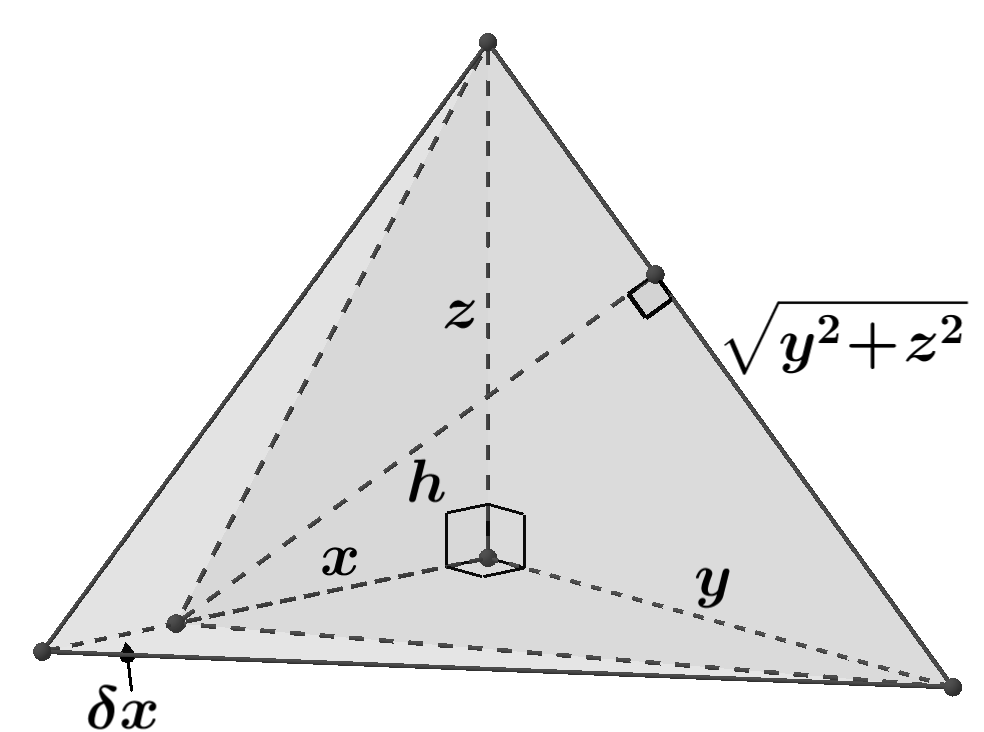}%
\includegraphics[width=43mm,keepaspectratio]{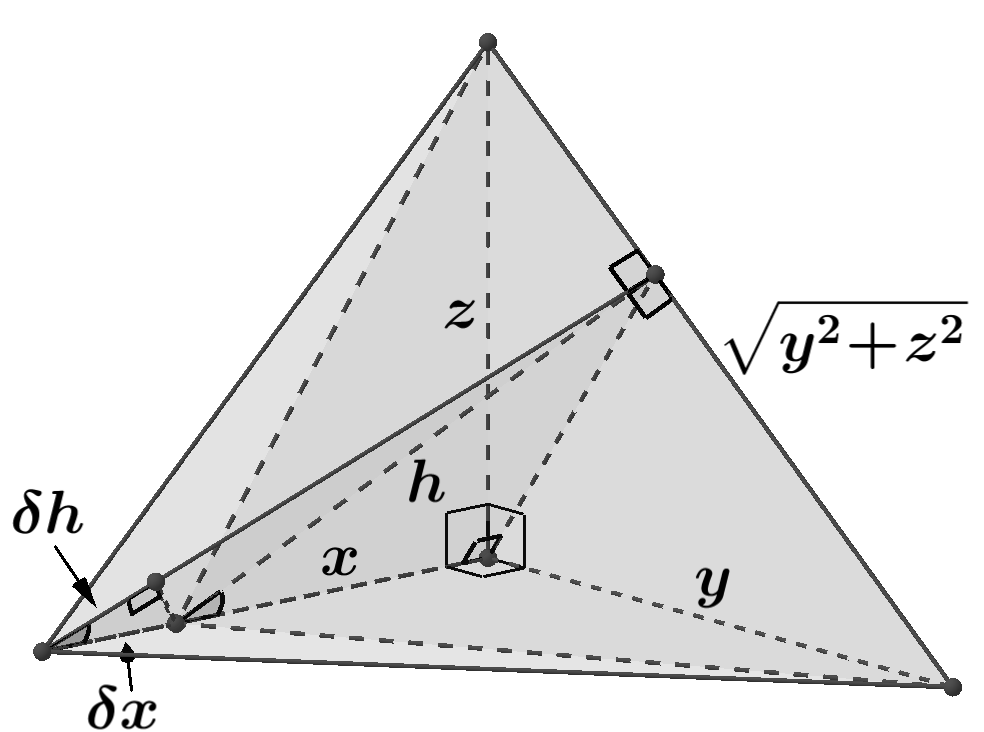}
}%

\vspace{4mm}

\noindent But what do we know of $\delta h$? First note that the foot of the $h$-length height is not affected by the deviation $\delta x$ since this $h$-length height and the $x$-length edge are in a plane orthogonal to the base of the $A$-area face. In this plane, we can check that at first order, the right triangle with $h$-length hypotenuse and $x$-length leg is similar to the one with $\delta x$-length hypotenuse and $\delta h$-length leg, so that
\begin{equation}
\delta h=\frac{x}{h}\delta x
\label{DeviationxGua}
\end{equation}

\noindent Combining this equation with result (\ref{DeviationAGua}), itself plugged in into eq. (\ref{PartialGua}) with $\delta y=0$, we have 
\begin{equation}
\frac{1}{2}\sqrt{y^2+z^2}\,\frac{x}{h}\delta x=\partial_x A \,\delta x
\end{equation}

\noindent Simplifying by $\delta x$ and using eq. (\ref{ExpressionA}) to get rid of $h$, we find a partial differential equation
\begin{equation}
A\partial_x A=\frac{1}{4}(y^2+z^2)x
\end{equation}

\noindent It can be integrated out to give the general solution
\begin{equation}
A^2(x,y,z)=\frac{1}{4}[(y^2+z^2)x^2+k(y,z)]
\end{equation}

\noindent with $k(y,z)$ an homogeneous and symmetric function of $y$ and $z$. Since $A(x,y,z)$ must itself be symmetric in $x$, $y$ and $z$ (i.e. we can make the same reasoning with a non-zero $\delta y$ or $\delta z$), $k(y,z)=y^2z^2$. Hence
\begin{equation}
A^2(x,y,z)=\left(\frac{xy}{2}\right)^2+\left(\frac{xz}{2}\right)^2+\left(\frac{yz}{2}\right)^2
\label{deGuaTheorem}
\end{equation}

\noindent known as de Gua's theorem \cite{deGua:1786}, first formulated by Descartes \cite{Descartes:1619}, which states that in any trirectangular tetrahedron, the square of the area of the face opposite to the right corner is equal to the sum of the squares of the areas of the other faces -- a three-dimensional generalization of the Pythagorean theorem.

\section{The Inscribed Angle}
\label{Angle}

\noindent Let us move on to the circle, and confront Newton to a simple problem. 

\noindent\parbox{86mm}{\parbox{43mm}{

\vspace{2mm}

\includegraphics[width=43mm,keepaspectratio]{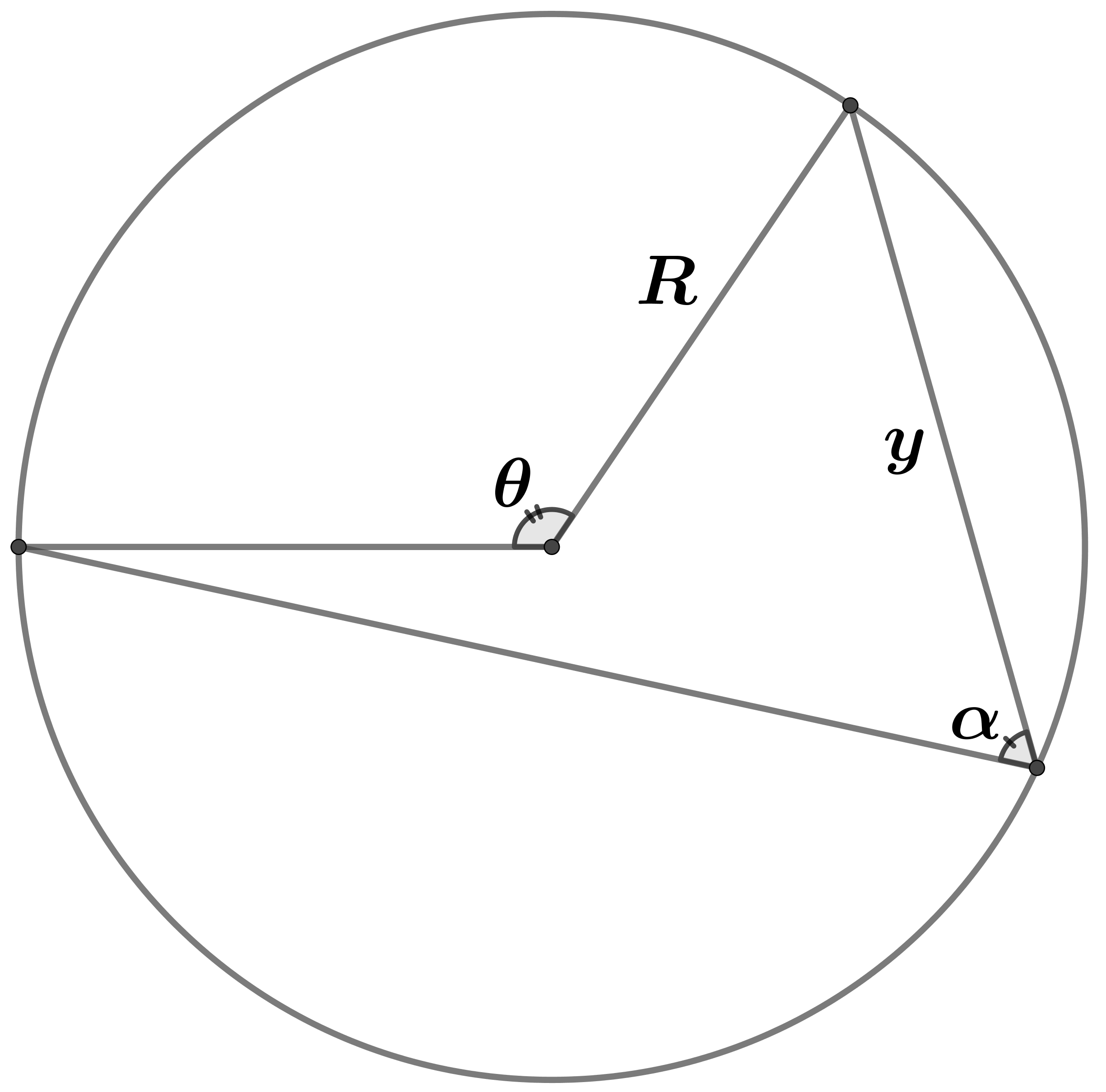}

\includegraphics[width=43mm,keepaspectratio]{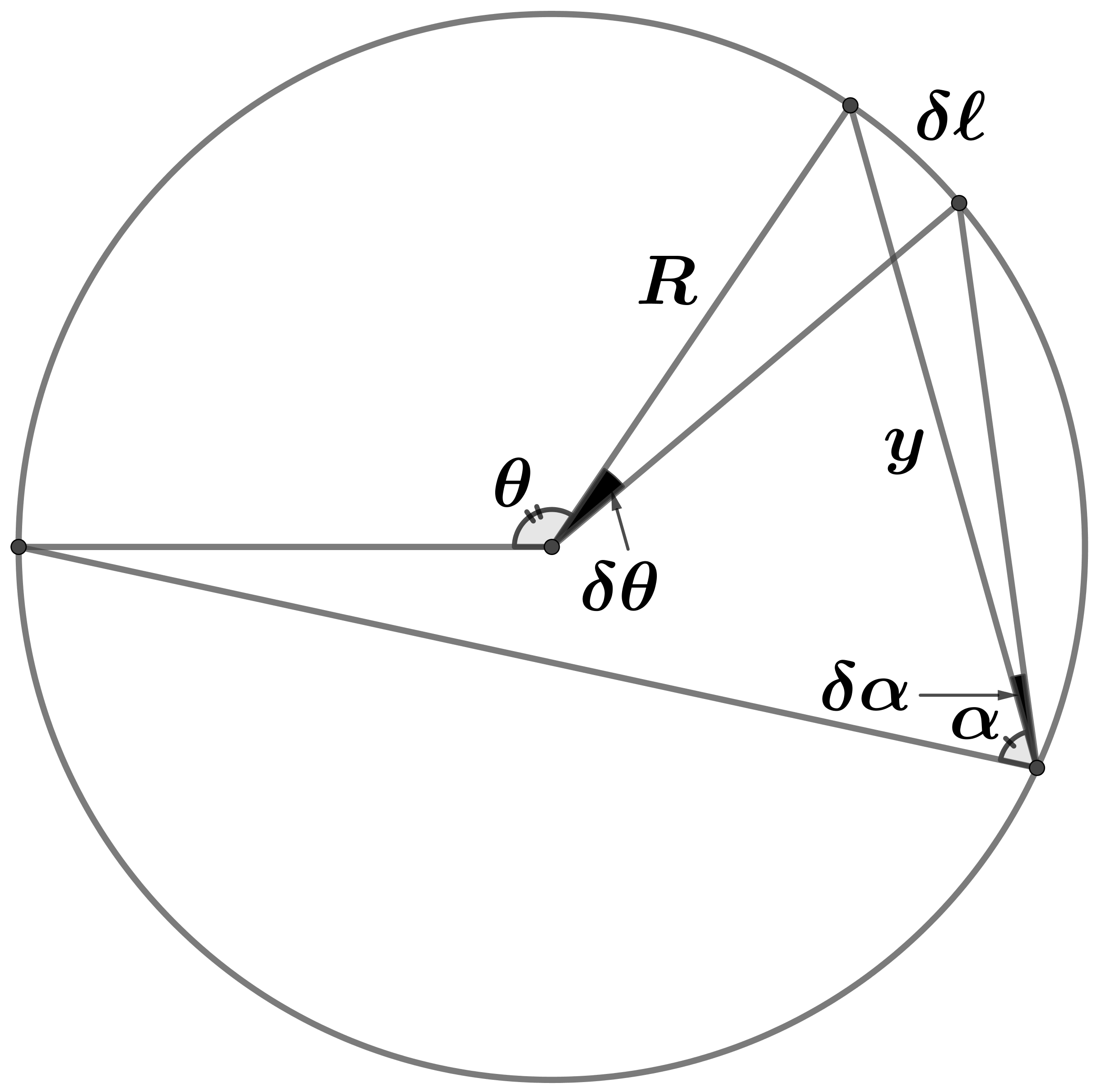}

}\hspace{2mm}%
\parbox{41mm}{

Assume that an arbitrary inscribed angle $\alpha$ is a smooth function of the central angle $\theta$ that intercepts the same arc on the circle, i.e. \vspace{-2mm}
\begin{equation}
\alpha=\alpha(\theta)
\end{equation}

\noindent After an infinitesimal deviation $\delta \theta$, at first order, 
\begin{equation}
\delta\alpha=\alpha'\delta\theta
\end{equation}

\noindent  At first order again, the small right triangle with one leg of length $y\delta\alpha$ and hypotenuse of length $\delta\ell$ is similar to the triangle of corresponding sidelengths $y/2$ and $R$ respectively. 
}}

\noindent\parbox{86mm}{
\includegraphics[width=43mm,keepaspectratio]{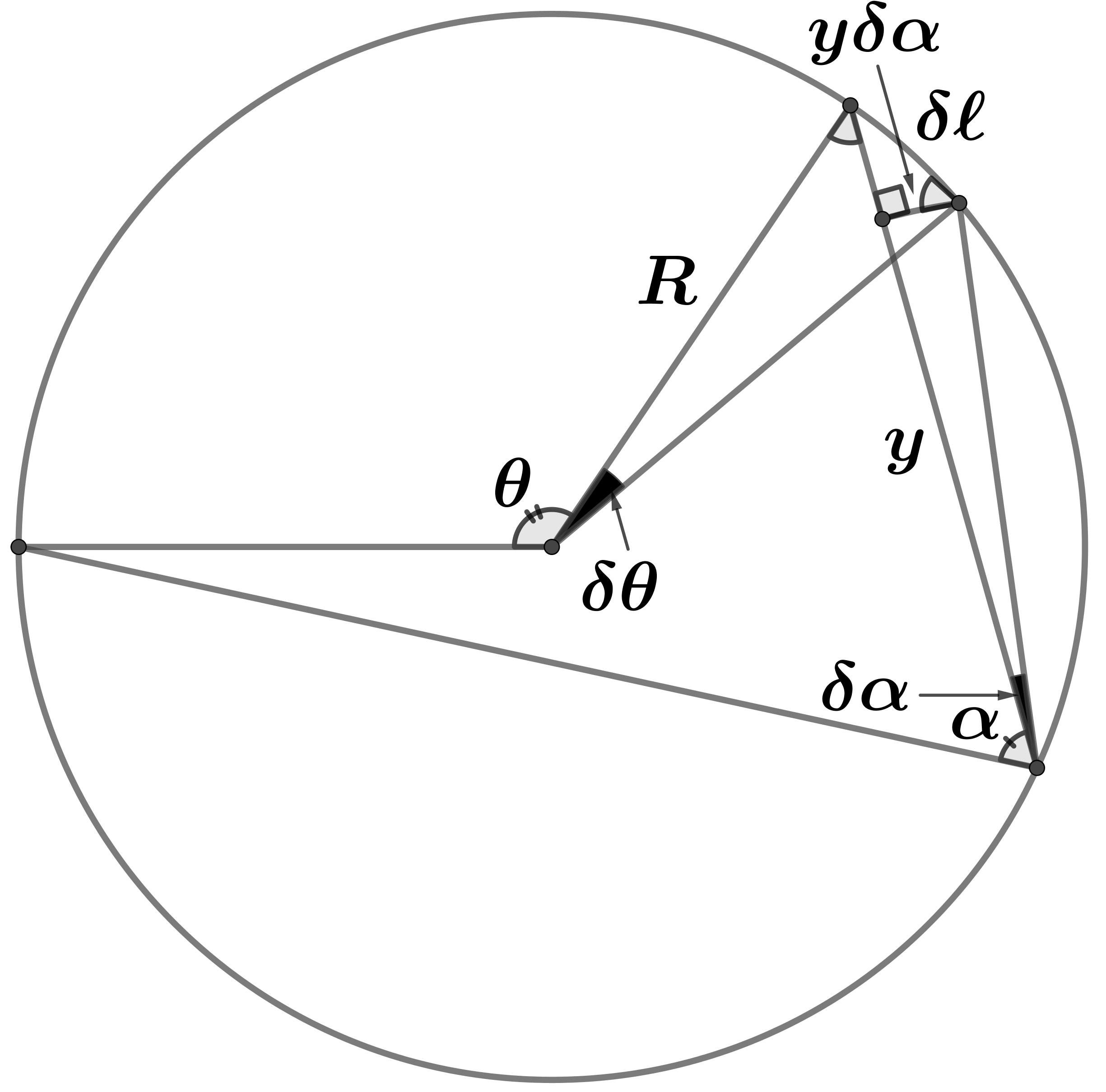}%
\includegraphics[width=43mm,keepaspectratio]{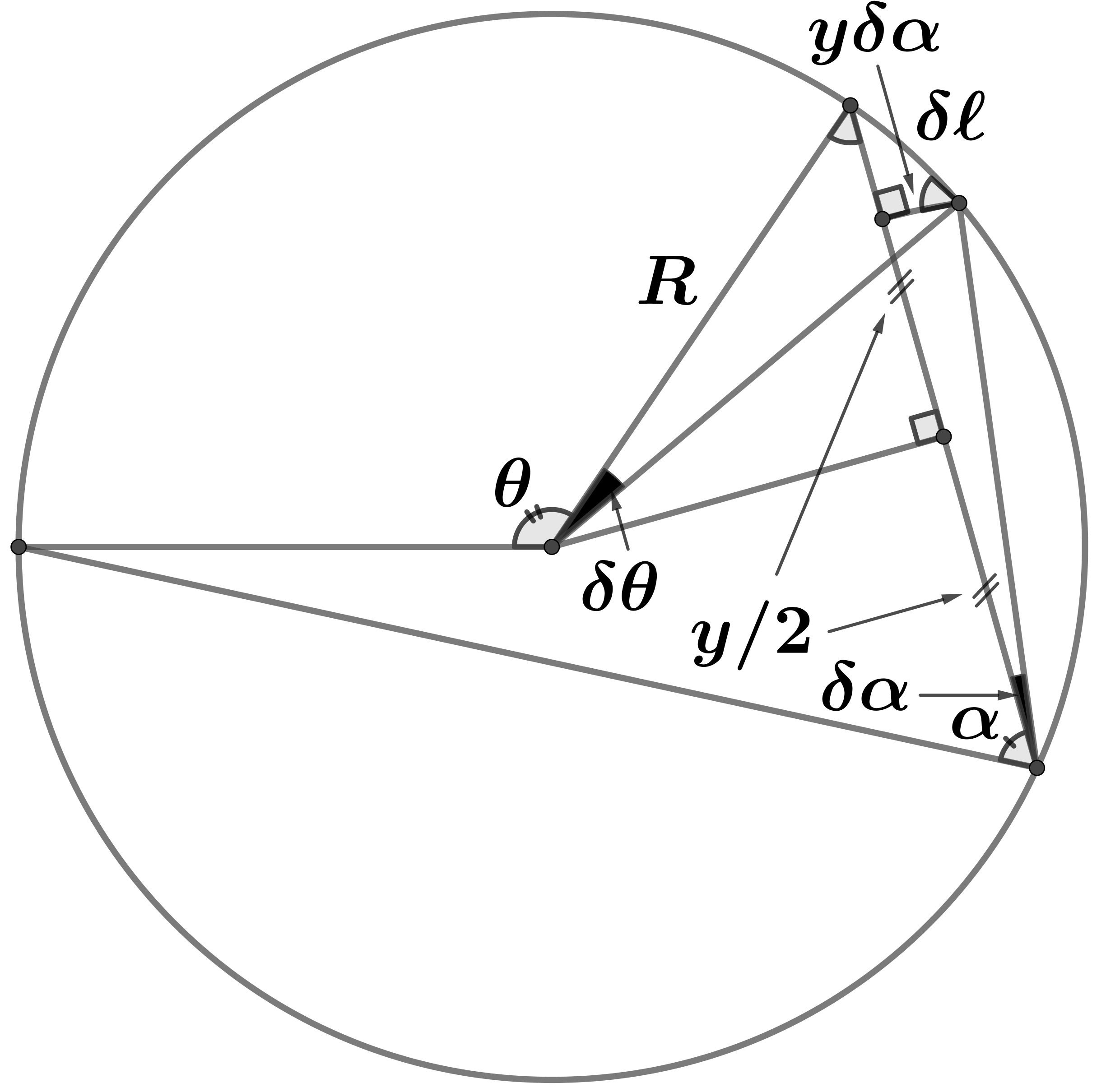}%
}%
\vspace{2.5mm}

\noindent Thus we have
\begin{equation}
y\delta\alpha=\frac{y/2}{R}\delta\ell
\end{equation}

\noindent Since $\delta \theta=\delta \ell/R$, this leads to the differential equation
\begin{equation}
\alpha'=\frac{1}{2}
\end{equation}

\noindent whose general solution is
\begin{equation}
\alpha(\theta)=\frac{\theta}{2}+k
\end{equation}

\noindent where $k$ is a constant that vanishes since $\alpha(0)=0$. 
Hence
\begin{equation}
\alpha(\theta)=\frac{\theta}{2}
\label{AngleTheorem}
\end{equation}

\noindent that is, the inscribed angle theorem.

\section{The Circumradius}
\label{Circumradius}

\noindent We could keep Newton in the circle and think to another question: of course, for any triangle, the circumradius length $R$ should be determined by the sidelengths $x$, $y$ and $z$. It was known to Euclid \cite[Book IV, Prop. V]{Euclid:300BC} and it is actually simple to prove it with the help of Heron's theorem. But let him play the game of finding it from scratch, that is, by postulating that
\begin{equation}
R=R(x,y,z)
\end{equation}

\vspace{2.5mm}
\noindent\parbox{86mm}{
\includegraphics[width=43mm,keepaspectratio]{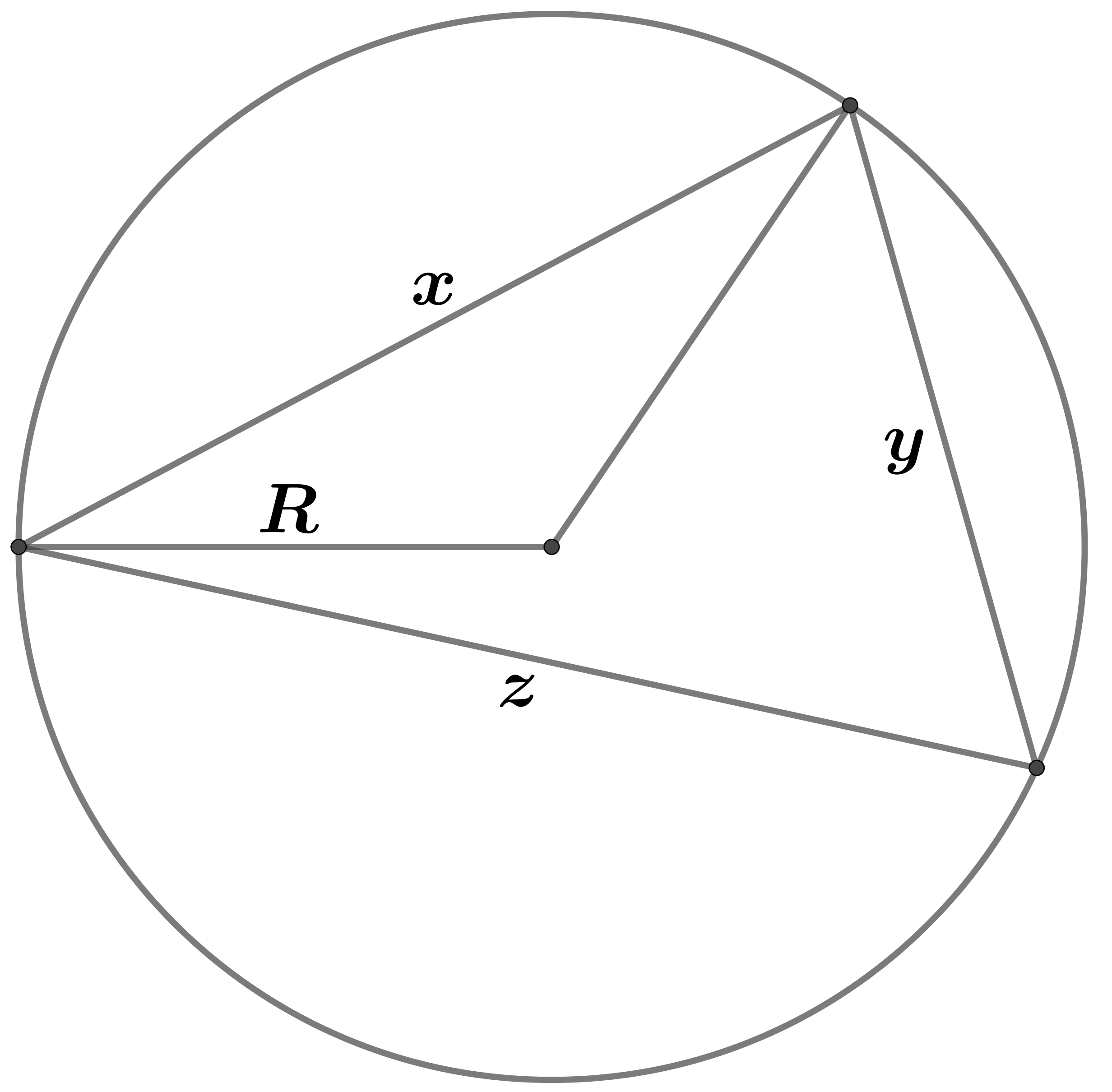}%
\includegraphics[width=43mm,keepaspectratio]{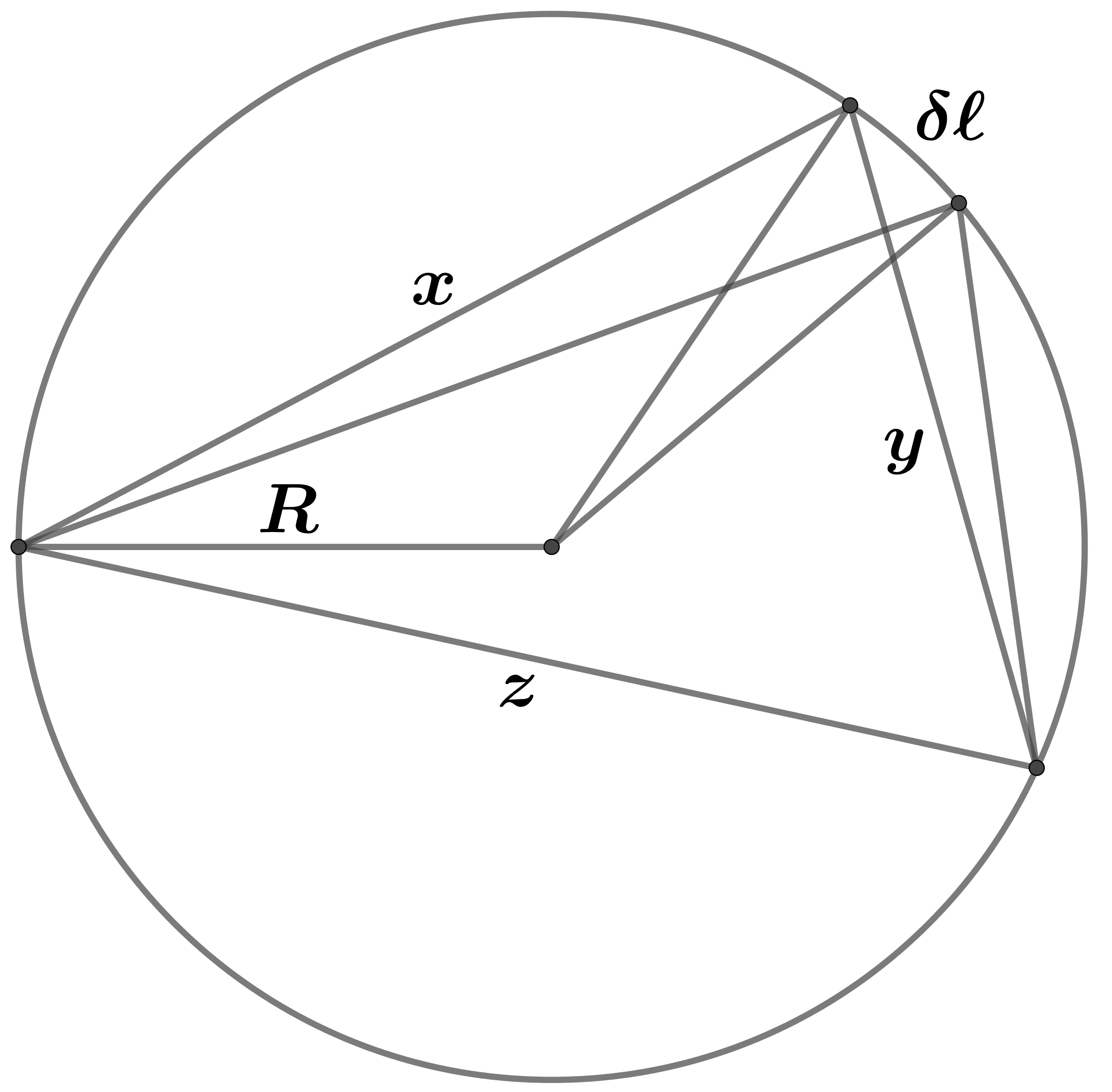}%
}%
\vspace{2.5mm}

\noindent If the $(x,y)$ vertex is slightly moved along the circumcircle, it generates infinitesimal deviations $\delta x$ and $\delta y$, while $\delta z=0$ and $\delta R=0$. At first order, we have: 
\begin{equation}
\partial_x R \,\delta x + \partial_y R \,\delta y=0
\label{PartialCircum}
\end{equation}

\noindent\parbox{86mm}{\parbox{41mm}{\noindent We need to go through some geometric considerations before proceeding: any angle in the triangle has the same magnitude as the ones between its opposite side and the tangent lines to the circumcircle from the two other vertices, since they intercept the same arc. }\hspace{2mm}%
\parbox{43mm}{\includegraphics[width=43mm,keepaspectratio]{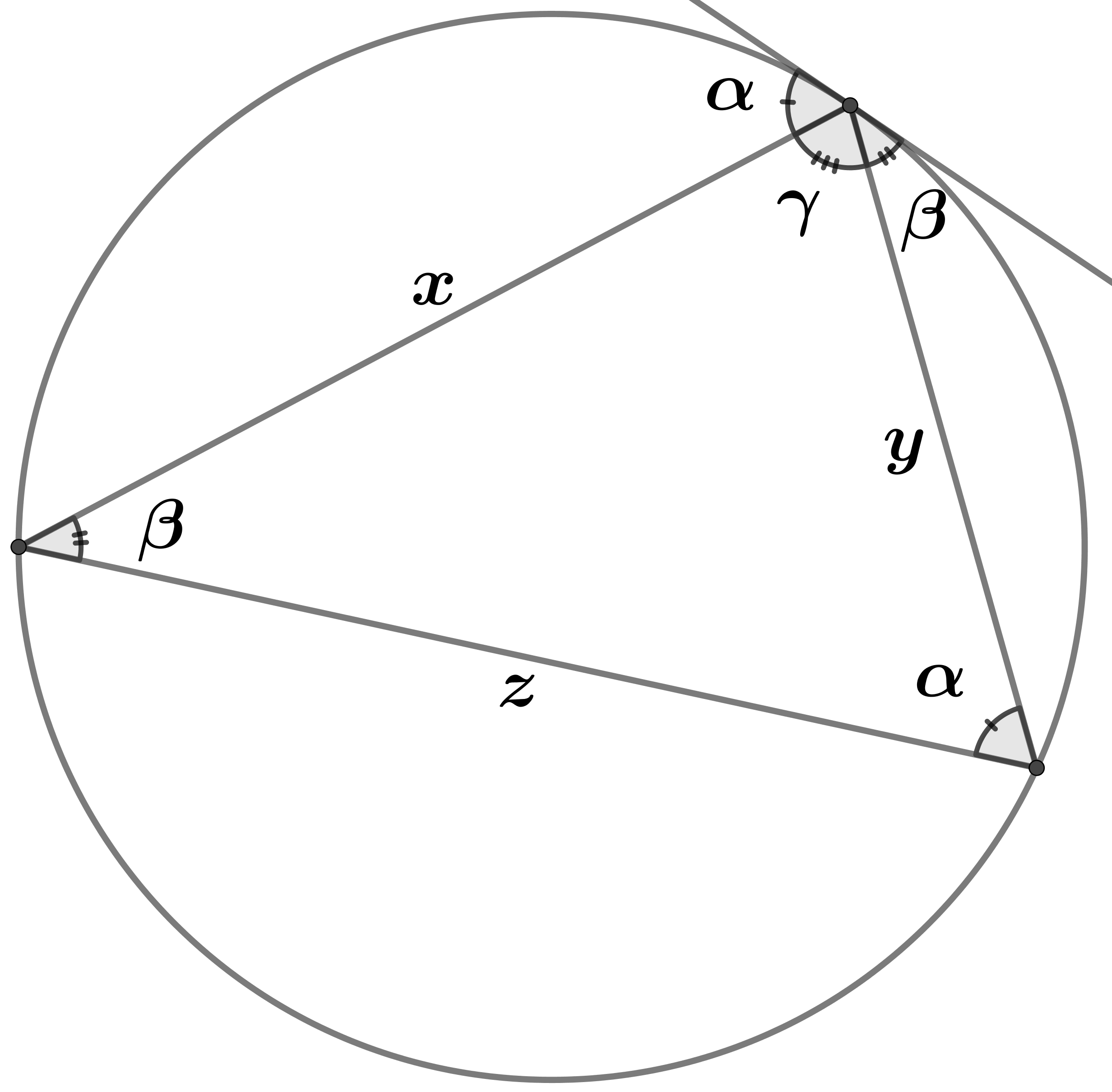}}%
}%

\vspace{2mm}

After the infinitesimal displacement of our vertex along the circumcircle, we consider the right triangles with $\delta \ell$-length hypotenuse and $\delta x$-length leg, and with $\delta \ell$-length hypotenuse and $-\delta y$-length leg, respectively, and derive
\begin{equation}
\delta x =\cos\alpha\,\delta \ell 
\qquad -\delta y=\cos\beta\,\delta \ell
\label{DeviationsCircum0}
\end{equation}

\vspace{2.5mm}
\noindent\parbox{86mm}{
\includegraphics[width=43mm,keepaspectratio]{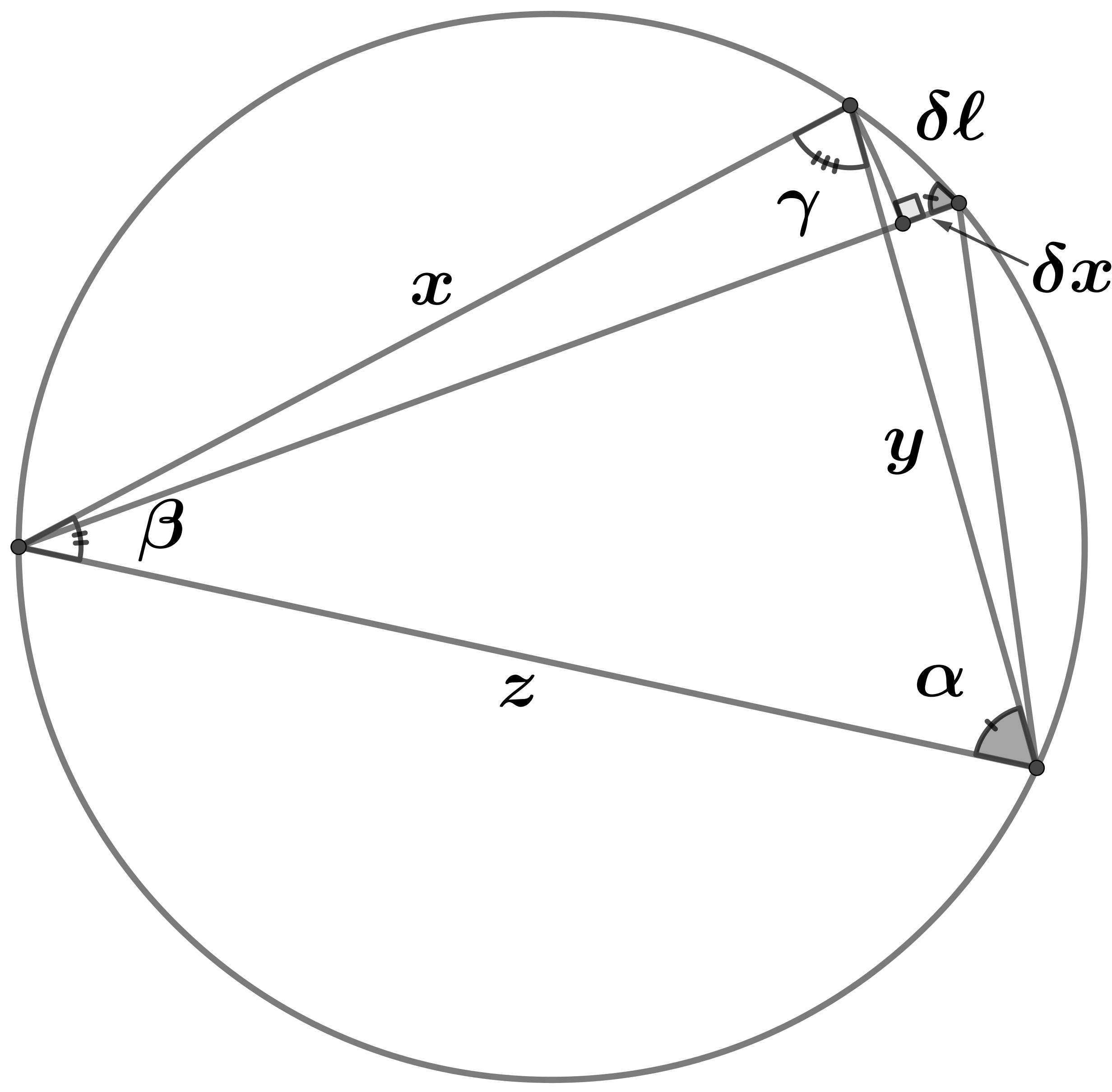}%
\includegraphics[width=43mm,keepaspectratio]{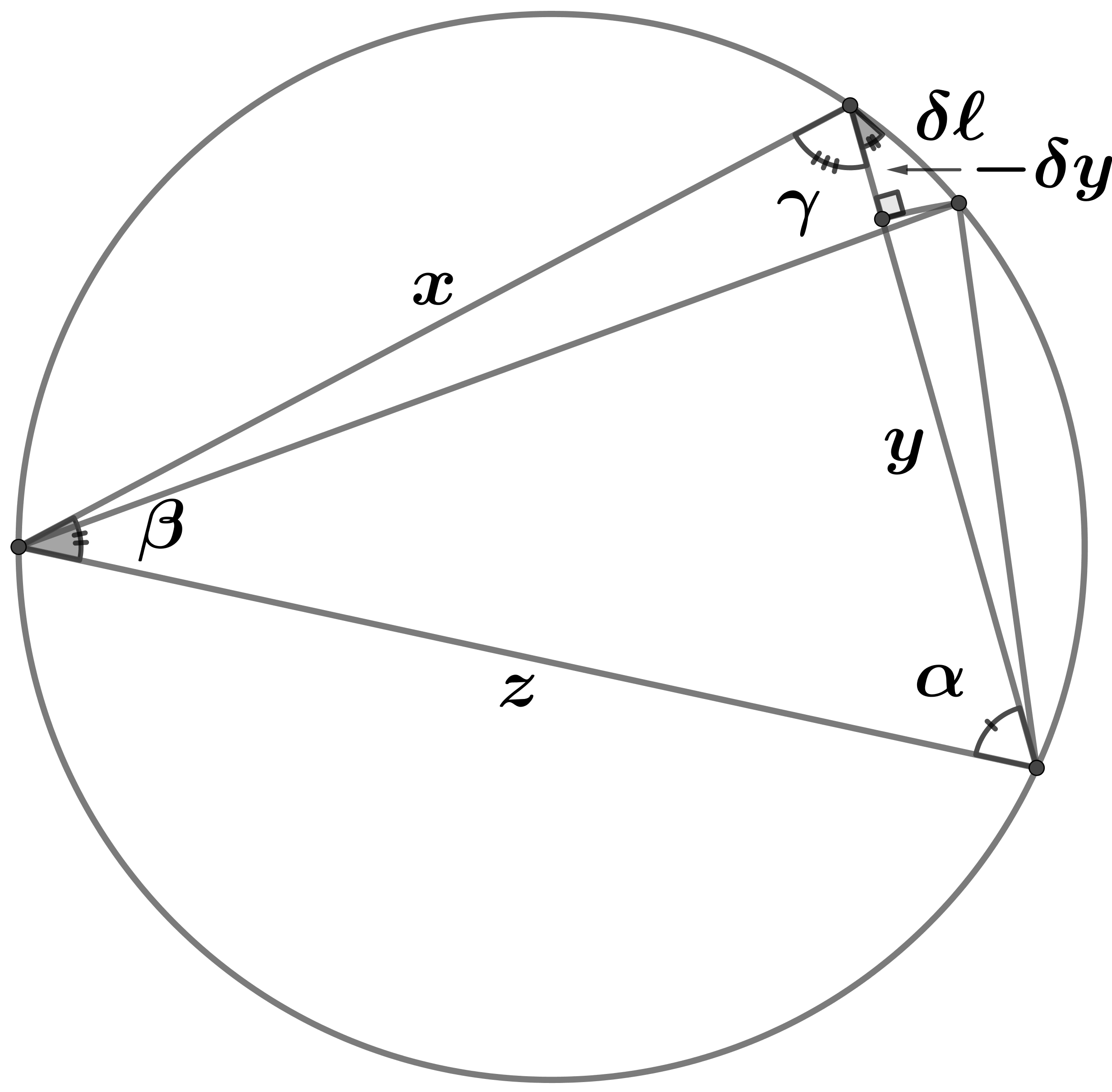}%
}%
\vspace{2.5mm}

\noindent Applying the al-Kashi theorem yields
\begin{equation}
\delta x =\frac{y^2-x^2+z^2}{2 y z}\delta \ell 
\qquad -\delta y=\frac{x^2-y^2+z^2}{2 x z}\delta \ell
\label{DeviationsCircum1}
\end{equation}

\noindent Plugging in results (\ref{DeviationsCircum1}) into eq. (\ref{PartialCircum}), simplifying by $\delta\ell/2z$ and isolating $\partial_y R$ gives the partial differential equation
\begin{equation}
\partial_y R=\frac{x}{y}\,\frac{y^2-x^2+z^2}{x^2-y^2+z^2}\,\partial_x R 
\label{PartialyCircum}
\end{equation}

\noindent A similar expression arises for $\partial_z R$ when the $(x,z)$ vertex is slightly moved along the circumcircle: 
\begin{equation}
\partial_z R=\frac{x}{z}\,\frac{z^2-x^2+y^2}{x^2-z^2+y^2}\,\partial_x R 
\label{PartialzCircum}
\end{equation}

\noindent Finally, we consider the infinitesimal scale transformation $R\mapsto R+\delta R$. Since each of the sidelengths increases in line with itself and $\delta R/R$, that is $\delta x= x/R\,\delta R$, etc., we have the partial differential equation
\begin{equation}
R=x\,\partial_x R + y\,\partial_y R + z\,\partial_z R
\label{ScaleCircum}
\end{equation}

\noindent which guarantees that each quantity at stake is taken into account with its appropriate dimension. Plugging in eq. (\ref{PartialyCircum}) and (\ref{PartialzCircum}) into the last one (\ref{ScaleCircum}) yields the partial differential equation
\begin{equation}
\frac{\partial_x R}{R}=\frac{x^4-y^4-z^4+2 y^2 z^2}{x[-x^4-y^4-z^4+2(y^2 z^2 + x^2 z^2 + x^2 y^2)]}
\end{equation}

\noindent which can be integrated out to give the general solution
\begin{equation}
R(x,y,z)
=\frac{x\,k(y,z)}{\sqrt{-x^4-y^4-z^4+2(y^2 z^2 + x^2 z^2 + x^2 y^2)}}
\end{equation}

\noindent with $k(y,z)$ a symmetric function of $y$ and $z$. Since $R(x,y,z)$ must itself be symmetric in $x$, $y$ and $z$, $k(y,z)=c\,yz$ where $c$ is a constant. We just have to compute the circumradius length of, say, an equilateral triangle of sidelength $1$, that is $R(1,1,1)=\sqrt{3}/3$, to find out that $c=1$. Furthermore, like for Heron's theorem, the argument of the square root can be factorized to give 
\begin{align}
&R(x,y,z)\\
&=\frac{x y z}{\sqrt{(x+y+z)(-x+y+z)(x-y+z)(x+y-z)}}
\nonumber
\end{align}

\noindent which is indeed the expression of the circumradius as a function of the sidelengths.

\section{The Law of Sines} 
\label{Sines}

\noindent Let us go back to the Thales theorem that Newton could have discovered by the grace of calculus if only he was born two millenia before: $y(x)=kx$ is already a nice result, that we have used all along this paper, but what about $k$? If he knew the angle $\gamma=\widehat{(x,y)}$ between the two sides, and one of the angles, say $\beta=\widehat{(z,x)}$, between the $x$-length side and the third side, of length $z$, he could have speculated that $y$ is a smooth function of $x$, $\beta$ and $\gamma$, that is \vspace{-2mm}
\begin{equation}
y=y(x,\beta,\gamma)
\end{equation}

\noindent\parbox{86mm}{\parbox{41mm}{\noindent It is possible to find the function $y$ with a simple rotation of one side around one of the adjacent vertices, like in sections \ref{Apollonius} to \ref{Terquem}, but the differential equations that arise are a bit complicated to solve: one is linear and the other one of the Bernouilli type. Even if it is not straightforward, it happens to be much simpler to work within the circumcircle. After an infinitesimal and clockwise displacement of the $\gamma$ vertex along the circumcircle, resulting in a small deviations $\delta y$, $\delta x$ and $\delta\beta$ while $\delta\gamma=0$, at first order,  }\hspace{2mm}%
\parbox{43mm}{

\includegraphics[width=43mm,keepaspectratio]{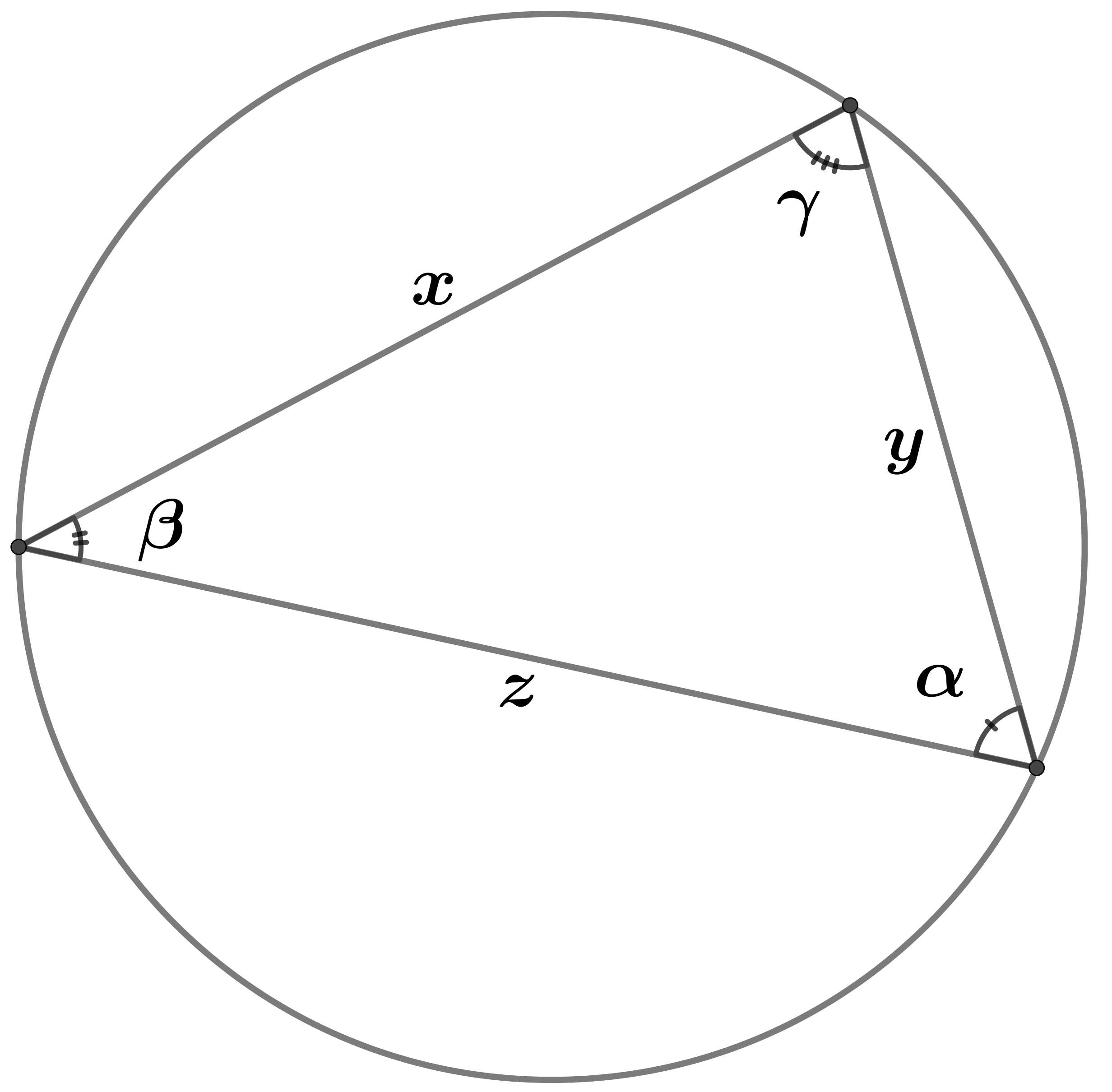}

\includegraphics[width=43mm,keepaspectratio]{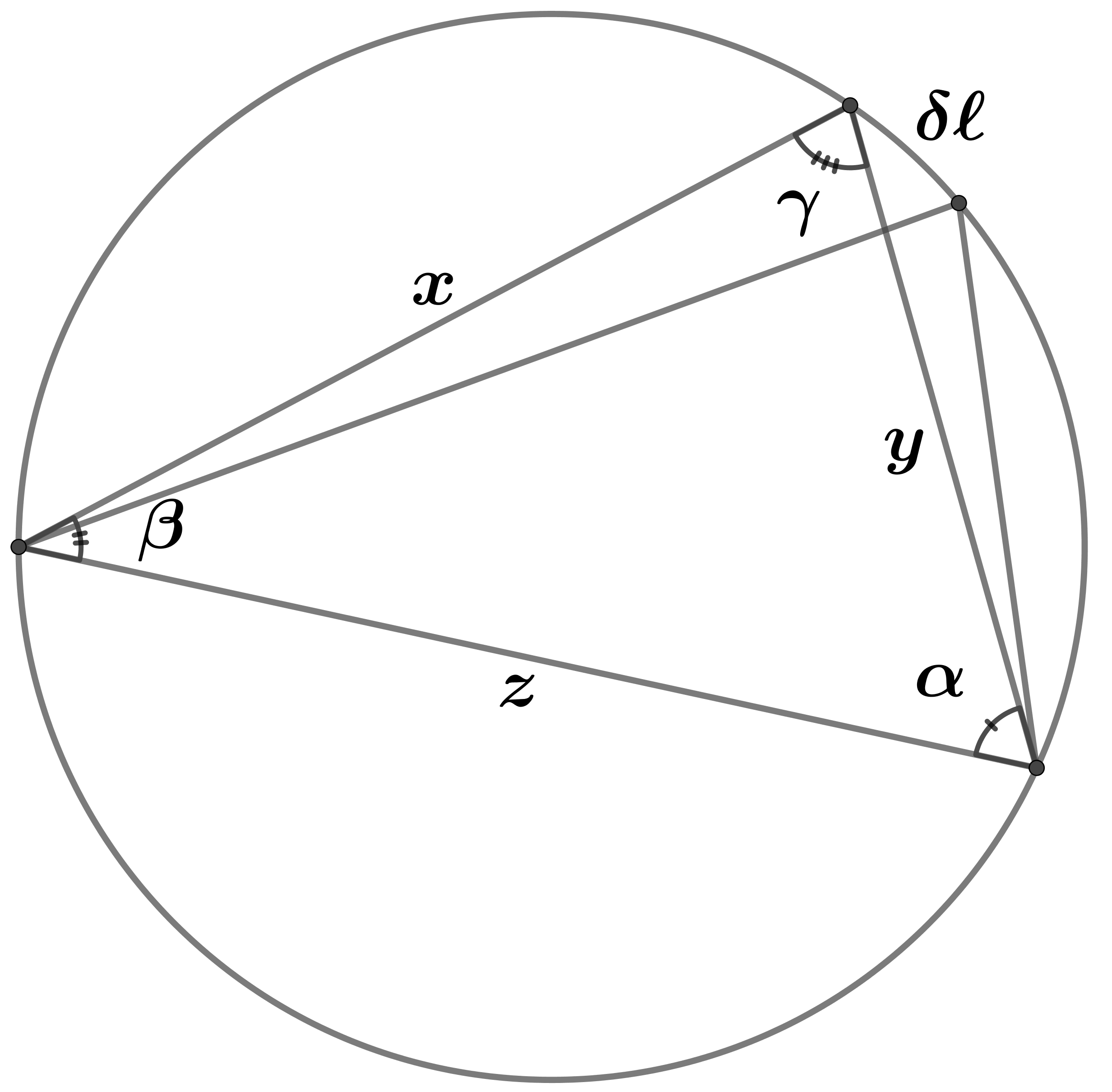}
}%
}%

\begin{equation}
\delta y=\partial_x y\,\delta x + \partial_{\beta} y\,\delta\beta
\label{PartialSines0}
\end{equation}

\vspace{1.5mm}
\noindent\parbox{86mm}{
\includegraphics[width=43mm,keepaspectratio]{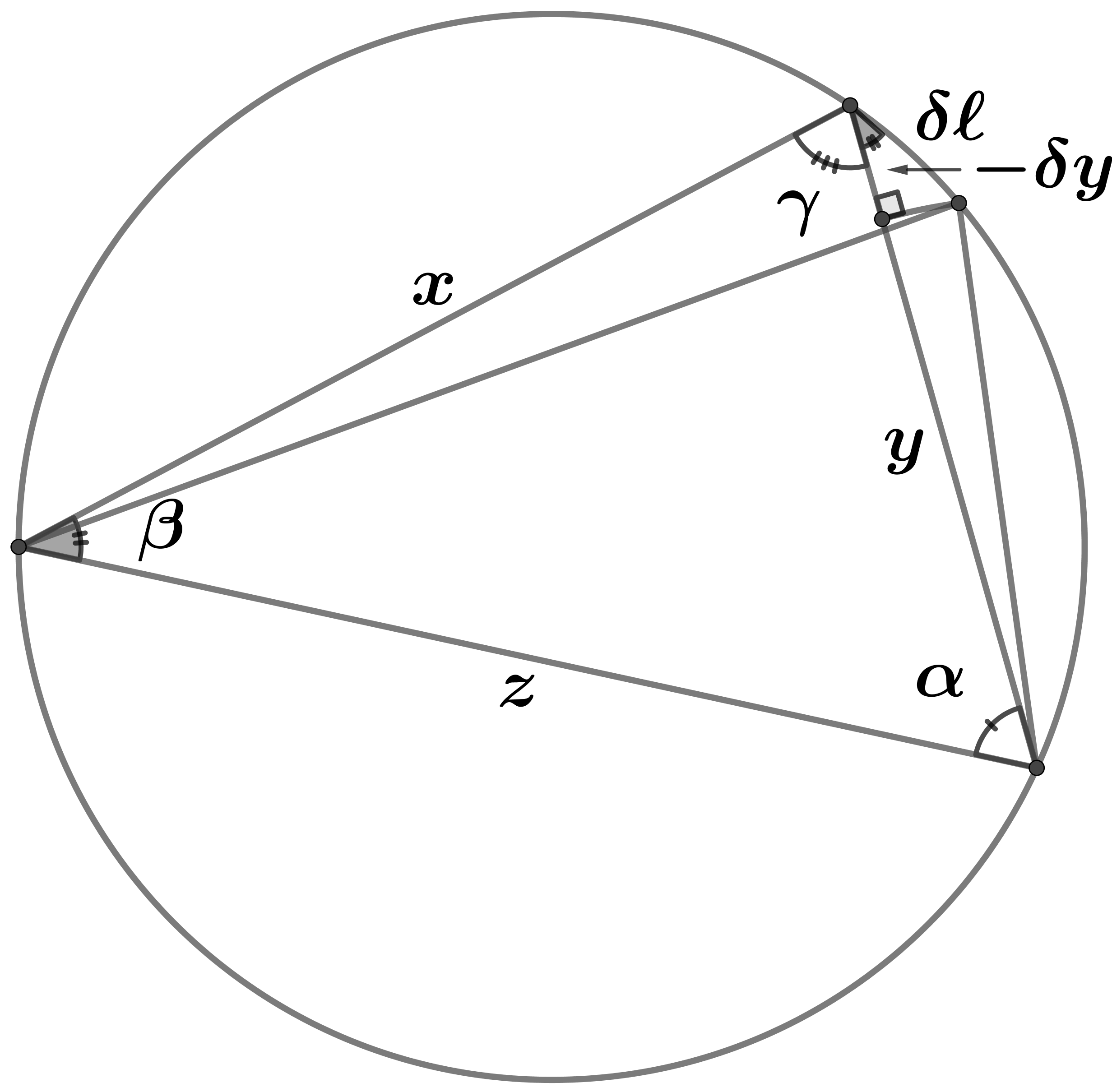}%
\includegraphics[width=43mm,keepaspectratio]{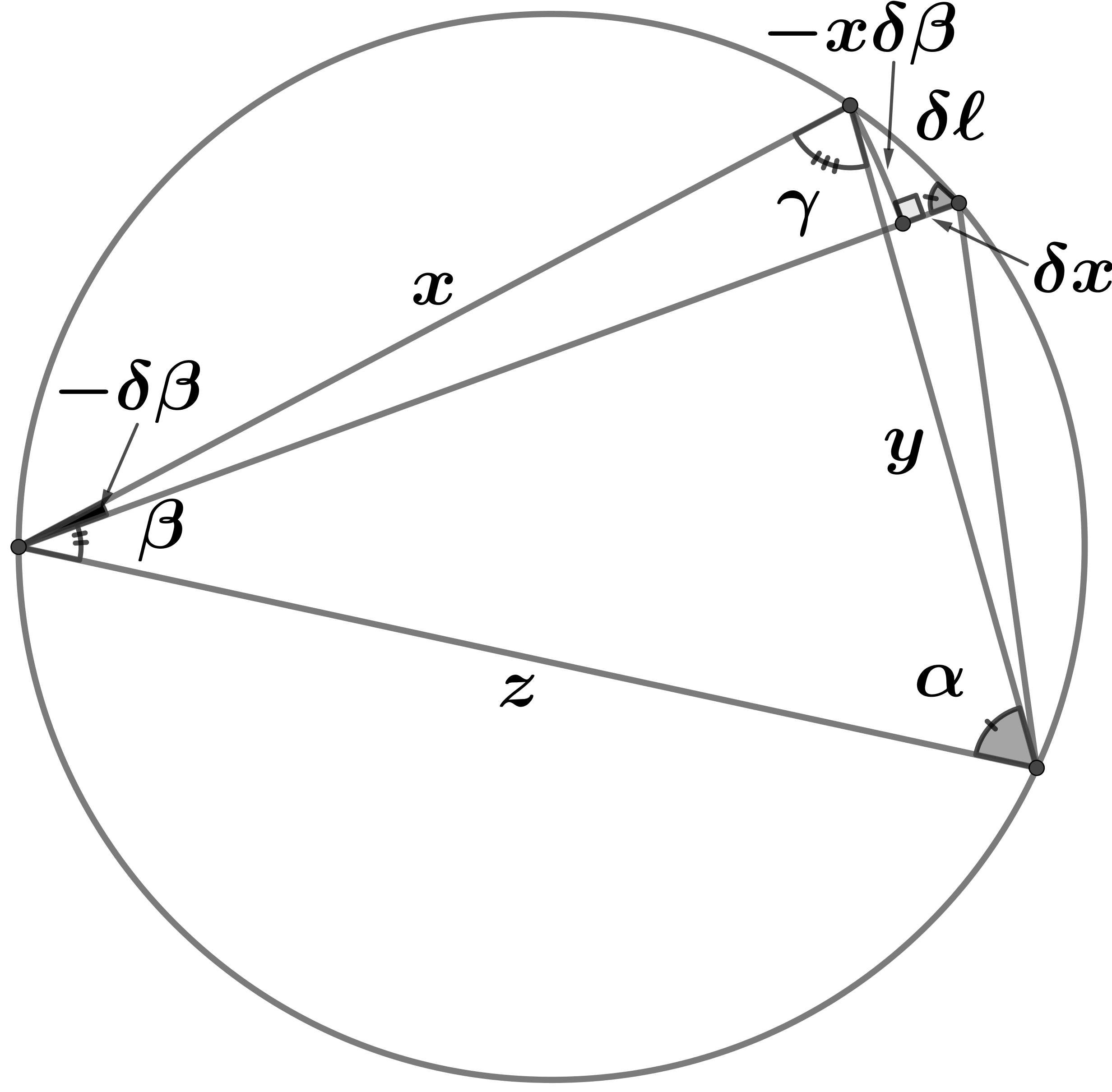}%
}%
\vspace{2.5mm}

\noindent The deviations $\delta x$ and $\delta y$ are the same as in section \ref{Circumradius}, eq. (\ref{DeviationsCircum0}). $-\delta\beta$ can be computed by considering $-x\delta\beta$ as the length of the leg opposite to $\alpha$ in the small right triangle of other leg of length $\delta x$ and hypotenuse of length $\delta\ell$. Using $\alpha=\pi-(\beta+\gamma)$ to get rid of $\alpha$, we have at first order
\begin{equation}
-\delta y\!=\!\cos\beta\,\delta \ell\quad\!\delta x\!=\!-\cos(\beta+\gamma) \delta\ell\quad\! -\delta\beta\!=\!\frac{\sin(\beta+\gamma)}{x}\delta\ell
\label{DeviationsSines1}
\end{equation}

\noindent Combining eq. (\ref{PartialSines0}) and (\ref{DeviationsSines1}) leads to the differential equation
\begin{equation}
x\cos\beta=x\cos(\beta+\gamma)\partial_x y+\sin(\beta+\gamma)\partial_{\beta}y
\label{PartialSines1}
\end{equation}

\noindent If we now consider an infinitesimal and clockwise displacement of the $\beta$ vertex along the circumcircle, resulting in deviations $\delta x$ and $\delta\gamma$ while $\delta y=0$ and $\delta\gamma=0$, at first order  \vspace{-2mm}
\begin{equation}
\partial_x y\,\delta x + \partial_{\gamma} y\,\delta\gamma=0
\label{PartialSines2}
\end{equation}

\noindent Using the analogous first-order similarities and the fact that $\alpha=\pi-(\beta+\gamma)$, we find
\begin{equation}
-\delta x=-\cos(\beta+\gamma)\,\delta \ell\qquad \delta\gamma=\frac{\sin(\beta+\gamma)}{x}\,\delta\ell
\label{DeviationsSines2}
\end{equation}

\noindent Combining eq. (\ref{PartialSines2}) and (\ref{DeviationsSines2}) leads to the differential equation
\begin{equation}
x\cos(\beta+\gamma)\partial_x y+\sin(\beta+\gamma)\partial_{\gamma}y=0
\label{PartialSines3}
\end{equation}

\noindent Finally, we consider the infinitesimal scale transformation $R\mapsto R+\delta R$, where $R$ is the circumradius. Since each of the sidelengths increases in line with itself and $\delta R/R$, that is $\delta x= x/R\,\delta R$, but $\delta\beta=0$ and $\delta\gamma=0$, we have the partial differential equation
\begin{equation}
y=x\,\partial_x y
\label{ScaleSines}
\end{equation}

\noindent Using eq. (\ref{ScaleSines}) allows to simplify eq. (\ref{PartialSines1}) and (\ref{PartialSines3}) to yield
\begin{align}
\partial_{\beta}y+\frac{y}{\tan(\beta+\gamma)}&=\frac{x\cos{\beta}}{\sin(\beta+\gamma)}
\label{PartialSines4}\\
\partial_{\gamma}y+\frac{y}{\tan(\beta+\gamma)}&=0
\label{PartialSines5}
\end{align}

\noindent They are both linear. Moreover, the second one is separable and admits the general solution
\begin{equation}
y(x,\beta,\gamma)=\frac{K(x,\beta)}{\sin(\beta+\gamma)}
\label{GeneralSines1}
\end{equation}

\noindent where $K(x,\beta)$ is a function of $x$ and $\beta$. But this solution must satisfy eq. (\ref{ScaleSines}), that is $K(x,\beta)=k(\beta)x$ where $k(\beta)$ is a function of $\beta$. The solution must also satisfy eq. (\ref{PartialSines4}), thus we find $k(\beta)=\sin\beta+C$ with $C$ a real constant. Finally, since for $\beta=0$, $y=0$, we have
\begin{equation}
y(x,\beta,\gamma)=\frac{x\sin\beta}{\sin(\beta+\gamma)}
\label{LawOfSines}
\end{equation}

\noindent which, using $\sin(\beta+\gamma)=\sin\alpha$, is the law of sines, proved by Nasir al-Din al-Tusi, Persian mathematician -- and architect -- regarded as a founder of trigonometry \cite{alTusi:1260}.

\section{Claudius Ptolemy} 
\label{Ptolemy}

\vspace{2mm}

\noindent\parbox{86mm}{\parbox{41mm}{\noindent We now take Newton back to Antiquity, in Roman Egypt, and make him look at an ancient problem: given a cyclic quadrilateral, determine the length of its diagonals as functions of the sidelengths $x$, $y$, $u$ and $v$, that is, for the $z$-length segments joining the $(v,x)$ and $(y,u)$ vertices: }\hspace{2mm}%
\parbox{43mm}{\includegraphics[width=43mm,keepaspectratio]{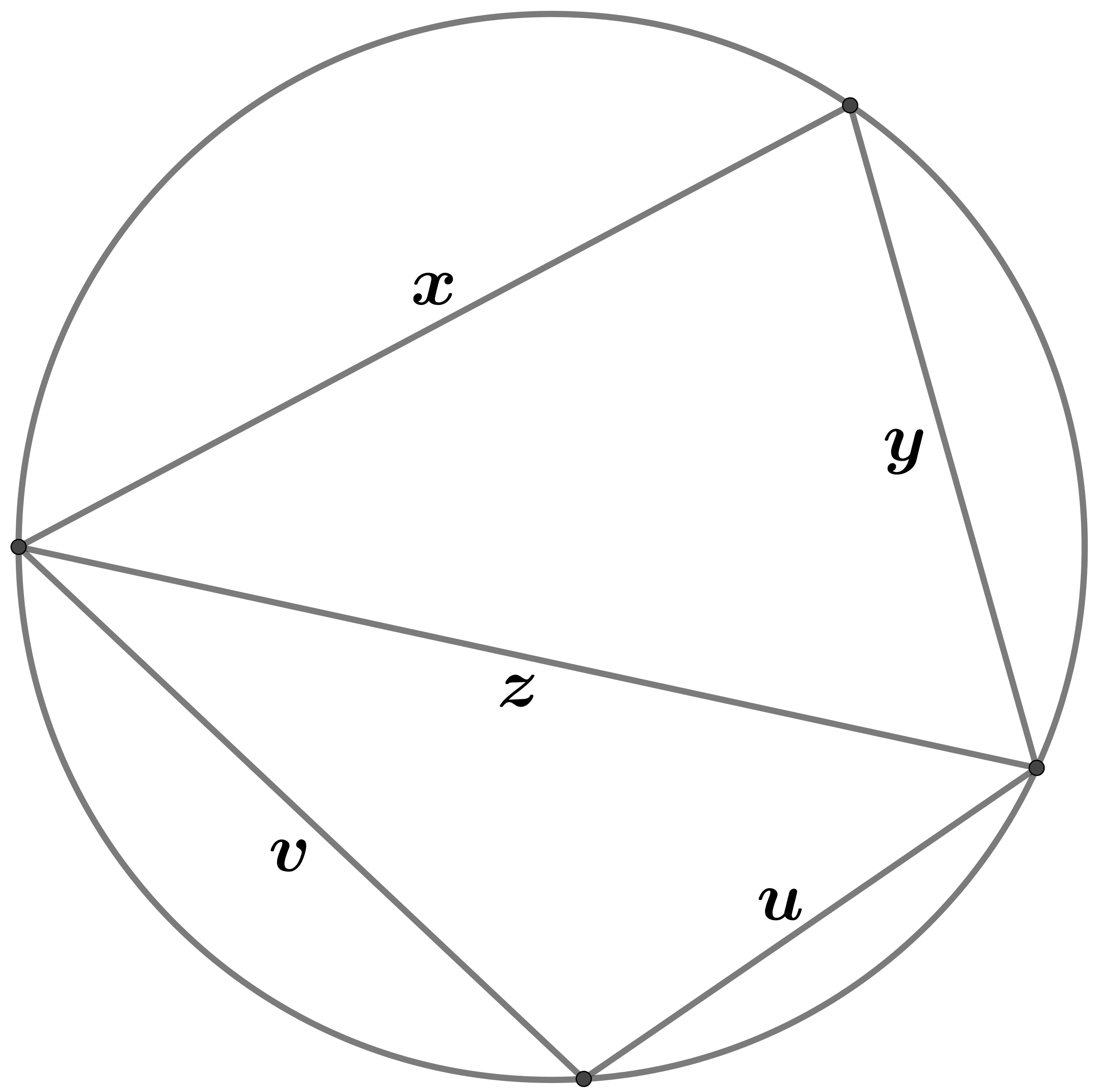}}%
}%
\begin{equation}
z=z(x,y,u,v)
\end{equation}

\noindent After a small displacement of the $(x,y)$ vertex along the circumcircle, resulting in deviations $\delta x$ and $\delta y$ while $\delta u=\delta v=0$ as well as $\delta z=0$, we have at first order
\begin{equation}
\partial_x z\,\delta x + \partial_{y} z\,\delta y=0
\label{PartialPtolemy0}
\end{equation}

\noindent The deviations $\delta x$, $\delta y$ are those of section \ref{Circumradius}, eq. (\ref{DeviationsCircum0}): 
\begin{equation}
\delta x =\cos\alpha\,\delta \ell 
\qquad -\delta y=\cos\beta\,\delta \ell
\label{DeviationsPtolemy0}
\end{equation}

\noindent where $\alpha$ and $\beta$ are the angles opposite to the $x$ and $y$-length sides of the $(x,y,z)$ triangle. 

\vspace{2.5mm}
\noindent\parbox{86mm}{
\includegraphics[width=43mm,keepaspectratio]{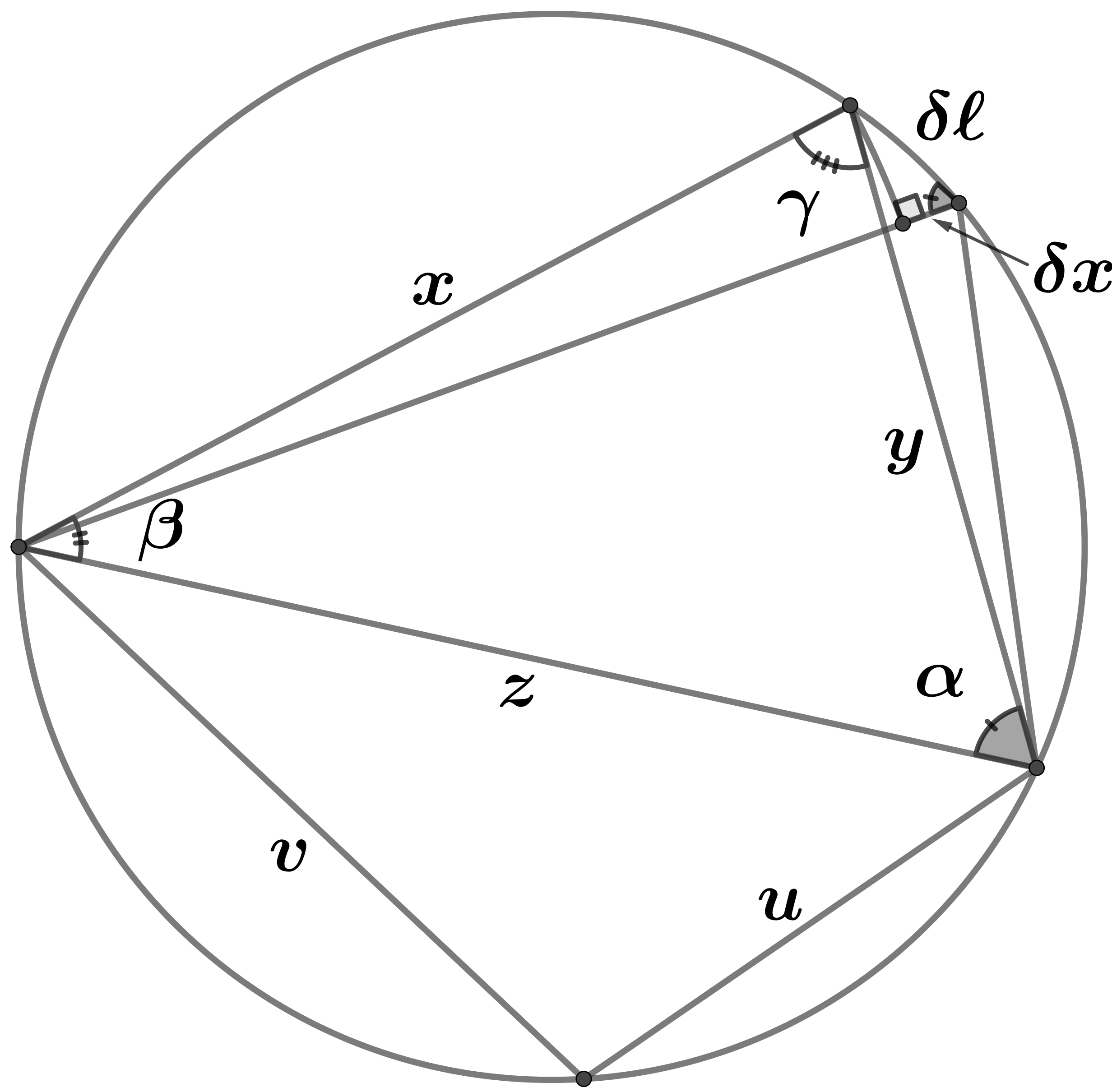}%
\includegraphics[width=43mm,keepaspectratio]{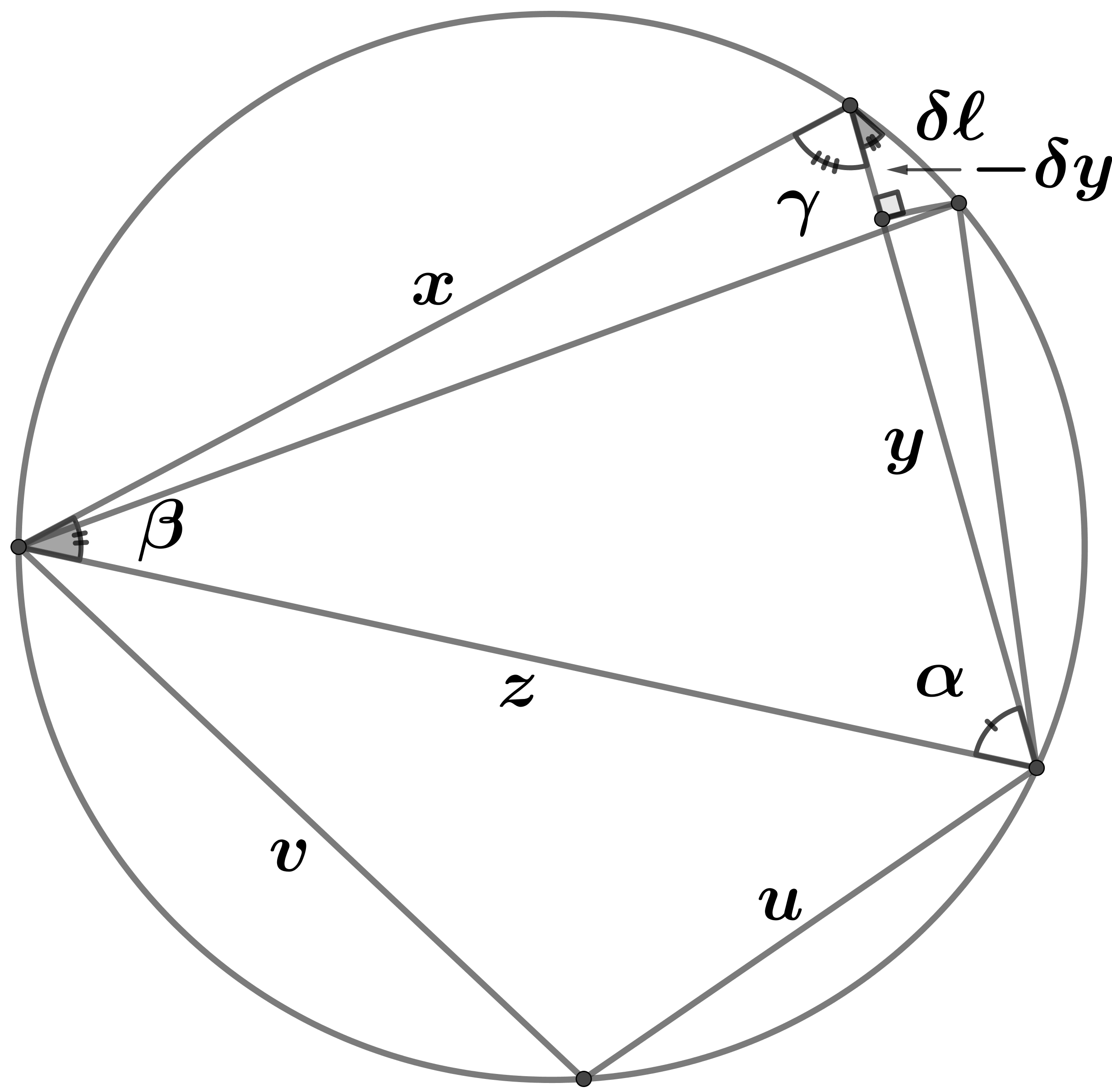}%
}%
\vspace{2.5mm}

\noindent Plugging in those results into eq. (\ref{PartialPtolemy0}) and considering the equivalent after a small displacement of the $(u,v)$ vertex, we find the two partial differential equations
\begin{align}
\label{PartialPtolemy1}
\cos\alpha\,\partial_x z -\cos\beta\,\partial_y z&=0\\
\cos\rho\,\partial_u z -\cos\phi\,\partial_v z&=0
\label{PartialPtolemy2}
\end{align}

\noindent\parbox{86mm}{\parbox{43mm}{

\vspace{2,5mm}

\includegraphics[width=43mm,keepaspectratio]{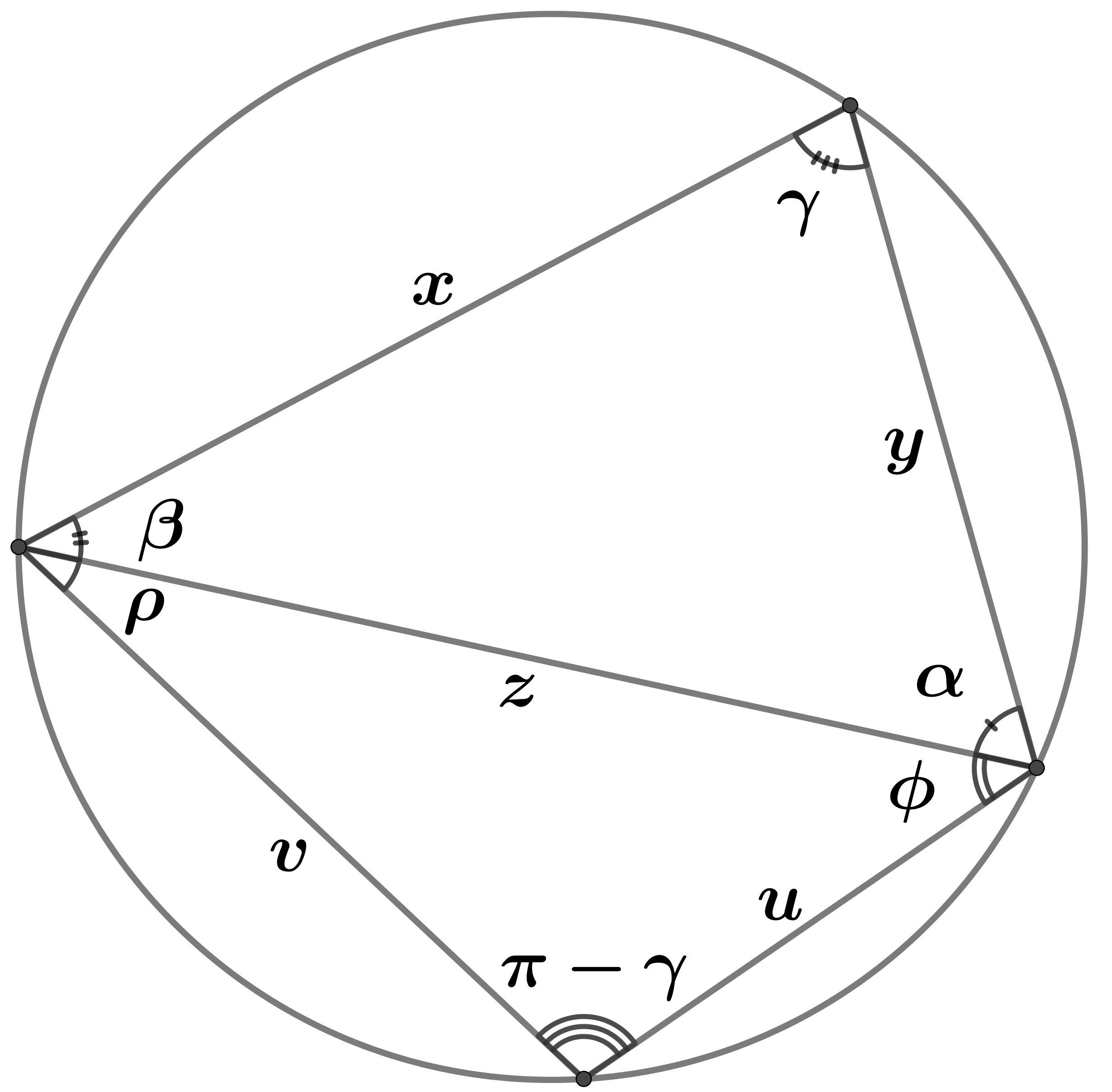}

\includegraphics[width=43mm,keepaspectratio]{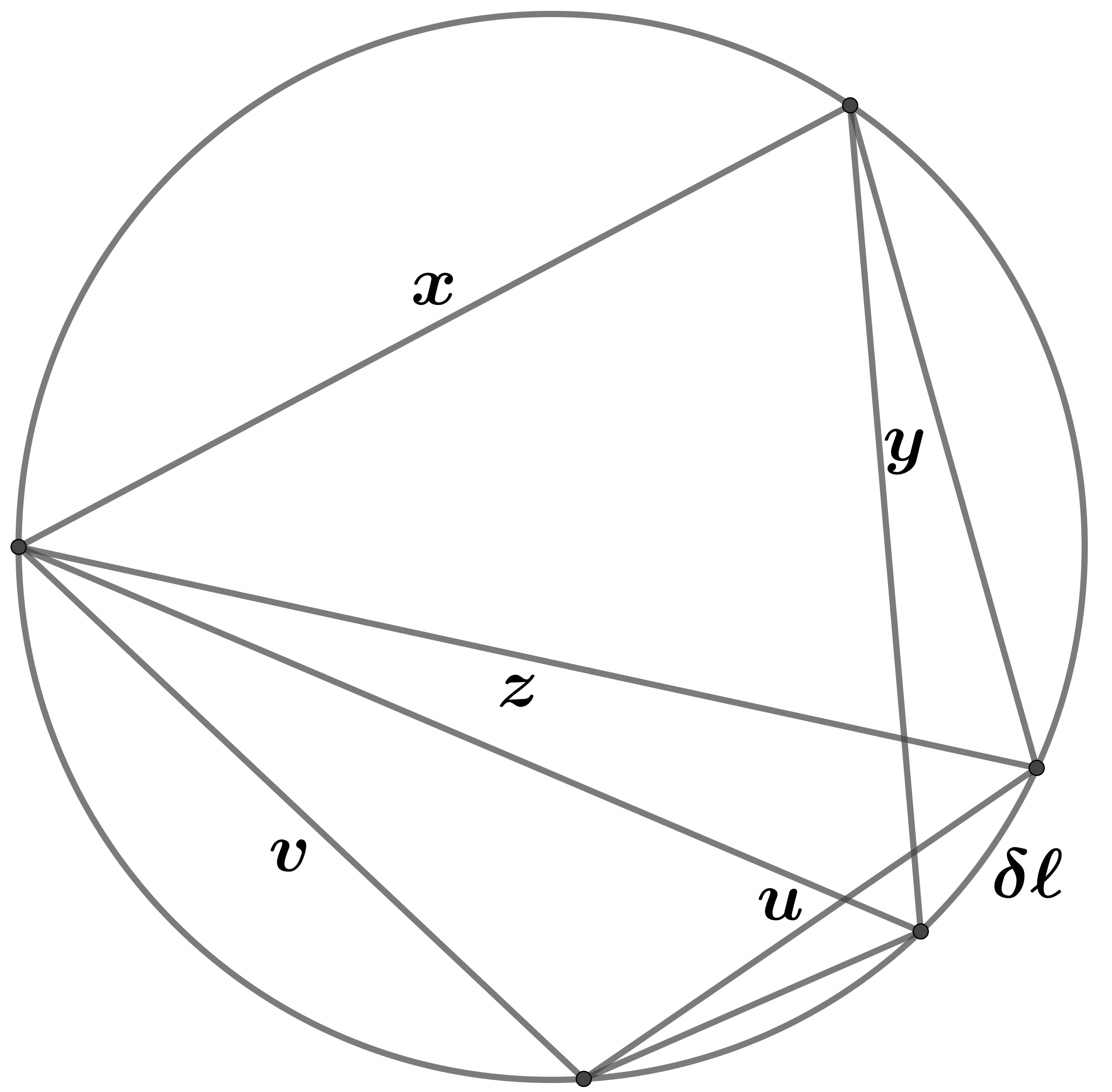}

\includegraphics[width=43mm,keepaspectratio]{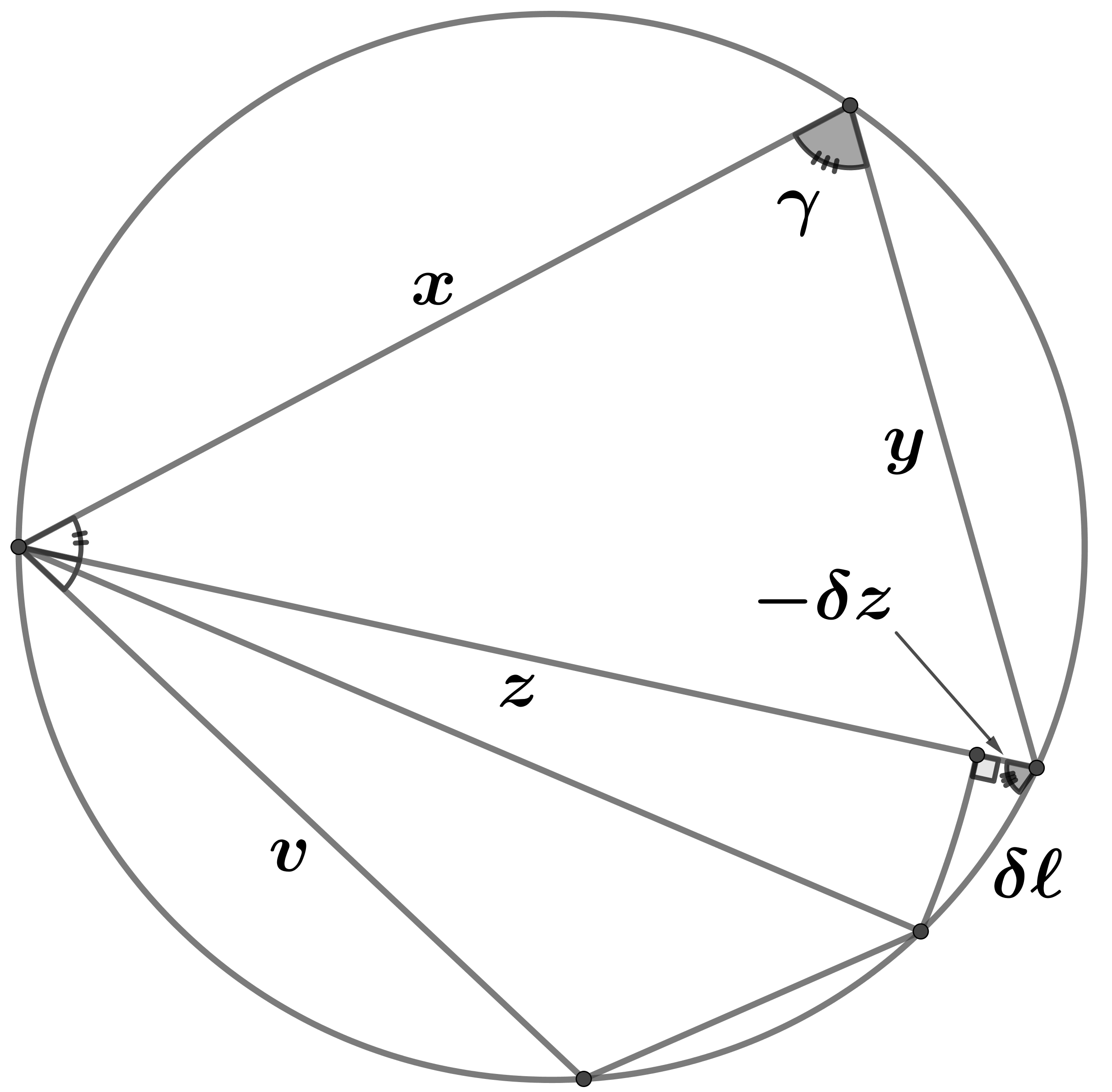}

\vspace{2,5mm}
}\hspace{2mm}%
\parbox{41mm}{
\noindent where $\rho$ and $\phi$ are the angles opposite to the $u$ and $v$-length sides of the $(u,v,z)$ triangle. 

\hspace{2mm} If we now slightly move the $(y,u)$ vertex along the circumcircle, with small deviations $\delta y$, $\delta u$ and $\delta z$, while $\delta x=\delta v=0$, we have at first order
\begin{equation}
\delta z=\partial_y z\,\delta y + \partial_{u} z\,\delta u
\label{PartialPtolemy3}
\end{equation}
\noindent The deviations $\delta y$, $\delta u$ and $\delta z$ are given by
\begin{align}
\delta y &=\cos\beta\,\delta \ell \nonumber\\
 -\delta u&=\cos\rho\,\delta \ell\nonumber\\
 -\delta z&=\cos \gamma \delta \ell
\label{DeviationsPtolemy1}
\end{align}
\noindent where $\delta y$ and $-\delta u$ are calculated analogously to $\delta x$ and $-\delta y$ above, while $-\delta z$ is to be seen as the length of one of the legs in a small right triangle of hypotenuse of length $\delta \ell$ and adjacent angle of same magnitude $\gamma$ as the one between the $x$ and $y$-length sides. 
}
}%

\noindent Inserting those results into eq. (\ref{PartialPtolemy3}) yields
\begin{equation}
\cos\gamma=\cos\rho\,\partial_u z -\cos\beta\,\partial_y z
\label{PartialPtolemy4}
\end{equation}

\noindent Moreover, from the infinitesimal scale transformation $R\mapsto R+\delta R$, a fourth differential equation arises:
\begin{equation}
z=x\,\partial_x z + y\,\partial_y z + u\,\partial_u z + v\,\partial_v z
\label{ScalePtolemy}
\end{equation}

\noindent Combining equations (\ref{PartialPtolemy1}), (\ref{PartialPtolemy2}), (\ref{PartialPtolemy4}) and (\ref{ScalePtolemy}), isolating $\partial_x z$ and using $x\cos\beta+y\cos\alpha=z$ and $u\cos\rho+v\cos\phi=z$ to simplify the expression, we have
\begin{equation}
\partial_x z=\frac{(\cos\rho\cos\phi-\cos\gamma)\cos\beta}{\cos\rho\cos\phi+\cos\alpha\cos\beta}
\end{equation}

\noindent Looking at the intercepted arcs, we observe that $\widehat{(u,v)}$ and $\gamma$ are supplementary. We then exploit $\gamma=\pi-(\alpha+\beta)=\rho+\phi$ and the angle sum identities to find
\begin{equation}
\partial_x z=\frac{\sin\rho\sin\phi\cos\beta}{\sin\rho\sin\phi+\sin\alpha\sin\beta}
\end{equation}

\noindent Thanks to the law of sines applied to both $(x,y,z)$ and $(u,v,z)$ triangles, and the al-Kashi theorem applied to $\cos\beta$, this equation finally becomes
\begin{equation}
2z\partial_x z=\frac{u v (x^2-y^2+z^2)}{(u v + x y) x}
\end{equation}

\noindent and turns out to be linear in $z^2$. Its general solution reads
\begin{equation}
z^2(x,y,u,v)=\frac{u v(x^2+y^2)+c(y,u,v)x }{u v + x y }
\end{equation}

\noindent where the function $c(y,u,v)$ must be symmetric in $u$ and $v$, while $z^2(x,y,u,v)$ must notably be symmetric in $x$ and $y$. Thus $c(y,u,v)=y(u^2+v^2)$. Hence, after factorization, 
\begin{equation}
z^2(x,y,u,v)=\frac{(u x + v y)(v x + u y)}{u v + x y }
\label{PtolemyTheorem}
\end{equation}

\noindent which, if multiplied by the square of the other diagonal length, yields the theorem of Ptolemy \cite{Ptolemy:150AD}. Eq. (\ref{PtolemyTheorem}) might have been first stated by Brahmagupta \cite[\S 28]{Brahma:628}.

\section{Brahmagupta} 
\label{Brahma}

\noindent A few centuries later in India, was solved the problem of expressing the area $A$ of any cyclic quadrilateral as a function of its sidelengths $x$, $y$, $u$ and $v$. Well it could have been solved through calculus too. Let us start from
\begin{equation}
A=A(x,y,u,v)
\end{equation}

\vspace{2.5mm}
\noindent\parbox{86mm}{
\includegraphics[width=43mm,keepaspectratio]{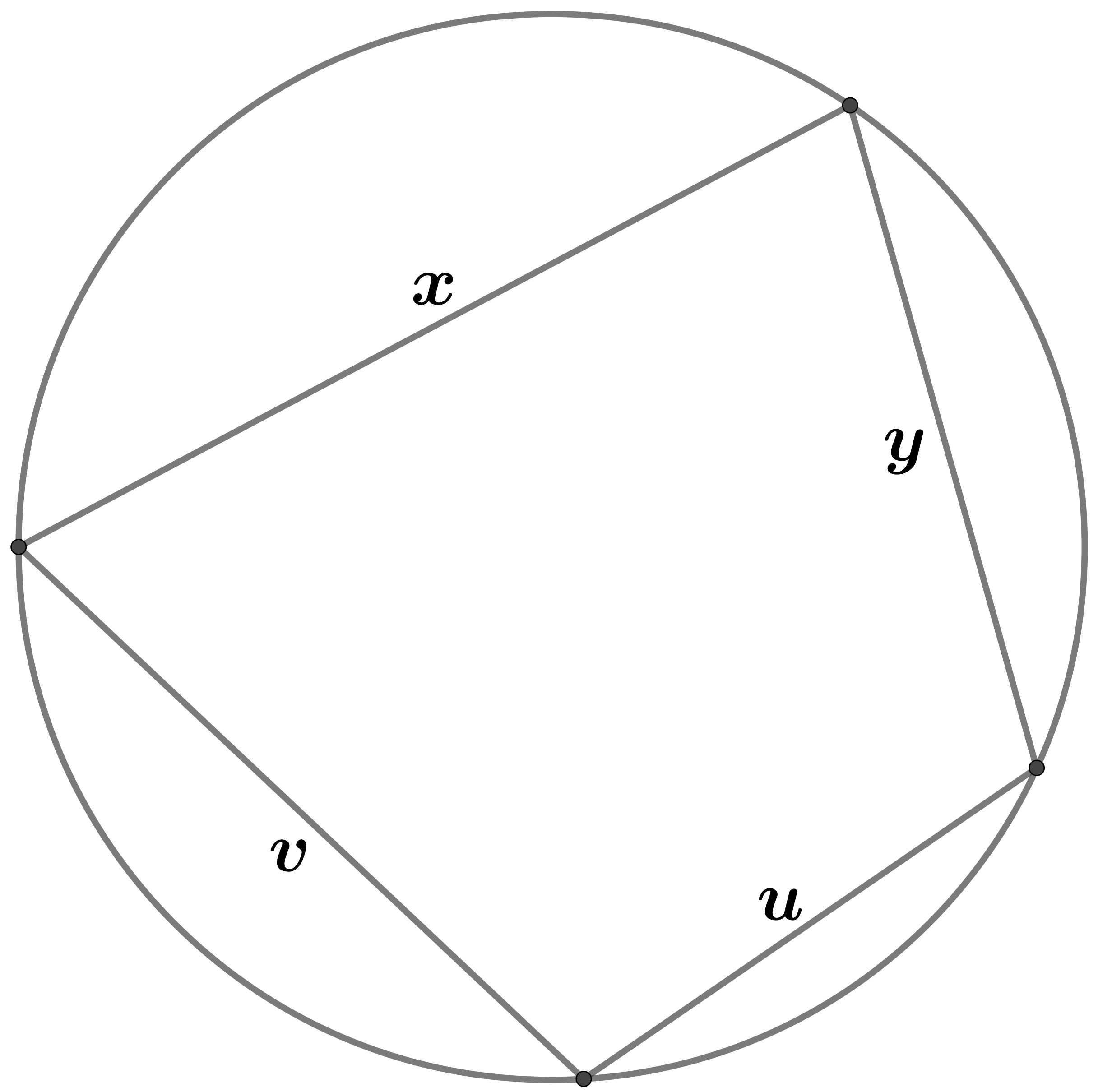}%
\includegraphics[width=43mm,keepaspectratio]{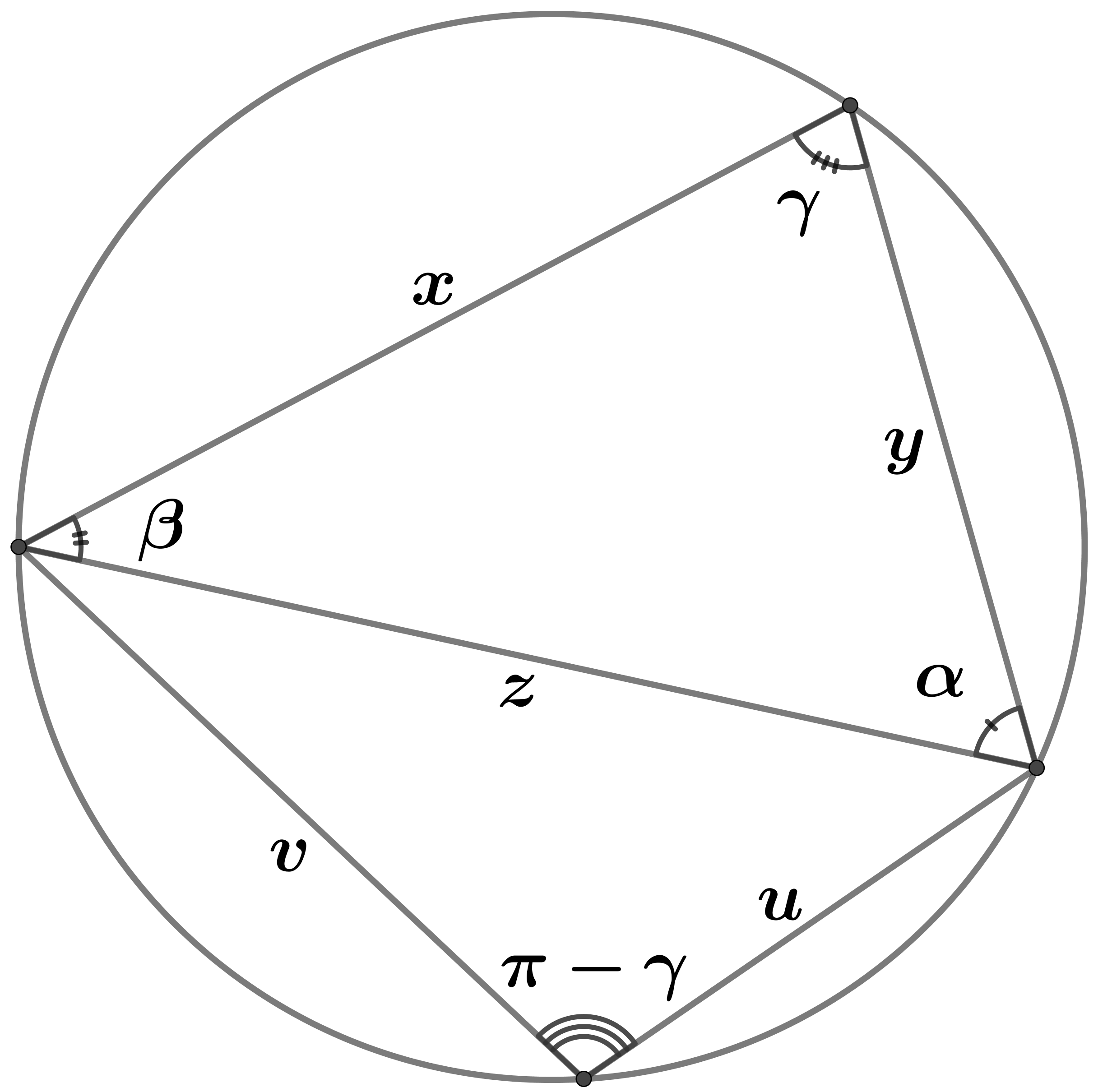}%
}%
\vspace{2.5mm}

\noindent Like in the last section, we draw the $z$-length diagonal, which divides the quadrilateral in two triangles of sides of length $x$, $y$, $z$, and $u$, $v$, $z$ respectively. Selecting the $y$-length side as the base of the $(x,y,z)$ triangle, we compute its area by multiplying $y$ and $x\sin\gamma$ divided by two, where $\gamma$ is the $\widehat{(x,y)}$ angle magnitude. But the intercepted arcs of the circumcircle tell us that in the other triangle, $\widehat{(u,v)}$ is the supplementary of $\widehat{(x,y)}$. Thus we have for their respective areas $A_{xyz}=xy\sin\gamma/2$ and $A_{uvz}=uv\sin\gamma/2$. This yields
\begin{equation}
\sin\gamma=\frac{2A}{xy+uv}
\label{SineBrahma}
\end{equation}

\noindent where $A=A_{xyz}+A_{uvz}$ is the total area of the quadrilateral. 

Let us go back to the project. Rather than moving the $(x,y)$ vertex along the circumcirle, which leads to complicated equations, we invite Newton to follow the reasoning developed in section \ref{Heron} for the Heron theorem and make an infinitesimal rotation of the $y$-length side around the $(y,z)$ vertex resulting in a small deviation $\delta x$, while $\delta y=0$ and $\delta z=0$. We have at first order
\begin{equation}
\delta A =\partial_x A\,\delta x
\label{PartialBrahma0}
\end{equation}

\vspace{2mm}

\noindent\parbox{86mm}{\parbox{41mm}{\noindent This slight departure from the cyclical nature of the quadrilateral will be of no consequence. The reasoning is valid in the triangle and as soon as the needed result is secured, we will get back into the circumcircle to proceed. Here $h$ will stand for the length of the height relative to the $z$-length side in the $(x,y,z)$ triangle, such that $A_{xyz}=zh/2$ and thus $\delta A=z\delta h/2$ since $\delta A_{uvz}=0$. We know a little more about trigonometry than in section \ref{Heron}, so that we can express, at first order, the deviations $\delta x$ and $\delta h$ in terms of $\gamma$ and $\alpha$, the $\widehat{(y,z)}$ angle magnitude: }\hspace{2mm}%
\parbox{43mm}{

\vspace{.65mm}
\includegraphics[width=43mm,keepaspectratio]{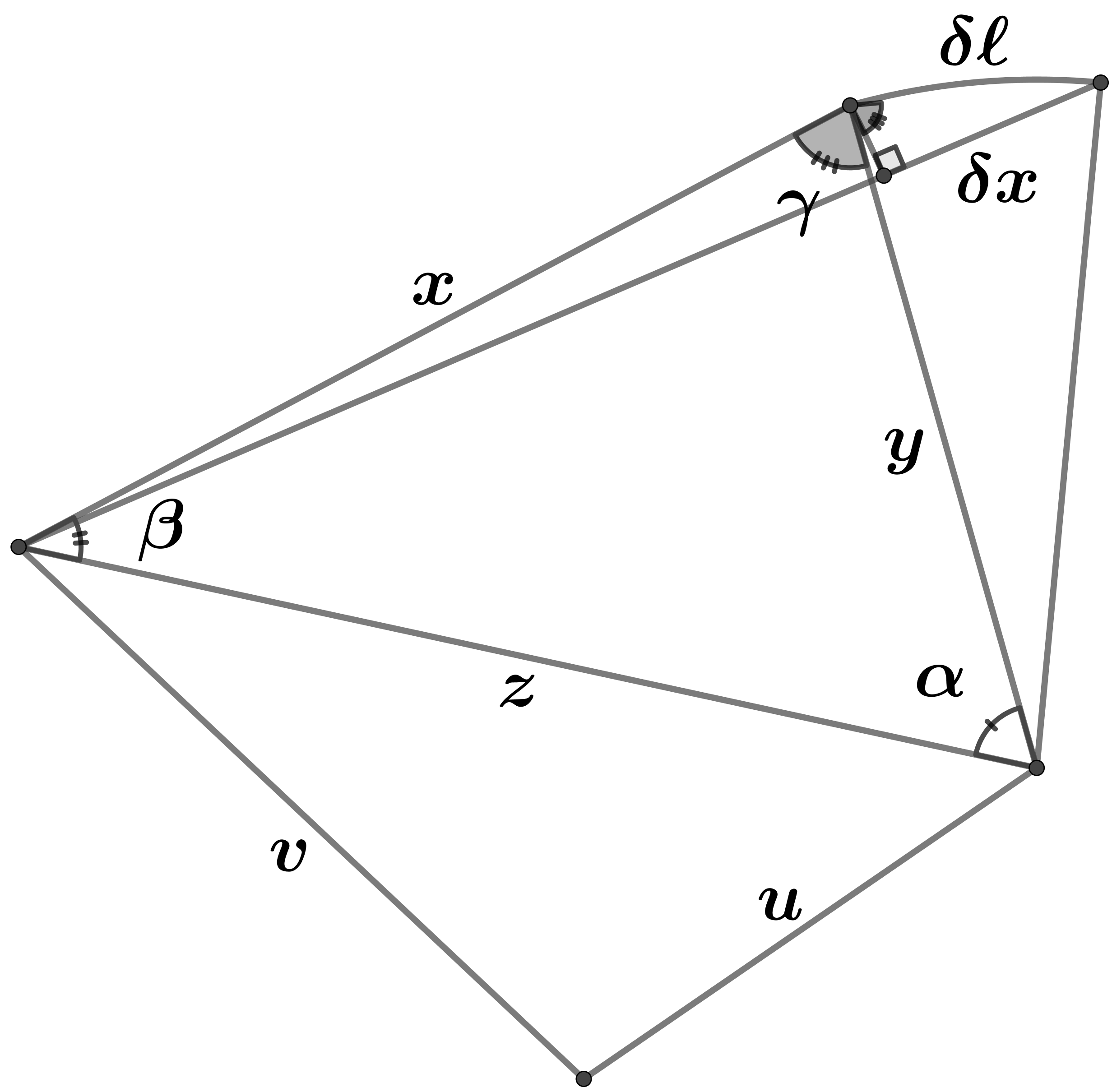}

\includegraphics[width=43mm,keepaspectratio]{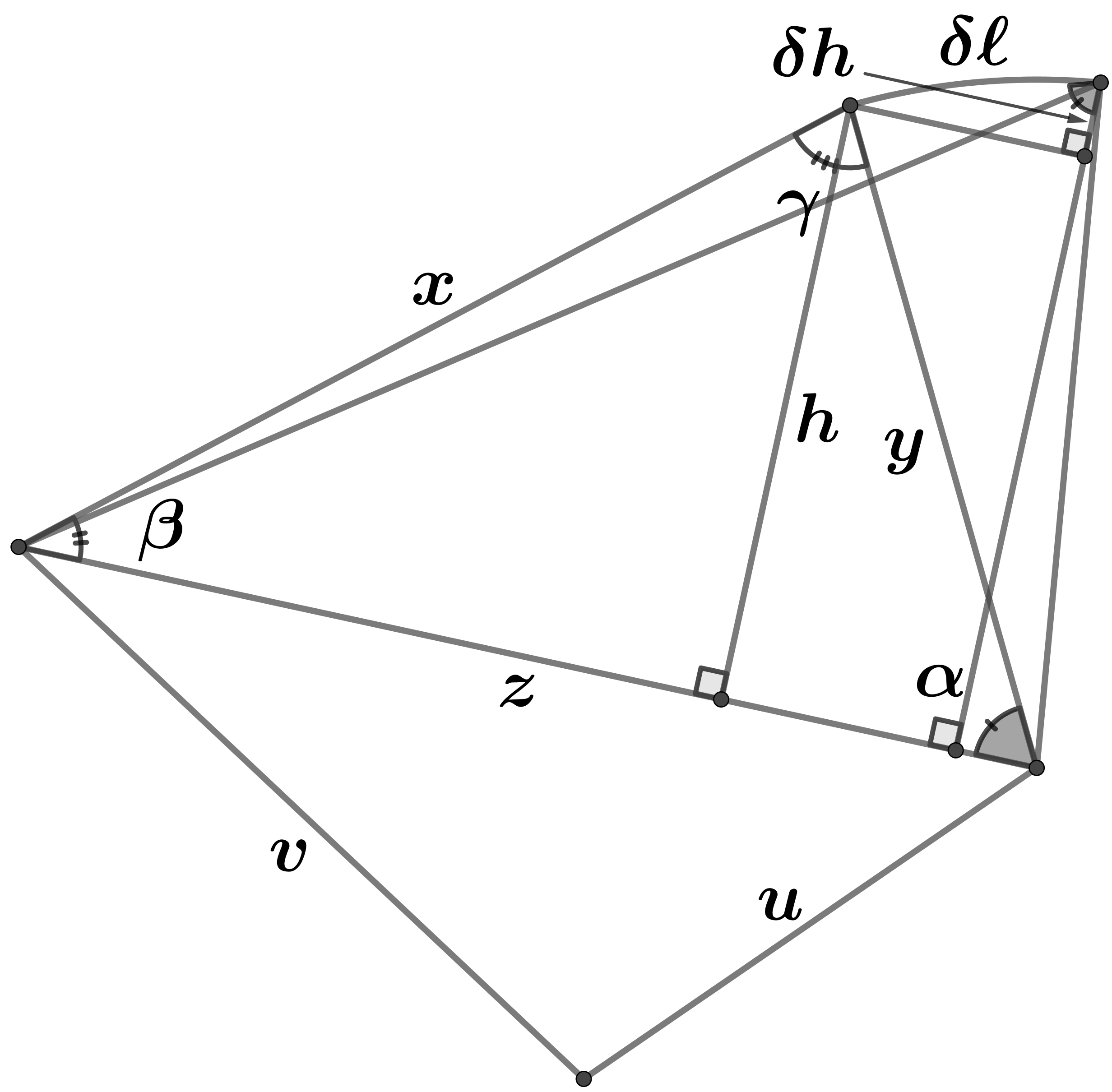}
\vspace{.65mm}
}%
}%

\begin{equation}
\delta x =\sin\gamma\,\delta \ell 
\qquad \delta h=\cos\alpha\,\delta \ell
\label{DeviationsBrahma0}
\end{equation}

\noindent where $\delta \ell$ is the distance travelled by the moved vertex. Plugging in these into eq. (\ref{PartialBrahma0}), using eq. (\ref{SineBrahma}) for the sine and al-Kashi theorem for the cosine, we find
\begin{equation}
A\partial_x A=\frac{xy+uv}{8y}[-x^2+y^2+z^2]
\label{PartialBrahma1}
\end{equation}

\noindent that is, applying eq. (\ref{PtolemyTheorem}) to express $z^2$ in terms of the sidelengths, 
\begin{equation}
A\partial_x A=\frac{1}{8}[-x^3+(y^2+u^2+v^2)x+2yuv]
\label{PartialBrahma2}
\end{equation}

\noindent which can be integrated out to give the general solution
\begin{align}
A^2(x,y,&\,u,v)\\
=\;&\!\frac{1}{16}[-x^4+2(y^2+u^2+v^2)x^2+8xyuv+k(y,u,v)]
\nonumber
\end{align}

\noindent where $k(y,u,v)$ is a homogeneous function of $y$, $u$ and $v$. $A(x,y,u,v)$ must be symmetric in $x$, $y$, $u$ and $v$. Consequently $k(y,u,v)=-y^4-u^4-v^4+2(y^2u^2+y^2v^2+u^2v^2)$. After factorization, we come across the Brahmagupta theorem \cite[\S 21]{Brahma:628}
\begin{align}
\label{BrahmaTheorem}
&A(x,y,u,v)\\
&=\sqrt{\frac{-x\!+\!y\!+\!u\!+\!v}{2}\;\frac{x\!-\!y\!+\!u\!+\!v}{2}\;\frac{x\!+\!y\!-\!u\!+\!v}{2}\;\frac{x\!+\!y\!+\!u\!-\!v}{2}}
\nonumber
\end{align}

\section{Leonhard Euler}
\label{Euler}

\noindent This one has a simple calculus part, but is a bit more trickier for the geometric principles. Imagine that Newton, or even Euler, would have been interested in determining the distance $d$ between the circumcenter and incenter of any triangle. Let him assume that $d$ depends on the inradius $r$ and the circumradius $R$, i.e.
\begin{equation}
d=d(r,R)
\label{EulerFunction_0}
\end{equation}

\vspace{2.5mm}
\noindent\parbox{86mm}{
\includegraphics[width=43mm,keepaspectratio]{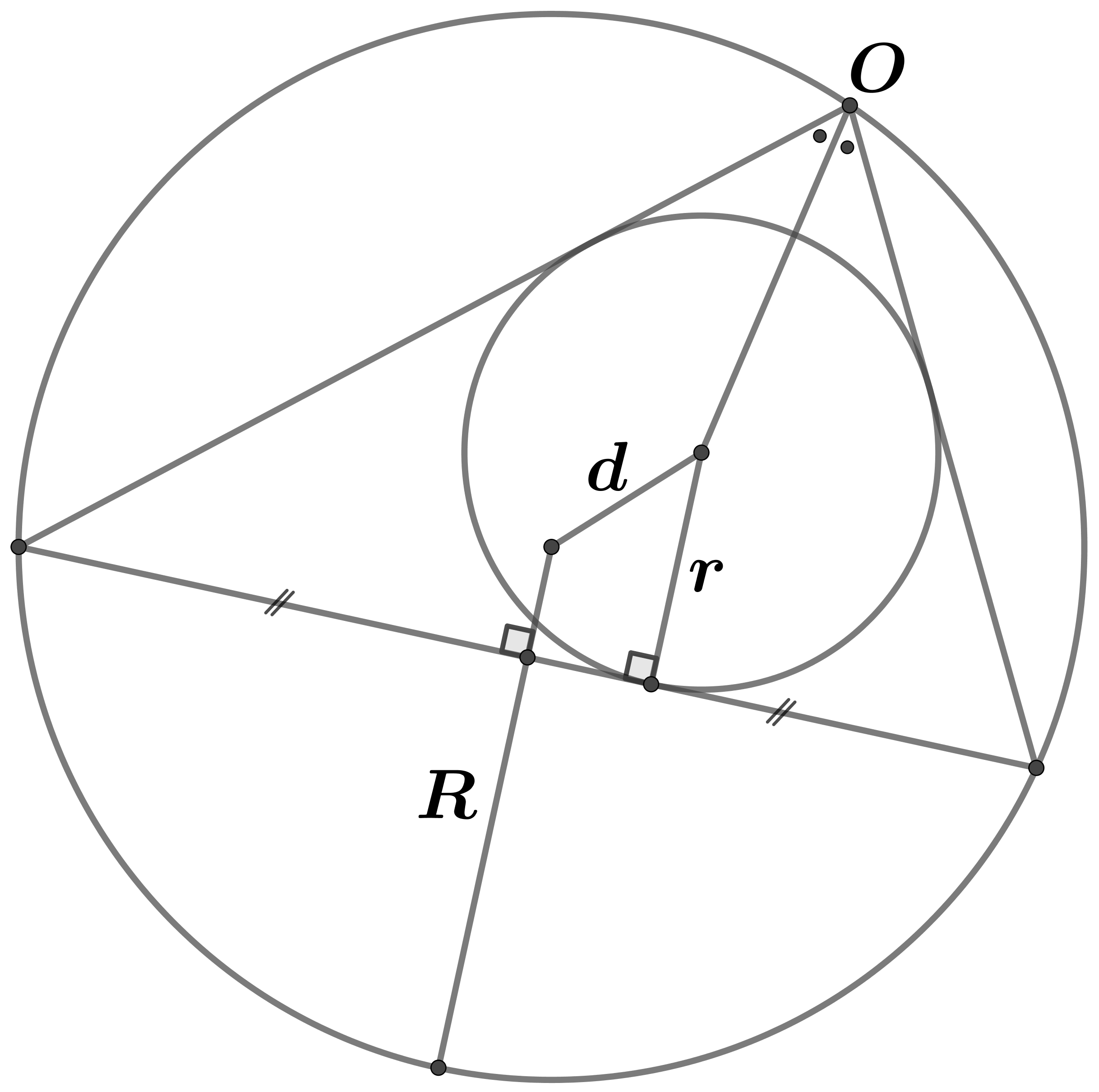}%
\includegraphics[width=43mm,keepaspectratio]{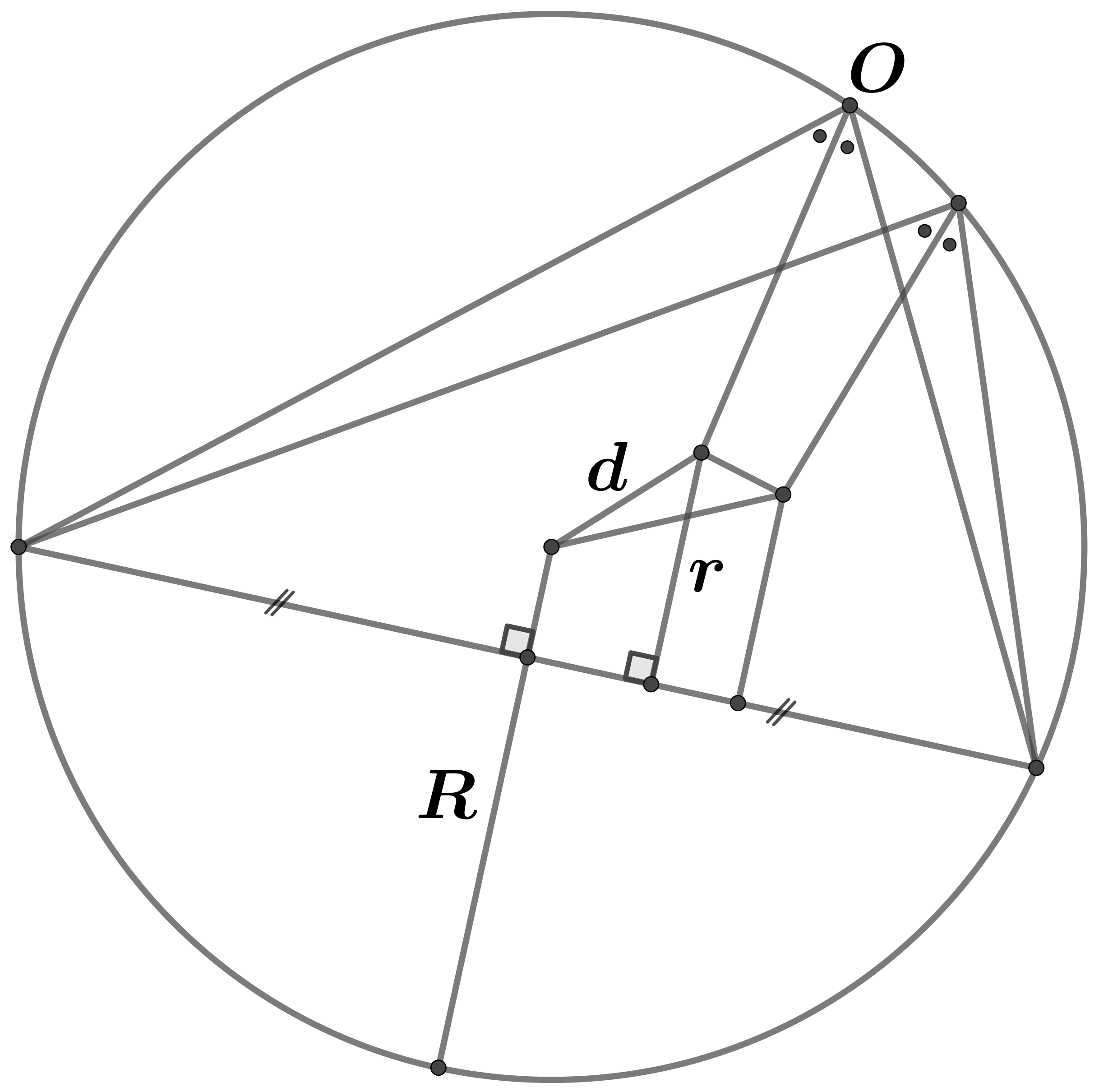}%
}%
\vspace{2.5mm}

\noindent After a slight displacement of one of the vertices, say $O$, along the circumcircle, resulting in a second triangle with small deviations $\delta d$ and $\delta r$, while $\delta R=0$, we at first order
\begin{equation}
\delta d=\partial_r d\,\delta r
\label{EulerDerivative_0}
\end{equation}

\noindent Newton could not ignore that the circumcenter is the intersection of the perpendicular side bisectors and that the incenter is the intersection of the internal angle bisectors. Some geometrical considerations need to be taken into account in order to proceed. 

\vfill
\vspace{2.5mm}
\noindent\parbox{86mm}{
\includegraphics[width=43mm,keepaspectratio]{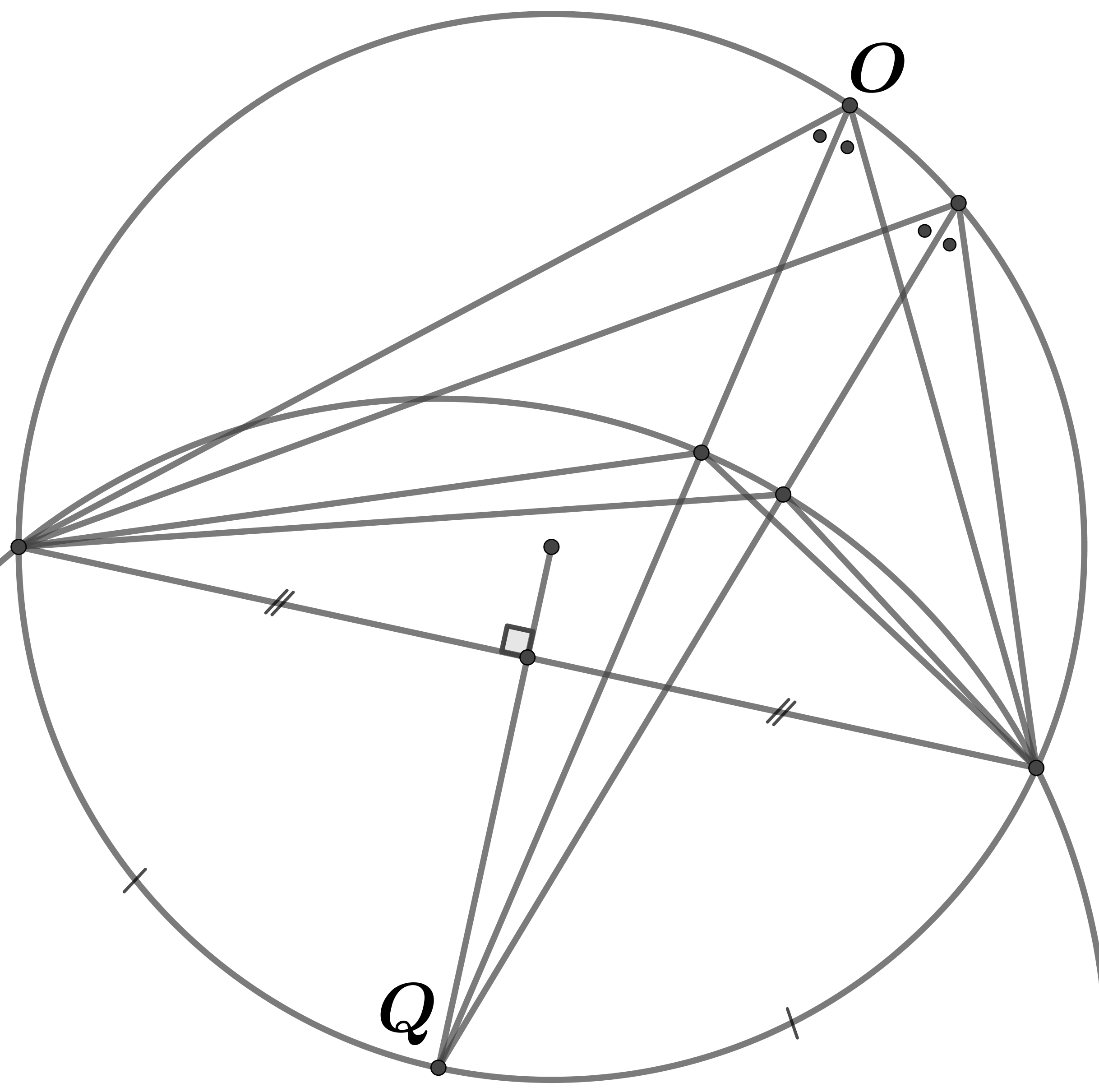}%
\includegraphics[width=43mm,keepaspectratio]{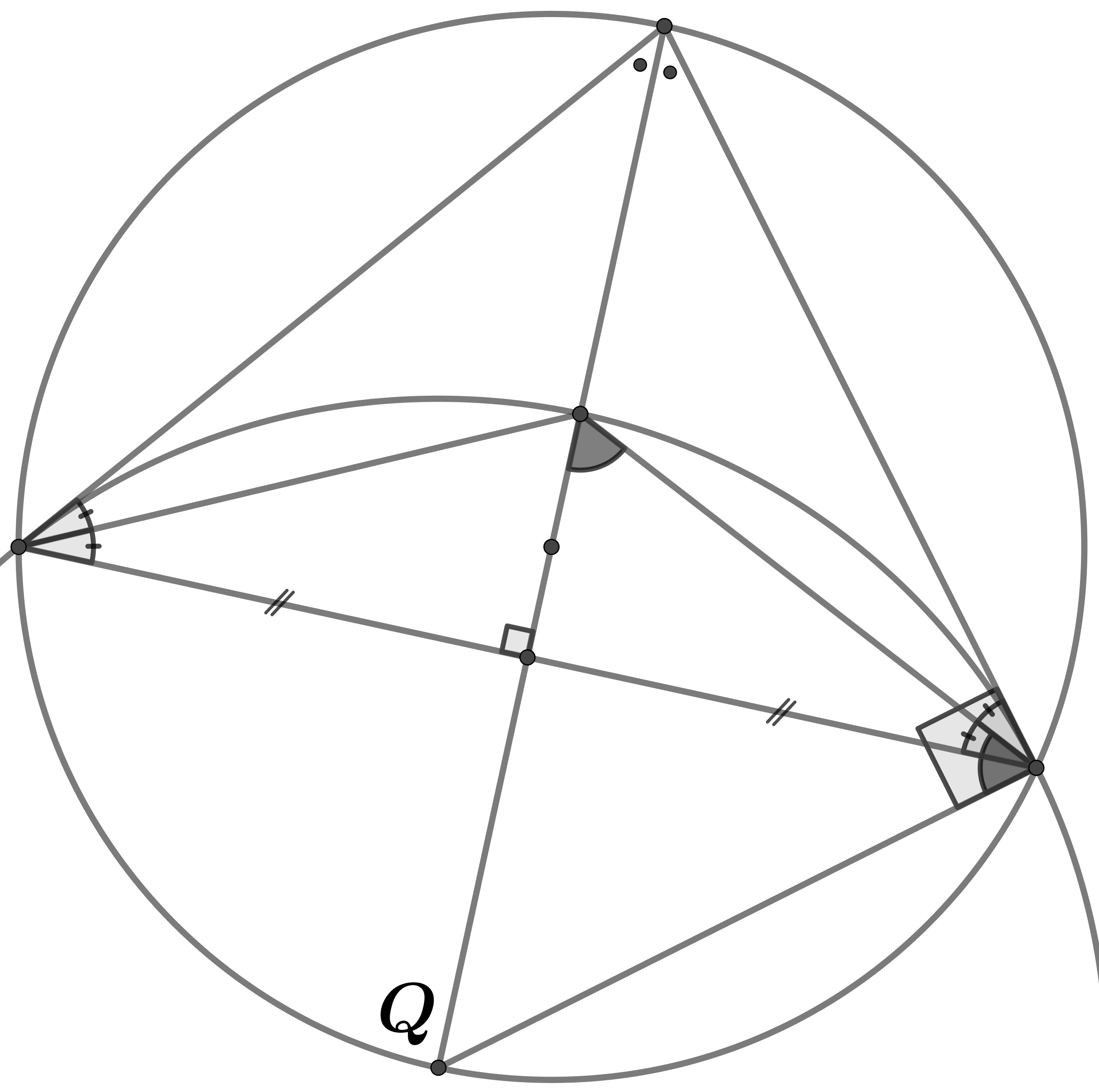}%
}

\newpage

First of all, the perpendicular side bisectors and the internal angle bisectors intersect on the circumcircle. To be convinced, just notice that the angles of the triangle are inscribed angles in the circumcircle. According to the inscribed angle theorem, their magnitude is determined by the intercepted arc. Hence, any angle bisector bisects the intercepted arc. So does the corresponding perpendicular side bisector. Let us call $Q$ the intersection of the angle bisector from $O$ and the perpendicular bisector of its opposite side. Obviously, this angle bisector still passes through $Q$ after the $O$ vertex has been moved along the circle. 

Secondly, when a vertex is moved along the circumcircle, the incenter moves along another circle whose center is precisely the intersection of the angle bisector from the moved vertex, and the perpendicular bisector of the opposite side, that is, $Q$. 

It is simple to show that the trajectory of the incenter is circular: since the angle of the moved vertex remains constant throughout the movement, so does the angle of the moved vertex of a new triangle formed by the two remaining fixed vertices and the incenter; the inscribed angle theorem again tells us that the incenter moves along a circle passing through the two fixed vertices. 

To see that the center of this circle is on the circumcircle -- and corresponds to $Q$ -- just move the incenter along its circular trajectory until it lays on the perpendicular bisector; in this symmetric configuration, a simple angle chase teaches us that the distance between the incenter and $Q$ is the same as between each one of the fixed vertices and $Q$, so that this distance must be the radius of the circle, and $Q$ its center. 

As a consequence, since it moves along a circle whose center is on the angle bisector from the moved vertex, at first order the incenter necessarily moves perpendicularly to this angle bisector. 

We can now go back to equation (\ref{EulerDerivative_0}) and try to find an expression for $\delta d$ and $\delta r$ by observing the way the initial triangle slightly moves to the second one. 

\vfill
\vspace{2.5mm}
\noindent\parbox{86mm}{
\includegraphics[width=86mm,keepaspectratio]{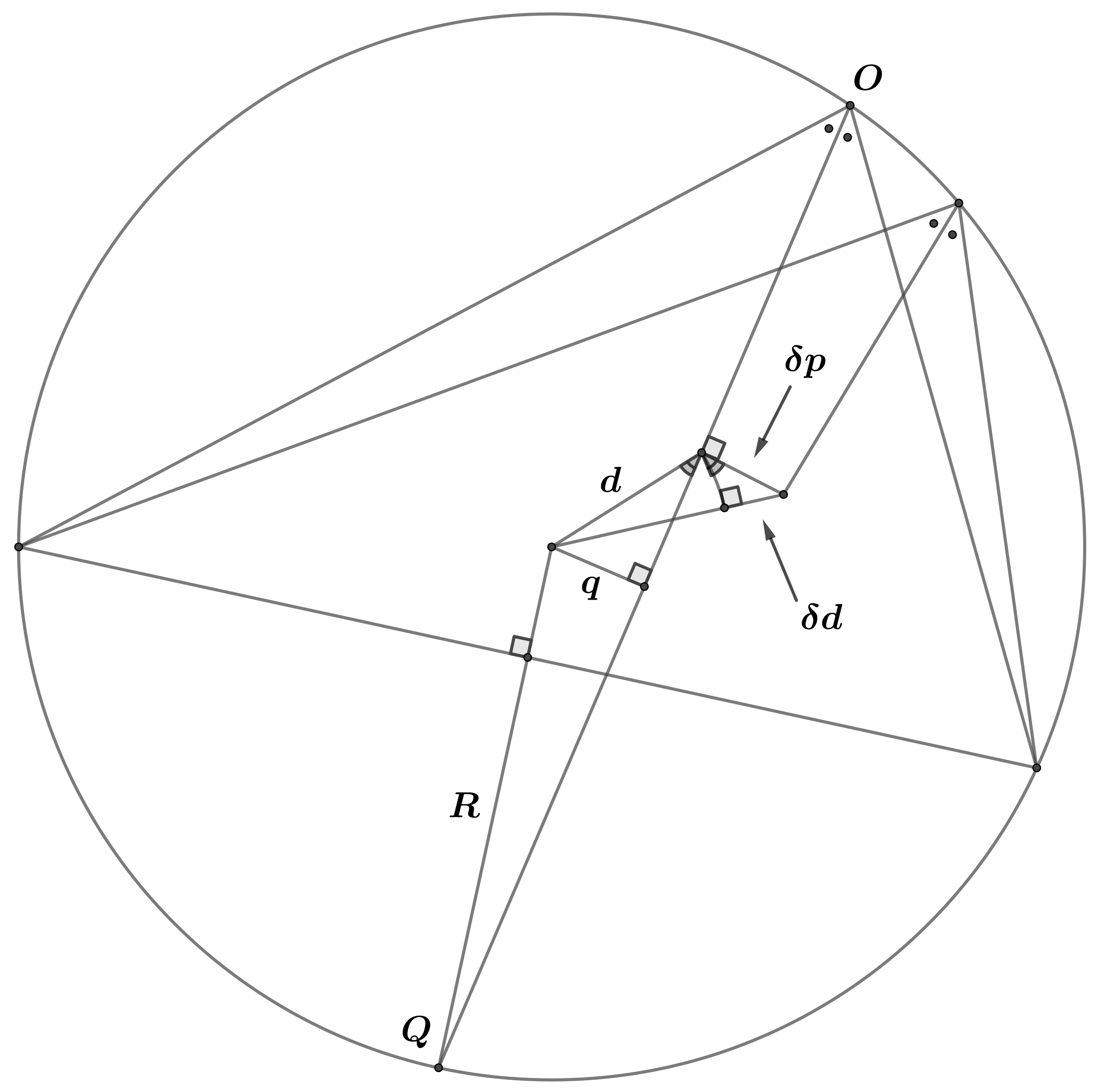}%
}

\newpage

Note that $\delta d$ (positive in the illustrative figure) is one of the legs of a right triangle whose hypotenuse length is the distance, say $\delta p$, between the incenters of the two triangles. Since at first order, the $\delta p$-length segment is perpendicular to the angle bisector from the moved vertex, and that the other leg is perpendicular to $d$-length segment, this right triangle is similar to another one of $d$-length hypotenuse and useful leg of length, say, $q$. Hence
\begin{equation}
\delta d=\frac{q}{d}\delta p
\label{EulerDeviations_0}
\end{equation} 

\vspace{2.5mm}
\noindent\parbox{86mm}{
\includegraphics[width=86mm,keepaspectratio]{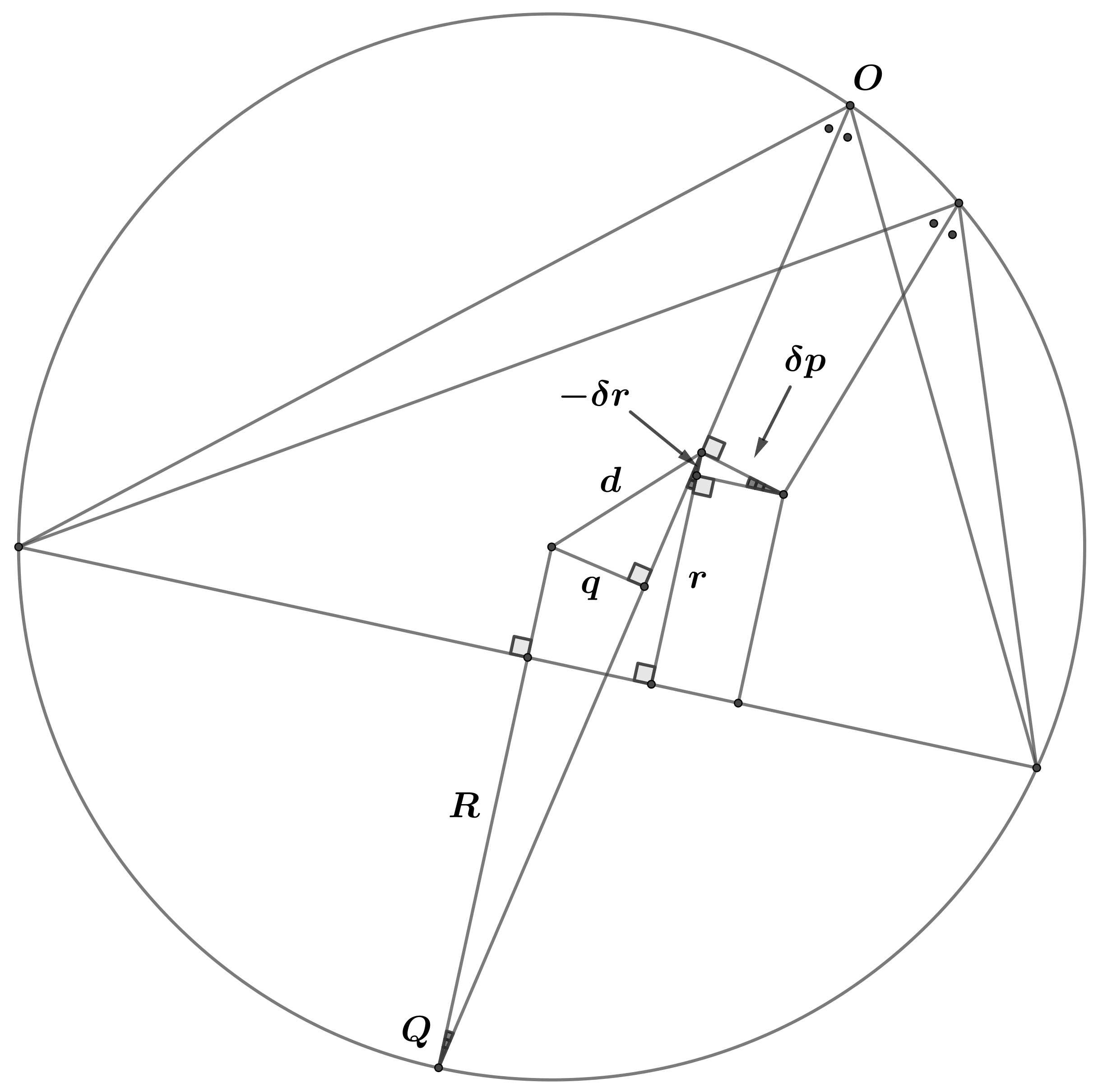}%
}%
\vspace{2.5mm}

\noindent For $\delta r$ (negative in the illustrative figure), we need to draw the inradius perpendicularly to the side opposite to the moved vertex, and observe the way it evolves. Here, it is $-\delta r$ that can be seen as the length of one of the legs of a right triangle with $\delta p$-length hypotenuse, which is similar to another one of hypotenuse of length $R$ and corresponding leg of length $q$. We have
\begin{equation}
\delta r=-\frac{q}{R}\delta p
\label{EulerDeviations_1}
\end{equation}

\noindent Dividing equations (\ref{EulerDeviations_0}) and (\ref{EulerDeviations_1}) member to member and using equation (\ref{EulerDerivative_0}) yields the differential equation
\begin{equation}
d\partial_r d=-R
\label{EulerEquaDif}
\end{equation}

\noindent which we integrate to obtain the general solution
\begin{equation}
d^2(r,R)=-2Rr+k(R)
\label{Euler_1}
\end{equation}

\noindent with $k(R)\ge 2Rr$ a real function of $R$. Since in an equilateral triangle, $d^2(R/2,R)=0$, we find $k(R)=R^2$. Hence 
\vspace{-2mm}
\begin{equation}
d^2(r,R)=R(R-2r)
\label{Euler_2}
\end{equation}

\noindent that is, Euler's theorem \cite{Euler:1767}, first proved by Chapple \cite{Chapple:1746}.

\section{The Angle Bisector} 
\label{BisPart}

\noindent We already have a section (\ref{Terquem}) dedicated to the angle bisector. But it was about its full length: from the vertex to the foot. What if Newton got now interested to the part between the vertex and the incenter? It should also be a smooth function of the sidelength $x$, $y$ and $z$, i.e.

\noindent\parbox{86mm}{\parbox{43mm}{

\includegraphics[width=43mm,keepaspectratio]{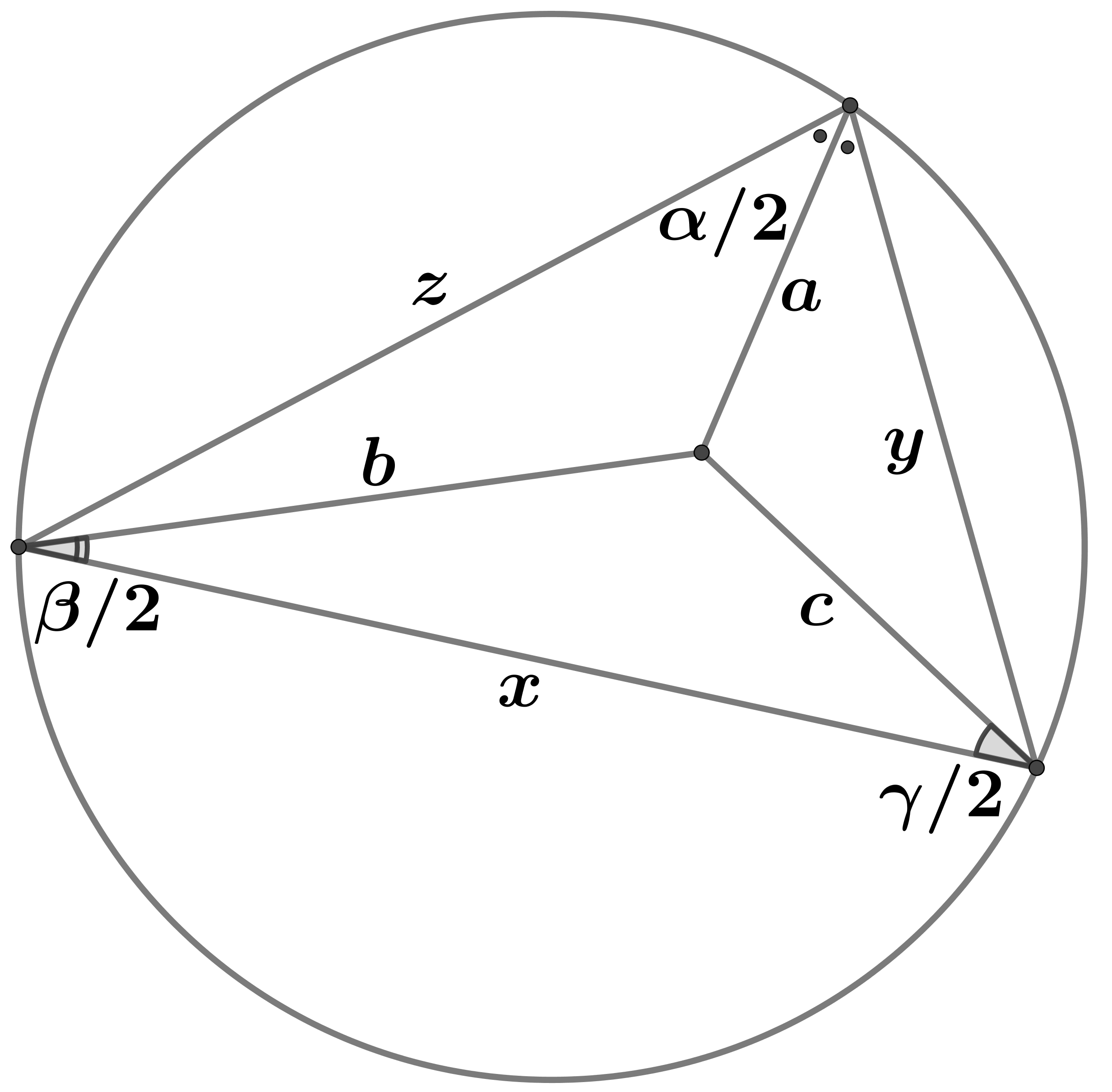}

\includegraphics[width=43mm,keepaspectratio]{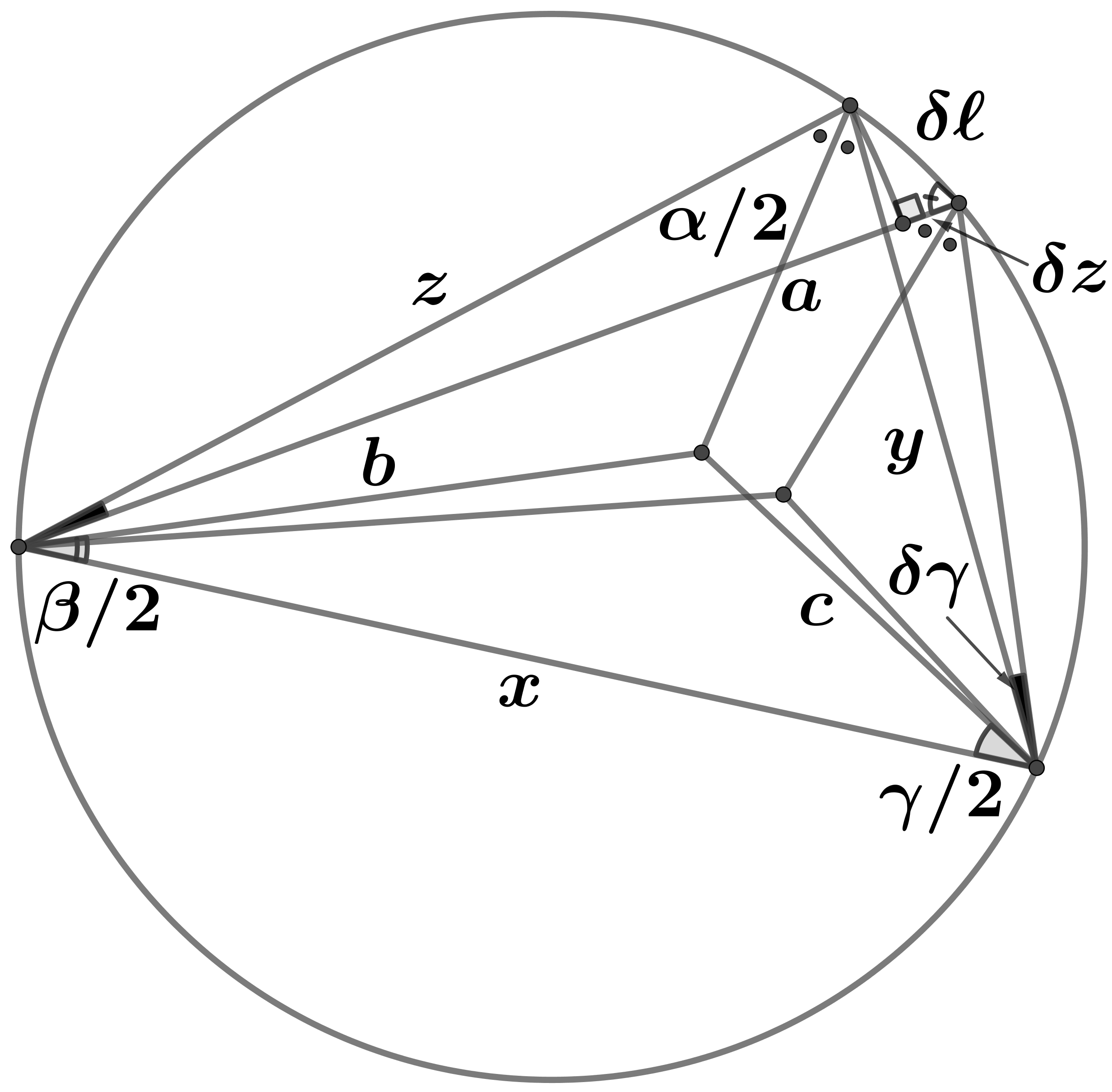}

}\hspace{2mm}%
\parbox{41mm}{
\begin{equation}
c=c(x,y,z)
\end{equation}

\noindent Since the position of the incenter matters, and that we know how it moves in the circumcircle, we will let Newton continue his considerations in this very circumcircle. After a small clockwise displacement of the $(z,y)$ vertex along the circumcircle, we have infinitesimal deviations $\delta y$ and $\delta z$, while $\delta x=0$ (in the illustrative figures, $x$ and $z$ have been swapped for aesthetic reasons). At first order,
\begin{equation}
\delta c= \partial_y c \,\delta y + \partial_z c \,\delta z
\label{PartialBisectorIn0}
\end{equation}
}
}

\noindent In the small right triangle with legs of length $-\delta c$ and $c\delta\gamma/2$ -- where $\gamma=\widehat{(x,y)}$ -- the angle opposite to the $c\delta\gamma/2$-length side has the same magnitude as the one between the $x$-length side and the $\beta=\widehat{(x,z)}$ angle bisector, that is $\beta/2$, since they both are inscribed angles intercepting the same arc of the circle along which the incenter travels. Thus we have 
\begin{equation}
\tan\frac{\beta}{2}=\frac{{c\delta\gamma}/{2}}{-\delta c}
\end{equation}

\vfill
\noindent\parbox{86mm}{
\includegraphics[width=86mm,keepaspectratio]{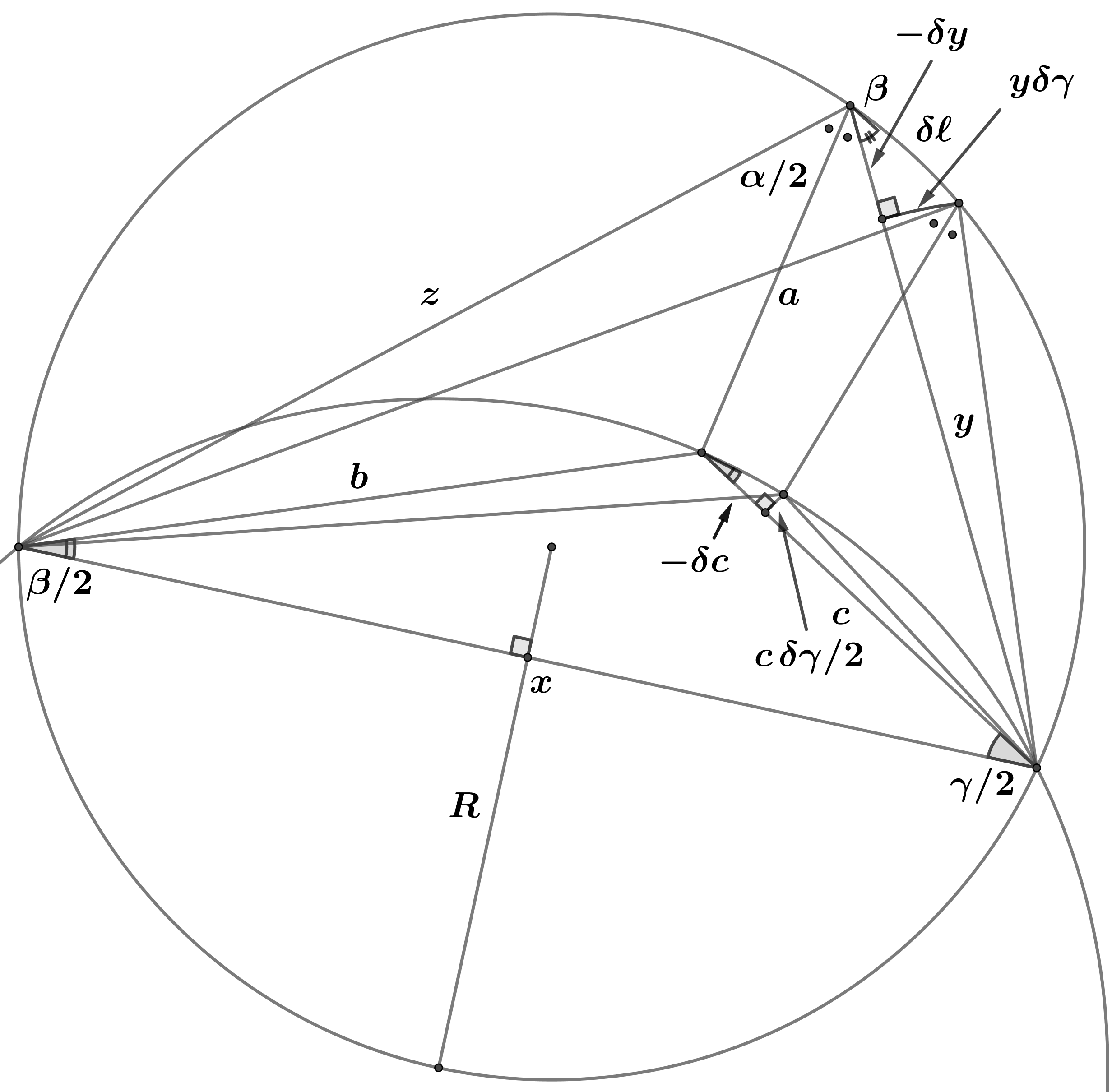}%
}

\newpage

\noindent Furthermore, based on the same reasoning as in the last section, we can connect $-\delta y$ and $\delta z$ to $\delta\ell$, as well as $y\delta\gamma=\sin\beta\delta\ell$, so that using eq. (\ref{tanhalfangle}) for the tangent, we find
\begin{equation}
-\delta y\!=\!\cos\beta\,\delta \ell\quad\!\delta z\!=\!\cos\gamma \delta\ell\quad\! \delta c\!=\!\frac{-c(1+\cos\beta)}{2y}\delta\ell
\label{DeviationsBisectorIn}
\end{equation}

\noindent Plugging in results (\ref{DeviationsBisectorIn}) into equation (\ref{PartialBisectorIn0}) gives a partial differential equation whose analogue can be obtained by considering a clockwise rotation of the $(x,z)$ vertex instead of the $(z,y)$ one. Both equations read
\begin{align}
\frac{c}{2y}(1+\cos\beta)&=\cos\beta\,\partial_y c -\cos\gamma\,\partial_z c\\
\frac{c}{2x}(1+\cos\alpha)&=\cos\alpha\,\partial_x c -\cos\gamma\,\partial_z c
\label{PartialBisectorIn1}
\end{align}

\noindent Combining those with the differential equation arising from the infinitesimal scale transformation $R\mapsto R+\delta R$, that is
\begin{equation}
c=x\,\partial_x c + y\,\partial_y c + z\,\partial_z c
\label{ScaleBisectorIn}
\end{equation}

\noindent and isolating $\partial_z c$, yields
\begin{equation}
\frac{\partial_z c}{c}=\frac{\cos\alpha + \cos\beta}{x\cos\beta\cos\gamma+y\cos\alpha\cos\gamma+z\cos\alpha\cos\beta}
\label{PartialBisectorIn2}
\end{equation}

\noindent which, using al-Kashi's theorem for the cosines and working things out a little bit, can be reexpressed as
\begin{equation}
\frac{\partial_z c}{c}=\frac{x+y}{(x+y+z)(x+y-z)}
\label{PartialBisectorIn3}
\end{equation}

\noindent and integrated out to give the general solution
\begin{equation}
c(x,y,z)
=k(x,y)\sqrt{\frac{x+y-z}{x+y+z}}
\label{GeneralBisectorIn}
\end{equation}

\noindent where $k(x,y)$ is a symmetric function of $x$ and $y$. In the particular case of a right triangle with hypotenuse of length $z$, $c=\sqrt{2}r$ where $r$ is the inradius length. Computing the area of this right triangle yields $r(x+y+z)=xy$. Substituting $r$ in the expression of $c$ and inserting the latter in eq. (\ref{GeneralBisectorIn}), we find $k(x,y)=\sqrt{xy}$. Hence
\begin{equation}
c(x,y,z)
=\sqrt{xy\,\frac{x+y-z}{x+y+z}}
\label{BisPartTheorem}
\end{equation}

\noindent that is, the distance between the incenter and the foot of the $\gamma=\widehat{(x,y)}$ angle bisector. Dividing $c(x,y,z)$ in this equation by the full length of the bisector $d(x,y,z)$ in eq. (\ref{TerquemTheorem}) gives the ratio $(x+y)/(x+y+z)$, meaning that the relative position/height of the incenter on the angle bisector (measured from its foot) is equal to the ratio of its corresponding sidelength and the perimeter of the triangle, that is $z/(x+y+z)$; a result that can easily be verified in barycentric coordinates. 

\section{The Inradius}
\label{Inradius}

\noindent If we were able to do it for the circumradius, section \ref{Circumradius}, we should not doubt of Newton's appetancy to do it for the inradius $r$! It would also be a smooth function of the sidelength $x$, $y$ and $z$, that is
\begin{equation}
r=r(x,y,z)
\end{equation}

\vspace{2.5mm}
\noindent\parbox{86mm}{
\includegraphics[width=43mm,keepaspectratio]{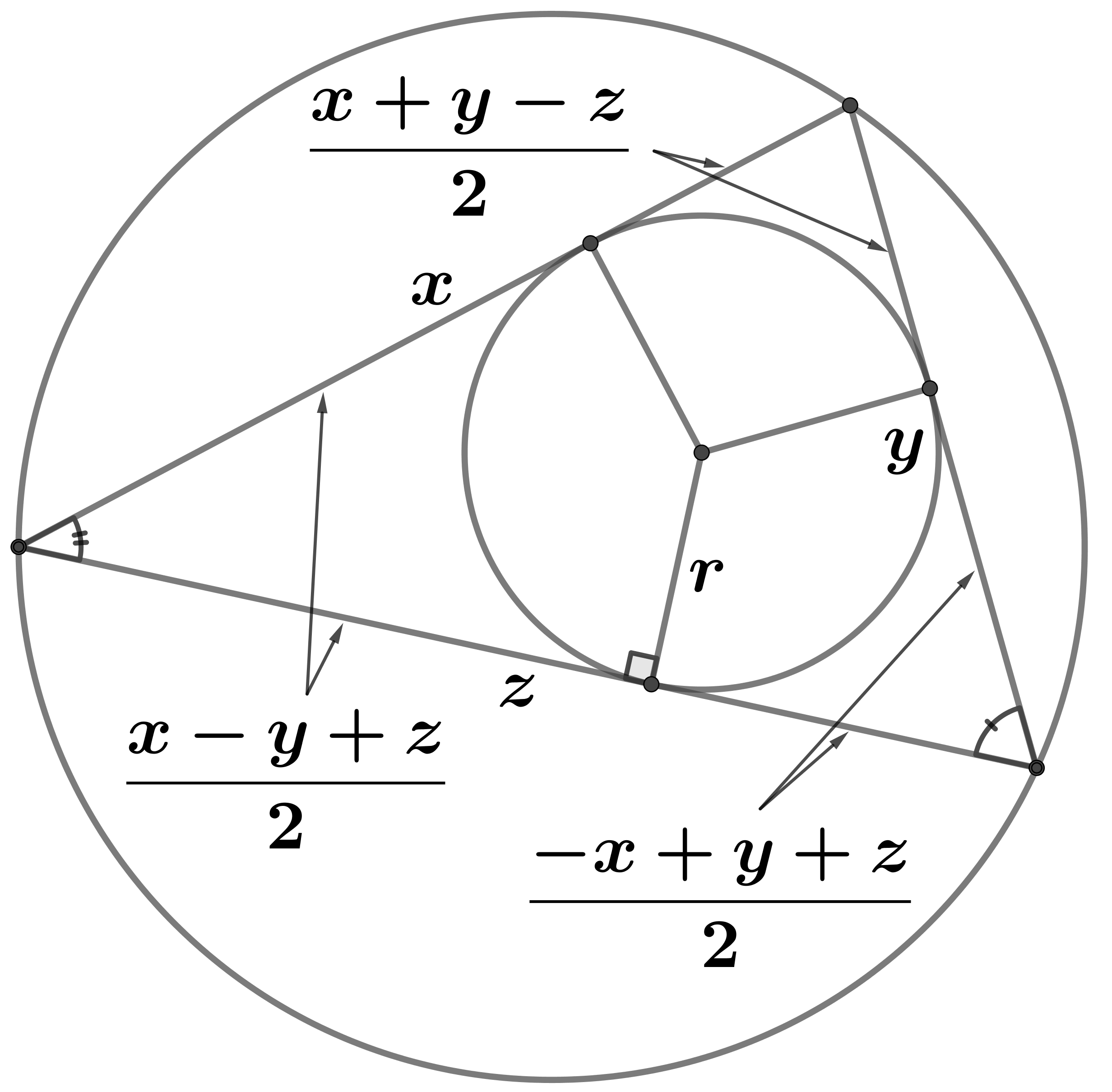}%
\includegraphics[width=43mm,keepaspectratio]{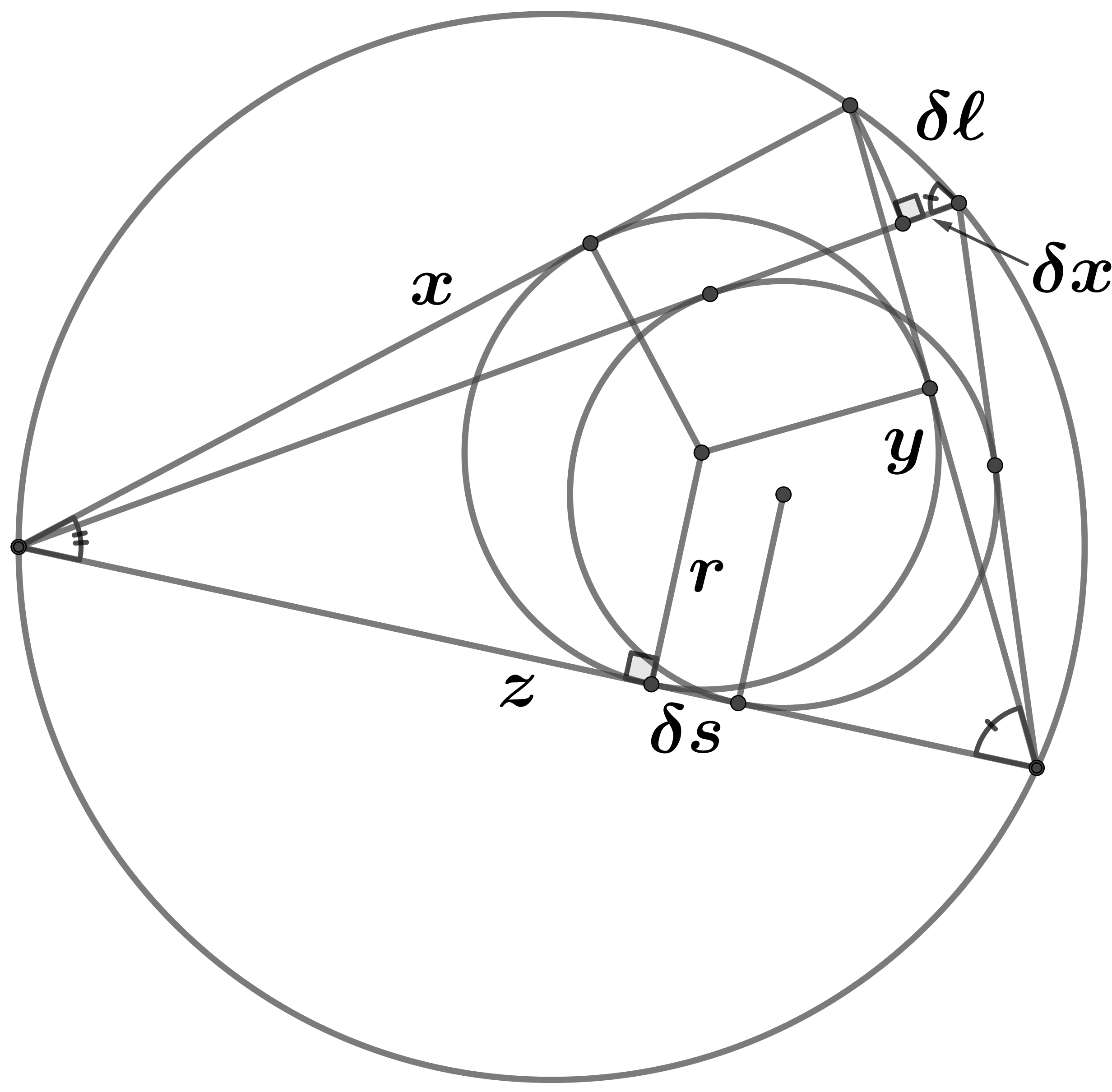}%
}%
\vspace{2.5mm}

\noindent Strange as it may seem, since the position of the incenter matters, the circumcircle looks again to be the best place to do the job for the inradius. After a small displacement of the $(x,y)$ vertex along the circumcircle, we have infinitesimal deviations $\delta x$ and $\delta y$, while $\delta z=0$. At first order,
\begin{equation}
\delta r= \partial_x r \,\delta x + \partial_y r \,\delta y
\label{PartialInradius0}
\end{equation}

\noindent $\delta x$ and $\delta y$ are determined in the same way as in the last sections -- see for instance eq. (\ref{DeviationsCircum0}) for the case of the circumradius:
\begin{equation}
\delta x =\cos\alpha\,\delta \ell 
\qquad -\delta y=\cos\beta\,\delta \ell
\label{DeviationsInradius0}
\end{equation}

\vspace{2.5mm}
\noindent\parbox{86mm}{
\includegraphics[width=86mm,keepaspectratio]{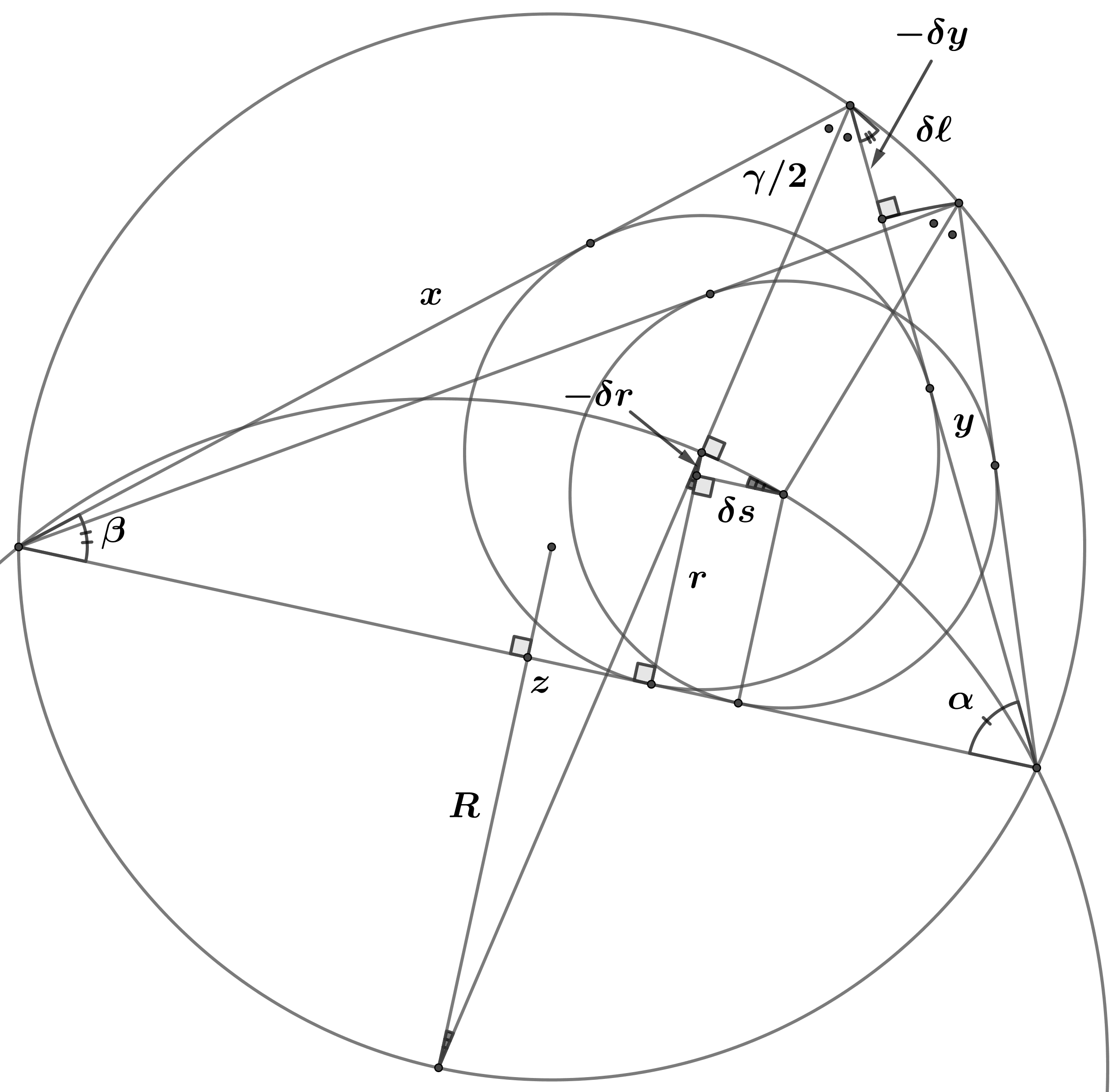}%
}%
\vspace{2.5mm}

\noindent For $\delta r$ we will refer to the geometrical considerations of section \ref{Euler} and consider the small right triangle of legs of lengths $-\delta r$ and $\delta s$ where $\delta s$ is the distance travelled by the foot of the inradius perpendicular to the $z$-length side. First of all, we check that this inradius foot divides the $z$-length side in two segments of lengths $(x-y+z)/2$ and $(-x+y+z)/2$ respectively, so that after the infinitesimal displacement it is moved of $\delta s=(\delta x -\delta y)/2$. Secondly, in this small right triangle, the angle opposite to the $-\delta r$-length leg is equal to the angle between the internal $\gamma$ angle bisector and the height from the same vertex (or any line perpendicular to the opposite $z$-length side), that is $\gamma/2-(\pi/2-\alpha)=(\alpha-\beta)/2$ since $\gamma=\pi-(\alpha+\beta)$. We have
\begin{align}
-\delta r &= \tan\frac{\alpha-\beta}{2}\,\frac{\delta x -\delta y}{2}
\nonumber\\
&=\frac{\sin\alpha-\sin\beta}{2}\,\delta\ell
\label{DeviationsInradius1}
\end{align}

\noindent thanks to eq. (\ref{DeviationsInradius0}) and the development of the tangent
\begin{equation}
\tan\frac{\alpha-\beta}{2}=\frac{\sin\alpha-\sin\beta}{\cos\alpha+\cos\beta}
\end{equation} 

\noindent Plugging in results (\ref{DeviationsInradius0}) and (\ref{DeviationsInradius1}) into equation (\ref{PartialInradius0}) gives a partial differential equation whose analogue can be obtained by considering a rotation of the $(z,x)$ vertex instead of the $(x,y)$ one. Both equations read
\begin{align}
\frac{\sin\beta-\sin\alpha}{2}&=\cos\alpha\partial_x r -\cos\beta\partial_y r\\
\frac{\sin\alpha-\sin\gamma}{2}&=\cos\gamma\partial_z r -\cos\alpha\partial_x r
\label{PartialInradius1}
\end{align}

\noindent Combining those with the differential equation arising from the infinitesimal scale transformation $R\mapsto R+\delta R$, that is
\begin{equation}
r=x\,\partial_x r + y\,\partial_y r + z\,\partial_z r
\label{ScaleInradius}
\end{equation}

\noindent and isolating $\partial_x r$, yields
\begin{align}
\hspace{-2.9mm}\partial_x r&\!-\!\frac{\cos\beta \cos\gamma}{x\cos\beta\cos\gamma+y\cos\alpha\cos\gamma+z\cos\alpha\cos\beta}\,r
\nonumber\\
&\!=\!\frac{y(\sin\beta-\sin\alpha)\cos\gamma+z(\sin\gamma-\sin\alpha)\cos\beta}{2(x\cos\beta\cos\gamma+y\cos\alpha\cos\gamma+z\cos\alpha\cos\beta)}\hspace{-1mm}
\label{PartialInradius2}
\end{align}

\noindent which, using the fundamental law for the sines, al-Kashi's theorem for the cosines and working things out a little bit, can be reexpressed as
\begin{align}
\hspace{-2.2mm}\partial_x r&\!-\!\frac{(x^2-y^2+z^2)(x^2+y^2-z^2)}{x(2(y^2 z^2 + x^2 z^2 + x^2 y^2)-x^4-y^4-z^4)}\,r
\nonumber\\
&\!=\!\frac{(y-x)(x^2+y^2-z^2)\!+\!(x-z)(x^2-y^2+z^2)}{2x\sqrt{2(y^2 z^2 + x^2 z^2 + x^2 y^2)-x^4-y^4-z^4}}\hspace{-1.7mm}
\label{PartialInradius3}
\end{align}

\noindent This equation is linear. Its general solution reads
\begin{equation*}
r(x,y,z)
=\frac{-x^3+(y+z)x^2-(y+z)(y-z)^2+x\,c(y,z)}{2\sqrt{2(y^2 z^2 + x^2 z^2 + x^2 y^2)-x^4-y^4-z^4}}
\label{GeneralInradius4}
\end{equation*}

\noindent where $c(y,z)$ is a homogeneous function of $y$ and $z$. Since $r(x,y,z)$ must be symmetric in $x$, $y$ and $z$, we have $c(y,z)=y^2+z^2+myz$ where $m$ is a real constant that can easily be determined by considering, for instance, an equilateral triangle with $x=y=z=2\sqrt{3}r$. We find $m=-2$. Hence, after factorization and simplification,
\begin{equation}
r(x,y,z)
=\sqrt{\frac{(-x+y+z)(x-y+z)(x+y-z)}{4(x+y+z)}}
\label{InradiusTheorem}
\end{equation}

\noindent that is, the expected expression of the inradius as the area (see section \ref{Heron}) divided by twice the perimeter.

\section{The Angle Bisectors Problem}
\label{BisProb}

\noindent There is a well-known problem about the internal angle bisectors of a triangle: given three real and positive numbers $a$, $b$, $c$, find the sidelength $x$, $y$, $z$ of the triangle that would admit those numbers as angle bisector lengths. A lot of things have been said about it. We would not be far from the truth if we say that there is a solution, but that it is the solution of a polynomial equation of too high a degree to be solved \cite{Terquem:1842}. See \cite{Altshiller:1953,Mironescu:1994,Dinca:2010} for a historical perspective and more recent attempt. We will just propulse Newton as a new challenger for a similar -- but not equivalent -- problem, that is identifying $a$, $b$, $c$ not as the full lengths of the angle bisectors, but of their part going from the vertex to the incenter. This problem has a solution, but as it is the solution of a third degree polynomial equation, we will not try to write it down. Let Newton state it this way: each of the sidelengths must be a smooth function of $a$, $b$, $c$, the length of the segments which bisect the angles $\alpha$, $\beta$ and $\gamma$ respectively, opposite to the sides of lengths $x$, $y$, $z$ respectively. For $z$, we have
\begin{equation}
z=z(a,b,c)
\end{equation}

\noindent\parbox{86mm}{
\includegraphics[width=43mm,keepaspectratio]{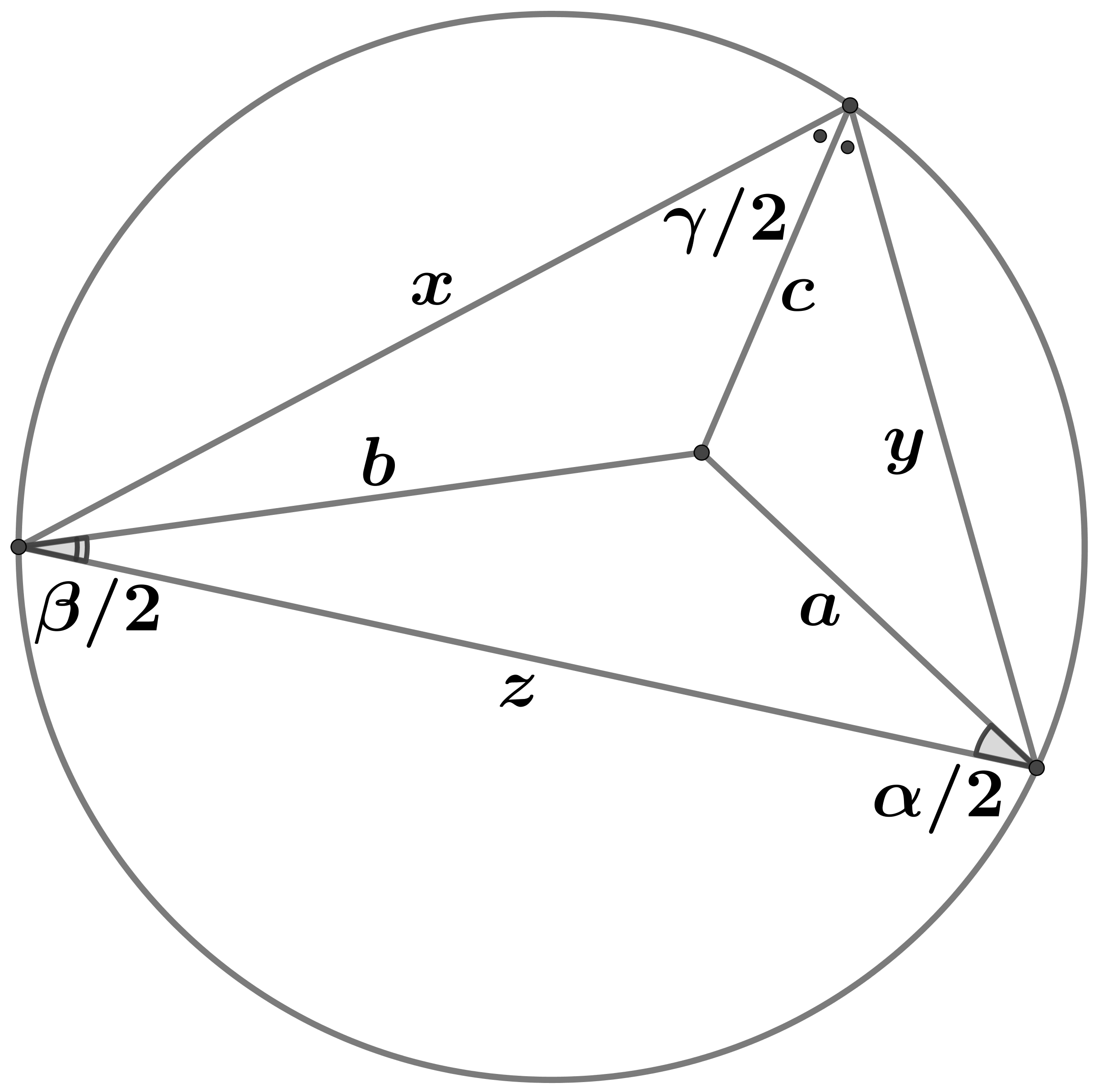}%
\includegraphics[width=43mm,keepaspectratio]{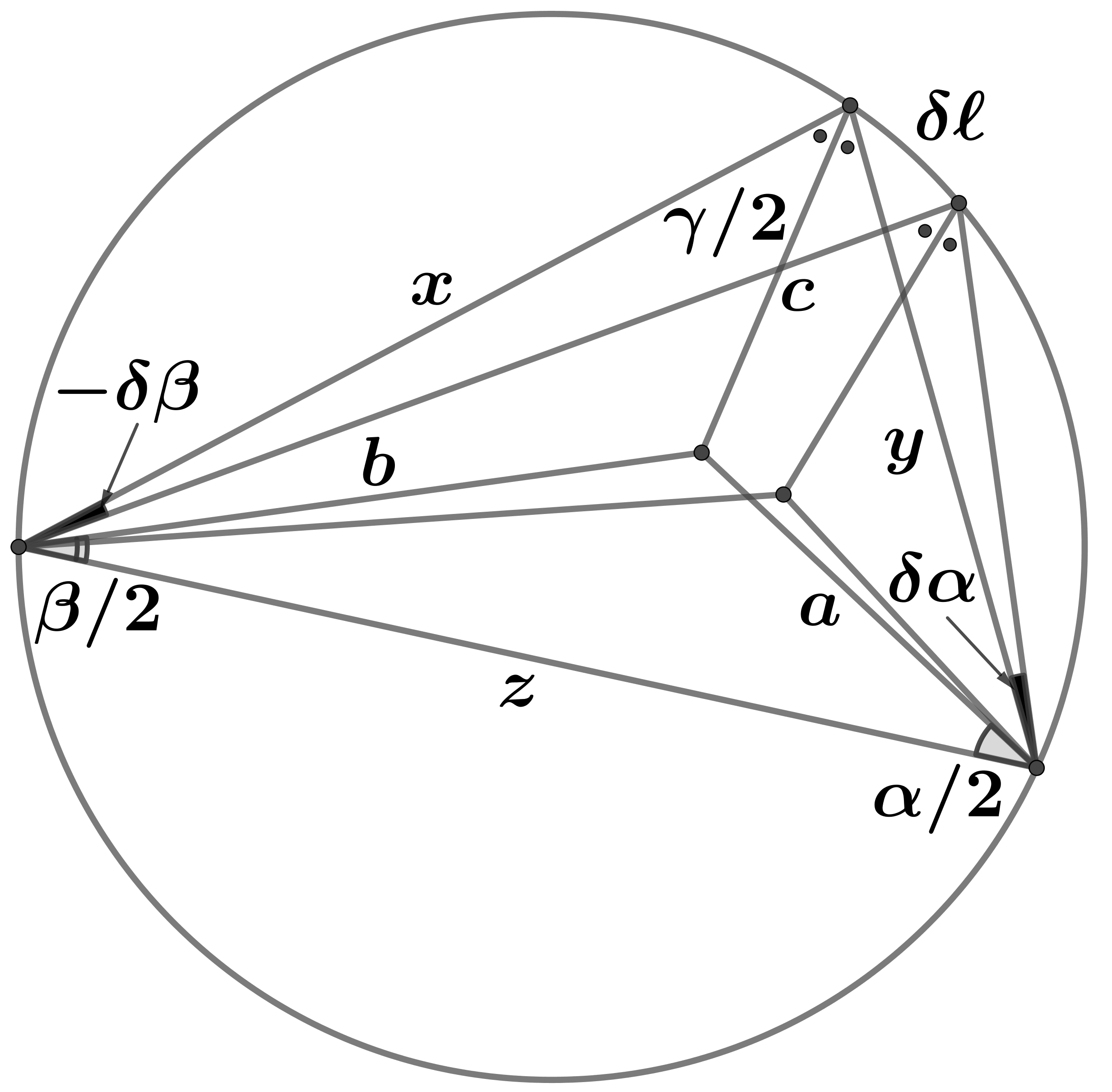}%
}%
\vspace{2.5mm}

\noindent Again, since the position and the movement of the incenter matter, the circumcircle is the best place to operate. After a small displacement of the $(x,y)$ vertex along the circumcircle, we have the infinitesimal deviations $\delta a$, $\delta b$ and $\delta c$, while $\delta z=0$. At first order,
\begin{equation}
\partial_a z \,\delta a + \partial_b z \,\delta b + \partial_c z \,\delta c = 0
\label{PartialBisectorProb0}
\end{equation}

\noindent If $\delta p$ is the distance between the incenters before and after the small displacement, the angle between this $\delta p$-length segment and the $-\delta a$-length is an inscribed angle of the circle along which the incenter travels, and has magnitude $\beta/2$. Analogously, the angle between the $\delta p$-length segment and the $\delta b$-length has magnitude $\alpha/2$. 

\noindent\parbox{86mm}{
\includegraphics[width=86mm,keepaspectratio]{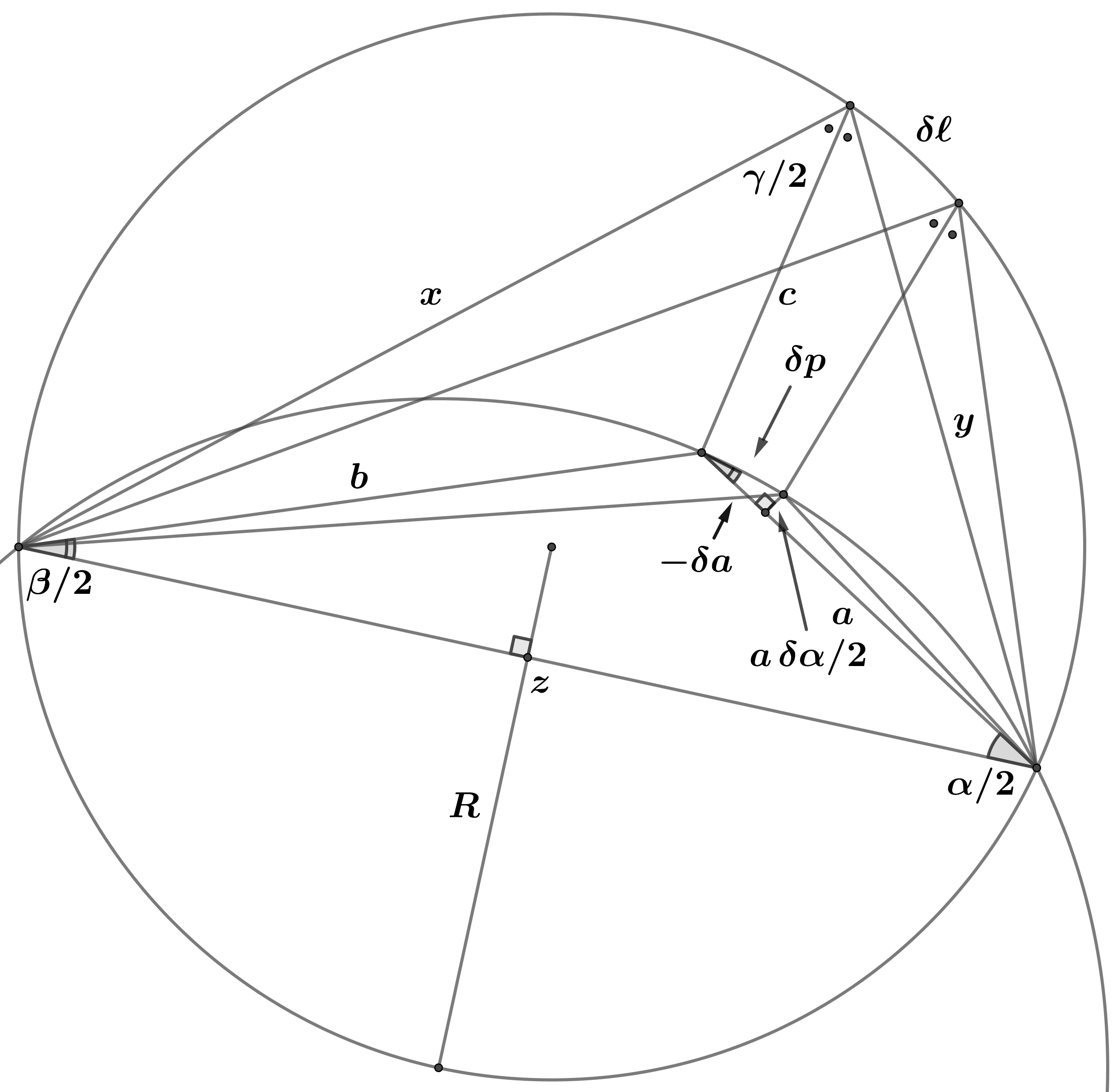}%
}%

\noindent Thus
\begin{equation}
\delta a =-\cos\frac{\beta}{2}\,\delta p
\qquad
\delta b =\cos\frac{\alpha}{2}\,\delta p
\label{DeviationsBisectorProb1}
\end{equation}

\noindent It is a bit more complicated for $\delta c$. We first observe that $-\delta c$ is the length of one of the legs of a small right triangle with hypotenuse of length $\delta \ell$, the distance travelled by the moved vertex along the cicrumcircle. The magnitude of the angle between those two small sides is given by $\gamma/2 + \beta$. Since $\gamma=\pi-(\alpha+\beta)$, the magnitude of the complemetary of this angle is $(\alpha-\beta)/2$. Hence 
\begin{equation}
\delta c =-\sin\frac{\alpha-\beta}{2}\,\delta \ell
\label{DeviationsBisectorProb2}
\end{equation}

\vspace{2.5mm}
\noindent\parbox{86mm}{
\includegraphics[width=86mm,keepaspectratio]{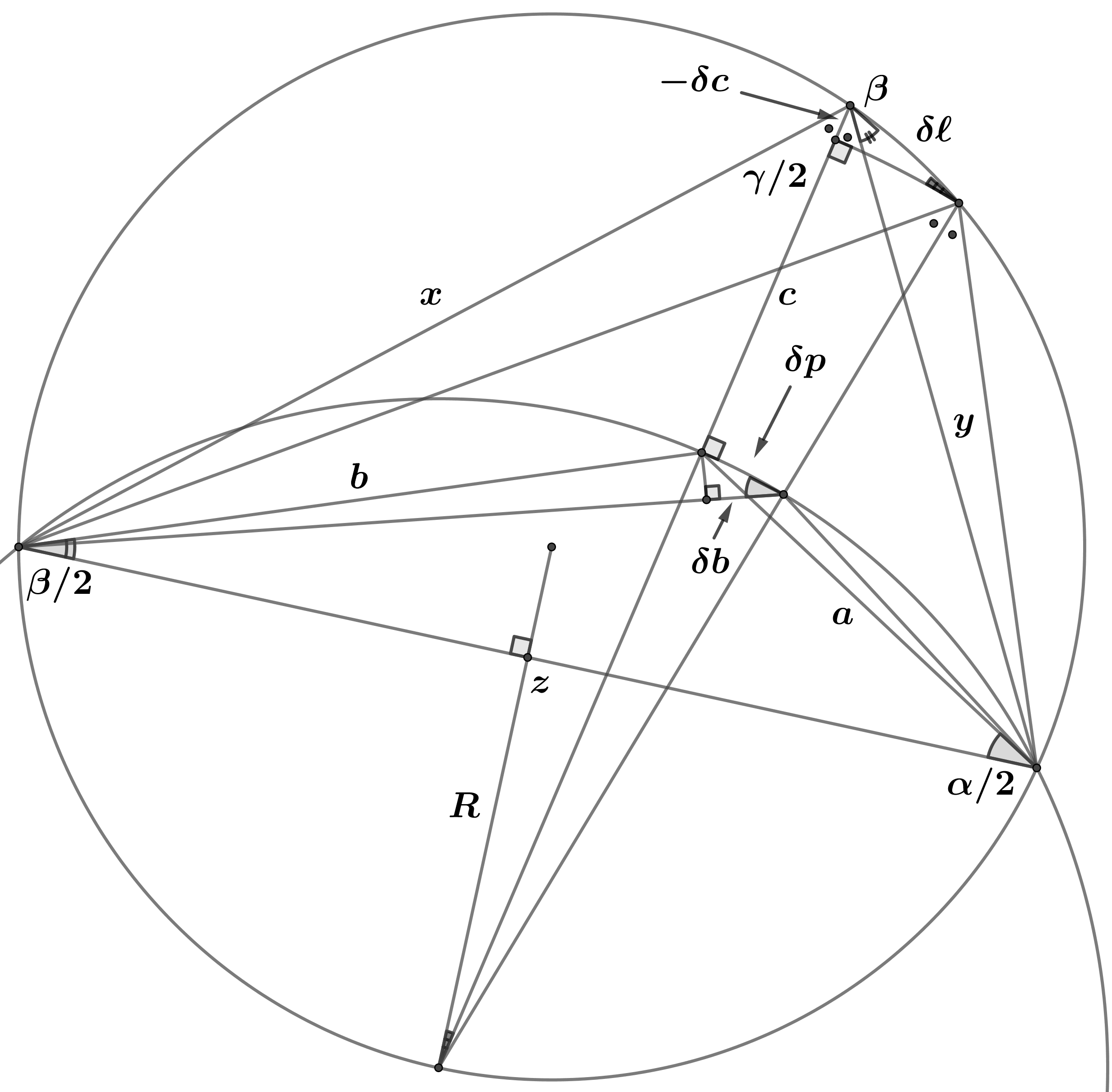}%
}%
\vspace{2.5mm}

\noindent To connect $\delta \ell$ to $\delta p$, just consider $\delta \ell$ as the product of the circumradius length $R$ by the corresponding central angle, that is twice the inscribed angle $\delta\alpha=-\delta\beta$ (since $\delta \gamma=0$). This inscribed angle $\delta\alpha$ is itself connected to $\delta p$ by trigonometry: $a\,\delta\alpha/2$ is the second leglength of the small right triangle of first leglength $-\delta a$ and hypotenuse $\delta p$. But the circumradius length is not supposed to appear in our equations; we can get rid of it in favour of $z$ by reminding that $\sin\gamma=z/(2R)$, as can be seen by moving the $(x,y)$ vertex along the circumcircle, thus preserving $\gamma$, until the $x$ or $y$-length side coincides with a diameter. These considerations can be summarized as follows:
\begin{equation}
\delta \ell = 2 R \,\delta \alpha
\qquad
 a\,\frac{\delta\alpha}{2}=\sin\frac{\beta}{2}\,\delta p
\qquad
2R=\frac{z}{\sin\gamma}
\end{equation}

\noindent Inserting those equations in eq. (\ref{DeviationsBisectorProb2}) yields
\begin{equation}
\delta c =-\sin\frac{\alpha-\beta}{2}\,\frac{z}{\sin\gamma}\,\frac{2}{a}\sin\frac{\beta}{2}\,\delta p
\label{DeviationsBisectorProb3}
\end{equation}

We need a little trigonometry to move forward. Applying the law of sines in at least two of the three triangles generated by the $a$, $b$ and $c$-length segments in the initial triangle, we find
\begin{equation}
a\sin\frac{\alpha}{2}=b\sin\frac{\beta}{2}=c\sin\frac{\gamma}{2}
\label{SinBisectorProb1}
\end{equation}

\noindent Furthermore, if we apply the al-Kashi theorem to the angle $\widehat{(a,b)}=\pi-(\alpha+\beta)/2=\pi/2+\gamma/2$, we have
\begin{equation}
\sin\frac{\gamma}{2}=\frac{z^2-a^2-b^2}{2ab}
\label{SinBisectorProb2}
\end{equation}

\noindent Inserting this in eq. (\ref{SinBisectorProb1}) gives
\begin{equation}
\sin\frac{\alpha}{2}=\frac{c}{a}\,\frac{z^2-a^2-b^2}{2ab}
\quad
\sin\frac{\beta}{2}=\frac{c}{b}\,\frac{z^2-a^2-b^2}{2ab}
\label{SinBisectorProb3}
\end{equation}

\noindent Applying the al-Kashi theorem to the angles $\widehat{(a,z)}=\alpha/2$ and $\widehat{(b,z)}=\beta/2$ respectively, we find 
\begin{equation}
\cos\frac{\alpha}{2}=\frac{z^2+a^2-b^2}{2az}
\qquad
\cos\frac{\beta}{2}=\frac{z^2-a^2+b^2}{2bz}
\label{CosBisectorProb1}
\end{equation}

\noindent Inserting those results in the development of 
$ \cos(\gamma/2)=\sin(\alpha/2+\beta/2)
=\sin(\alpha/2)\cos(\beta/2)+\cos(\alpha/2)\sin(\beta/2)$
yields
\vspace{-5mm}
\begin{equation}
\cos\frac{\gamma}{2}
=\frac{cz}{ab}\,\frac{z^2-a^2-b^2}{2ab}
\label{CosBisectorProb2}
\end{equation}

\noindent Note that submitting $\sin(\gamma/2)$ in eq. (\ref{SinBisectorProb2}) and $\cos(\gamma/2)$ in eq. (\ref{CosBisectorProb2}) to the fundamental law leads to a polynomial equation of degree $3$ in $z^2$, that is, the complicated but existing solution to our problem. 

But in this quest, we were following Newton on another path: it is now time to exploit the expressions of the small deviations in eq. (\ref{DeviationsBisectorProb1}) and (\ref{DeviationsBisectorProb3}) -- developing $\sin(\alpha/2-\beta/2)$ and using $\sin\gamma=2\sin(\gamma/2)\cos(\gamma/2)$ in the latter -- and plug in them into eq. (\ref{PartialBisectorProb0}) to finally obtain our first partial differential equation
\begin{equation}
\frac{z^2-a^2+b^2}{2b}\,\partial_a z
-\frac{z^2+a^2-b^2}{2a}\,\partial_b z
+c\,\frac{b^2-a^2}{ab}\,\partial_c z=0
\label{PartialBisectorProb1}
\end{equation}

\noindent Combining it with the differential equation arising from the infinitesimal scale transformation $R\mapsto R+\delta R$, that is
\begin{equation}
z=a\,\partial_a z + b\,\partial_b z + c\,\partial_c z
\label{ScaleBisectorProb}
\end{equation}

\noindent and isolating $\partial_a z$ and $\partial_b z$, yields
\begin{align}
\partial_a z&=\frac{z^2+a^2-b^2}{2az}-\frac{c}{z}\,\frac{z^2-a^2+b^2}{2az}\,\partial_c z
\label{PartialBisectorProb2a}\\
\partial_b z&=\frac{z^2-a^2+b^2}{2bz}-\frac{c}{z}\,\frac{z^2+a^2-b^2}{2bz}\,\partial_c z
\label{PartialBisectorProb2b}
\end{align}

That was the easy part. We need a third partial differential equation to complete the system and be able to solve it. The analogue of eq. (\ref{PartialBisectorProb1}) can be obtained by considering a small displacement of the $(y,z)$ vertex along the circumcircle, instead of the $(x,y)$ one, but it is more complicated since now $\delta z\neq 0$. At first order 
\begin{equation}
\delta z=\partial_a z \,\delta a + \partial_b z \,\delta b + \partial_c z \,\delta c
\label{PartialBisectorProb3}
\end{equation}

\noindent Following a similar reasoning as previously we find
\begin{align}
\delta a &=-\sin\frac{\beta-\gamma}{2}\,\frac{z}{\sin\gamma}\,\frac{2}{b}\sin\frac{\gamma}{2}\,\delta p\\
\delta b &=-\cos\frac{\gamma}{2}\,\delta p\\
\delta c &=\cos\frac{\beta}{2}\,\delta p\\
\delta z &=-\cos\gamma\,\frac{z}{\sin\gamma}\,\frac{2}{b}\sin\frac{\gamma}{2}\,\delta p
\label{DeviationsBisectorProb4}
\end{align}

\noindent where the latter derives from $\delta z =-\cos\gamma\,\delta\ell$, by analogy with what was done, for instance, to determine $\delta y$ in eq. (\ref{DeviationsCircum0}). Using the same trigonometry as for the previous case, plus $\cos\gamma=\cos^2(\gamma/2)-\sin^2(\gamma/2)$, and plugging in the results into eq. (\ref{PartialBisectorProb3}), we have
\begin{align}
&\frac{c^2z^2-a^2b^2}{ab^2c}\,\frac{z^2-a^2-b^2}{2a}
\nonumber\\
&=\left[\frac{cz}{ab}\,\frac{z^2-a^2-b^2}{2b}-\frac{ab}{cz}\,\frac{z^2-a^2+b^2}{2b}\right]\,\partial_a z
\nonumber\\
&+\frac{cz}{ab}\,\frac{z^2-a^2-b^2}{2a}\,\partial_b z
-\frac{z^2-a^2+b^2}{2z}\,\partial_c z
\label{PartialBisectorProb4}
\end{align}

\noindent All that remains is to replace, in this equation, $\partial_a z$ and $\partial_a z$ by their expression in eq. (\ref{PartialBisectorProb2a}) and (\ref{PartialBisectorProb2b}) respectively. After a tedious calculation we finally obtain the differential equation
\begin{align}
&c\,\partial_c z=\\
&\frac{z(z^2-a^2-b^2)(-z^4+2(a^2+b^2)z^2-(a^2-b^2)^2)}{z^6-3(a^2+b^2)z^4+3(a^2-b^2)^2z^2-(a^2+b^2)(a^2-b^2)^2 }
\nonumber
\end{align}

\noindent which
can be integrated out to give the general solution
\begin{equation}
\frac{\sqrt{((a+b)^2-z^2)(z^2-(a-b)^2)}}{z(z^2-a^2-b^2)}=k(a,b)\, c
\label{GeneralBisectorProb}
\end{equation}

\noindent where $k(a,b)$ is a symmetric function of $a$ and $b$. Here $z(a,b,c)$ si implicitely given as a root of a polynomial equation of degree $3$ in $z^2$... And this equation happens to be the same as the one we have mentionned earlier if we take $k(a,b)=\pm 1/(ab)$.

\section{The Polynomials} 
\label{Polynomials}

\noindent One of the lessons that we can draw from this last attempt is that the method cannot evade the sometimes complicated question of finding the roots of a polynomial. Newton surely was interested in polynomials. Could he have tried to find their roots through calculus? Let us first look at a polynomial $P(x)$ of degree $2$ with real coefficients $a$, $b$ and $c$. Each of its roots, if any, satisfies the equation 
\begin{equation}
ax^2+bx+c=0
\label{Polynomial2Roots}
\end{equation}

\noindent Assume that the root $x$ of our polynomial is a smooth function of the coefficients: 
\begin{equation}
x=x(a,b,c)
\label{Polynomial2Function}
\end{equation}

\noindent After an infinitesimal deviation of those coefficients, we have, at first order:  
\begin{equation}
\delta x=\partial_a x \,\delta a + \partial_b x \,\delta b + \partial_c x \,\delta c 
\label{PartialPolynomial2}
\end{equation}

\noindent But eq. (\ref{Polynomial2Roots}) must also remain valid, so that, at first order,
\begin{equation}
x^2\delta a + x\delta b +\delta c=-(2ax+b)\delta x
\end{equation}

\noindent If we consider first $\delta c\neq 0$ with $\delta a=\delta b=0$, then $\delta b\neq 0$ with $\delta a=\delta c=0$, and finally $\delta a\neq 0$ with $\delta b=\delta c=0$, we have repsectively
\begin{equation*}
\delta c=-(2ax+b)\delta x
\;\;
\delta b =-\frac{2ax+b}{x}\delta x
\;\;
\delta a =-\frac{2ax+b}{x^2}\delta x
\end{equation*}

\noindent that we plug in into eq. (\ref{PartialPolynomial2}) to obtain three partial differential equations
\begin{equation}
\partial_c x=\frac{-1}{2ax+b}
\quad
\partial_b x=\frac{-x}{2ax+b}
\quad
\partial_a x=\frac{-x^2}{2ax+b}
\end{equation}

\noindent The general solutions of these equations are respectively
\begin{align}
ax^2+bx+c&=k_c(a,b)
\\
(2ax+b)^2&=b^2+k_b(a,c)
\\
(2ax+b)^2-b^2&=k_a(b,c)a
\end{align}

\noindent where $k_c(a,b)$, $k_b(a,c)$, and $k_a(b,c)$ are arbitrary functions. Using the two last solutions leads to $k_b(a,c)=k_a(b,c)a= k(c) a$ where $k(c)$ is an arbitrary function; comparing this to the first solution yields $4k_c(a,b)=K$ and $k(c)=K-4c$, with $K$ a constant that must vanish since $x$ is a root of the polynomial; which brings us back to eq. (\ref{Polynomial2Roots}). 

For a polynomial of degree $n$ with real coefficients, the roots satisfy the polynomial equation
\begin{equation}
\sum_{k=0}^n a_k x^k=0
\label{PolynomialnRoots}
\end{equation}

\noindent Assume that the root $x$ is a smooth function of the coefficients: 
\begin{equation}
x=x(a_0,\ldots,a_n)
\label{PolynomialnFunction}
\end{equation}

\noindent After an infinitesimal deviation of those coefficients, we have, at first order:  
\begin{equation}
\delta x=\sum_{k=0}^n\partial_{a_k} x \,\delta a_k
\label{PartialPolynomialn0}
\end{equation}

\noindent If we replace $a_k$ by $a_k+\delta a_k$ and $x$ by $x+\delta x$ in eq. (\ref{PolynomialnRoots}), we have at first order again
\begin{equation}
\sum_{k=0}^n x^k\delta a_k =-\sum_{i=0}^n i a_i x^{i-1}\,\delta x
\end{equation}

\noindent Now, by successively keeping only one non-zero coefficient deviation $\delta a_k$ for $k=0,\ldots,n$, we find the $(n+1)$ relations 
\begin{equation}
x^k\delta a_k=-\sum_{i=0}^n i a_i x^{i-1}\,\delta x
\end{equation}

\noindent that can be plugged in into eq. (\ref{PartialPolynomialn0}) to yield the $(n+1)$ partial differential equations
\begin{equation}
\left(\sum_{i=0}^n i a_i x^{i-1}\right)\partial_{a_k}x +x^k=0
\label{PartialPolynomialn1}
\end{equation}

\noindent for $k=0,\ldots,n$. The general solution of this system should provide the roots of the polynomial! Each of these equations can be expressed as an ordinary first-order differential equation that turns out to be exact. For each $k$, we can indeed see our polynomial as a function $u_k$ of $x$ and $a_k$, i.e.
\begin{equation}
u_k(x, a_k)=\sum_{i=0}^n a_i x^i
\end{equation}

\noindent Its partial derivatives
\begin{equation}
\partial_x u_k=\sum_{i=0}^n i a_i x^{i-1}\quad\mbox{and}\quad \partial_{a_k} u_k=x^k
\end{equation}

\noindent are naturally the coefficients of eq. (\ref{PartialPolynomialn1}) so that for each $k$, we have the total differential $du_k=0$, and $u_k(x, a_k)$ must be equal to an arbitrary function of all the $a_i$ for $i\neq k$. Equating those $(n+1)$ arbitrary functions forces them to be constant and eventually to vanish since $x$ is a root of the polynomial. We end up with eq. (\ref{PolynomialnRoots}). In other words, we have found a set of partial differential equations equivalent to the polynomial equation, but the solution of this set is naturally given in the implicit form of the polynomial equation itself, as we could have expected! Newton, at least, would have\ldots

\section{Archimedes of Syracuse}
\label{Archimedes}

\noindent We don't need to speculate on Newton's aspirations, to pretend to discover that the method can be used to compute the area of a circle and the volume of a sphere as functions of their radius length $r$ by observing the way they behave under the transformation $r\mapsto r+\delta r$. 

For the circle area $A(r)$, we have $\delta A=A'\delta r$, and it is easy to calculate $\delta A=2\pi r \delta r$, so that $A'=2\pi r$ and $A(r)=\pi r^2+k$ where $k$ is a constant equal to $0$ since $A(0)=0$. 

We cannot do it for the perimeter $P(r)$ of the circle since the intermediate result $\delta P=2\pi\delta r$ that leads to $P'=2\pi$ and $P(r)=2\pi r$, precisely relies on the fact that $P=2\pi r$, that is, the definition of $\pi$. 

If $V_s(r)$ is the volume of the sphere, then $\delta V_s=V_s'\delta r$. By imparting a small thickness $\delta r$ to the spherical envelope of area $4\pi r^2$, we calculate $\delta V_s=4\pi r^2 \delta r$, so that $V_s'=4\pi r^2$ and $V_s(r)=4\pi r^3/3+k$ where $k$ is a constant equal to $0$ since $V_s(0)=0$. 

We could do it for the area of the sphere $A_s(r)$, but we did not find a way to compute $\delta A_s=8\pi r\delta r$ without using integral calculus and/or, at least, the cylindrical projection that Archimedes used to directly compute $A_s(r)=4\pi r^2$ \cite{Archimedes:225BC}.

\section{A Simple Equation} 
\label{Equation}

\noindent The equations (\ref{ScaleCircum}), (\ref{ScaleSines}), (\ref{ScalePtolemy}), (\ref{ScaleBisectorIn}), (\ref{ScaleInradius}) and (\ref{ScaleBisectorProb})  that were used to determine the circumradius, the law of sines, etc. could be formulated in a more general way -- still specific to Euclidean Geometry: 
\begin{equation}
f(x_1,...,x_p)=\frac{1}{n}\sum_{i=1}^p n_i x_i \partial_{x_i} f
\label{ScaleGeneral}
\end{equation}

\noindent where $f$ is a metric quantity (length, area, volume,\ldots) and $n$ its dimension ($1,2,3,\ldots$), while $x_i$ are the metric quantities on which $f$ depends, and $n_i$ their own dimensions. Every theorem of Euclidean Geometry expressed as a smooth function, every length, surface, volume, hypersurface formula should satisfy equation (\ref{ScaleGeneral}) which insures that all metric quantities involved in it are dimensionally correctly mixed together. It is naturally true for all the theorems re-derived in this article. 

Besides the six cases mentionned above, we show how we could have used equation (\ref{ScaleGeneral}) to derive Pythagoras and Heron's theorems. In the first subsection we observe how this equation specifically acts by restricting a large class of solutions. The second subsection is devoted to the study of a $n=2$ case, since Heron deals with the area of triangles.

\subsection{Pythagoras}

\noindent Imagine that instead of linearily increasing $x$ (and then $y$) of an infinitesimal $\delta x$, we would have chosen to perform an infinitesimal rotation of the hypotenuse around, say, the $(y,z)$ vertex. Here $\delta z=0$, we would then have
\begin{equation}
\partial_x z \,\delta x + \partial_y z \,\delta y=0
\label{PartialPyth3}
\end{equation}

\noindent where the small deviations $\delta x$ and $-\delta y$ happen to be the leglengths of a right triangle similar to the initial one -- since its hypotenuse is, at first order, perpendicular to the hypotenuse of the initial triangle --, yielding
\begin{equation}
\delta y=\frac{x}{y}\delta x
\end{equation}
\noindent\parbox{83mm}{
\includegraphics[width=43mm,keepaspectratio]{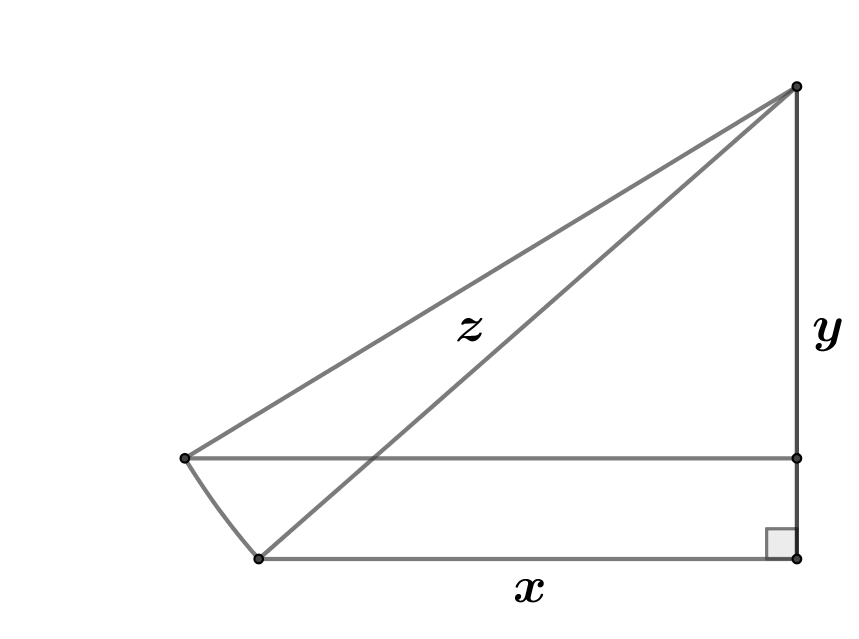}%
\includegraphics[width=43mm,keepaspectratio]{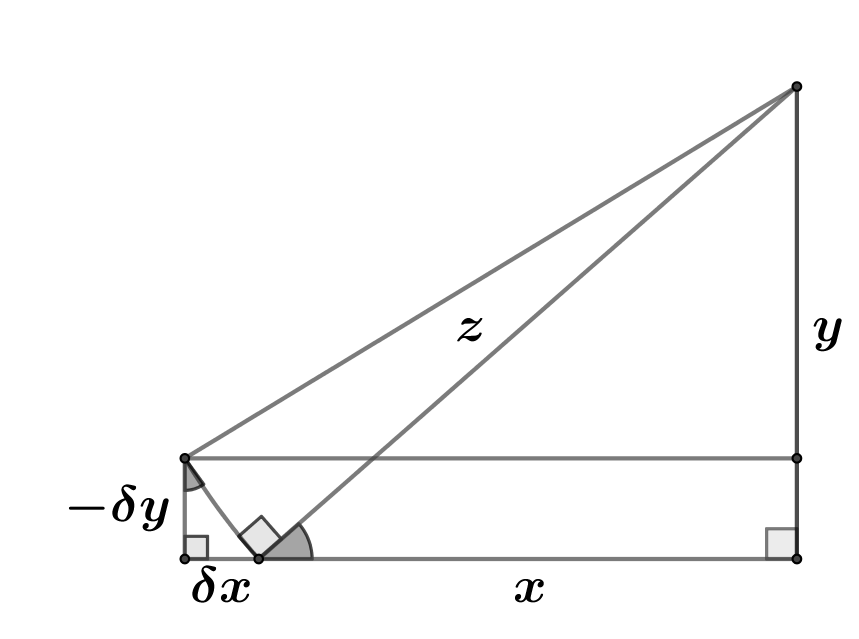}
}

\noindent Plugging in this equation into equation (\ref{PartialPyth3}) gives
\begin{equation}
y\partial_x z - x \partial_y z=0
\label{PartialPyth4}
\end{equation}

\noindent which admits a large class of solutions of the form $z(x,y)=f(x^2+y^2)$. To complete the derivation of the Pythagorean theorem, we need to consider the way $z$ behaves, at first order, under a scale transformation, that is, $R\mapsto R+\delta R$, where $R$ is the circumradius (here $R=z/2$) or any other length to be used as a scale. We have
\begin{equation}
z=x\,\partial_x z + y\,\partial_y z
\label{ScalePyth}
\end{equation}

\noindent which corresponds to equation (\ref{ScaleGeneral}). Joining equations (\ref{PartialPyth4}) and (\ref{ScalePyth}) in a system and isolating $\partial_x z$ yields
\begin{equation}
\frac{\partial_x z}{z}=\frac{x}{x^2+y^2}
\label{PartialPyth5}
\end{equation}

\noindent and its general solution
\begin{equation}
z(x,y)=k(y)\sqrt{x^2+y^2}
\end{equation}

\noindent with $k(y)$ a function of $y$. But the equation (\ref{PartialPyth5}) has its analogous for $\partial_y z$, generating by integration a function $k(x)$ which must be equal to $k(y)$, that is, to a real and positive constant $k$. Considering a triangle degenerated in a segment implies that $k=1$.

\subsection{Heron}

\noindent We still postulate that
\begin{equation}
A=A(x,y,z)
\end{equation}

\noindent but we now consider the triangle in its circumcircle. If the $(x,y)$ vertex is slightly moved along the circumcircle, it generates the infinitesimal deviations $\delta x$ and $\delta y$, while $\delta z=0$. At first order, we have: 
\begin{equation}
\delta A=\partial_x A \,\delta x + \partial_y A \,\delta y
\label{PartialHeron3}
\end{equation}

\noindent The deviations $\delta x$ and $\delta y$ are the same as in the circumcircle section, equation (\ref{DeviationsCircum1}) while $\delta A=z\delta h/2$. To determine $\delta h$, we observe that it is one of the leglengths of a small right triangle of hypotenuse of length $\delta\ell$, and that this right triangle is similar to another one of hypotenuse of length $R$ and leg with corresponding length $d_z$, the distance between the circumcenter and the height of length $h$, that is, the distance between the middle of the $z$-length side and the foot of the $h$-length height, or else

\noindent\parbox{86mm}{\parbox{43mm}{

\includegraphics[width=43mm,keepaspectratio]{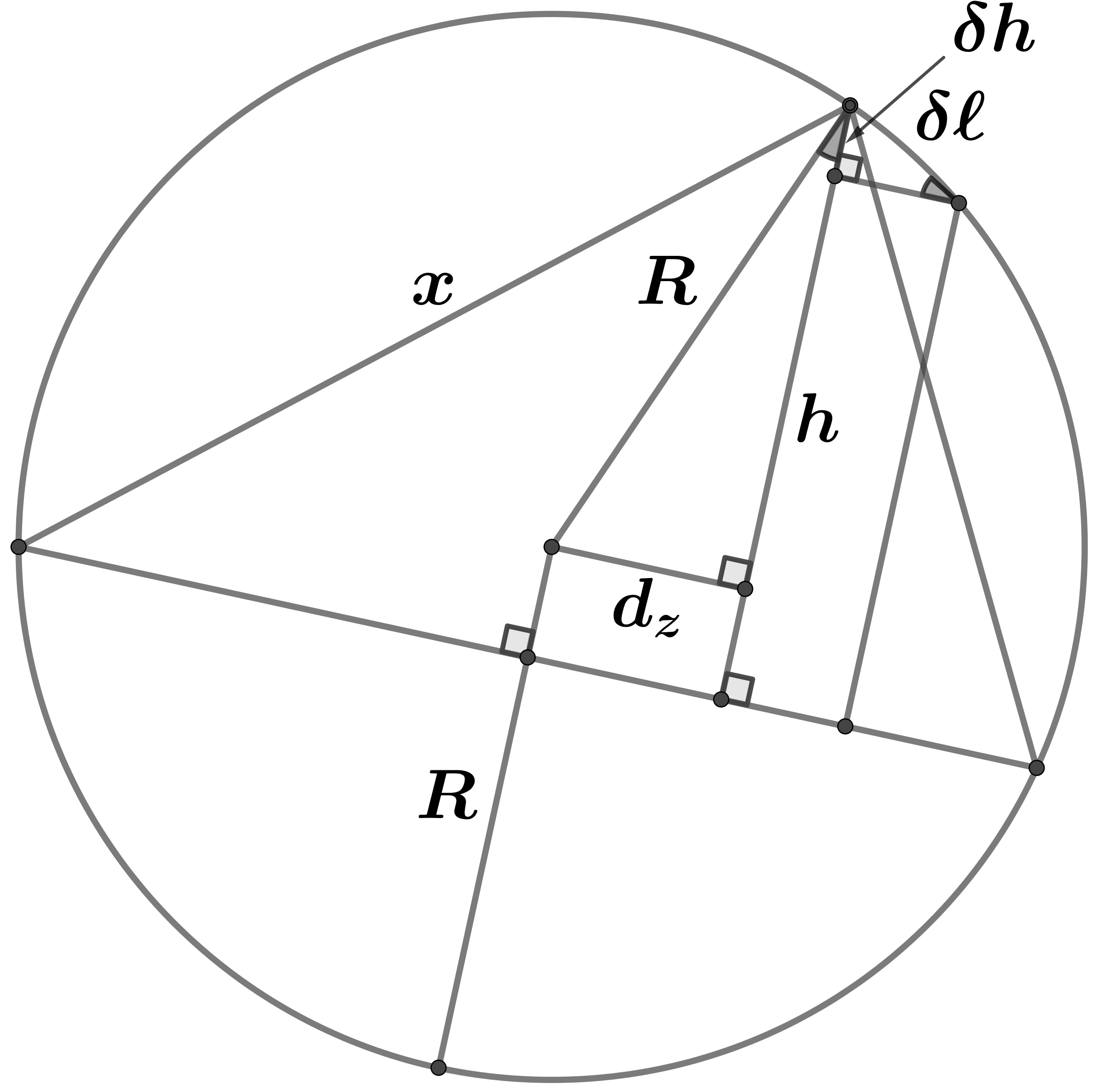}

\vspace{1mm}

\includegraphics[width=43mm,keepaspectratio]{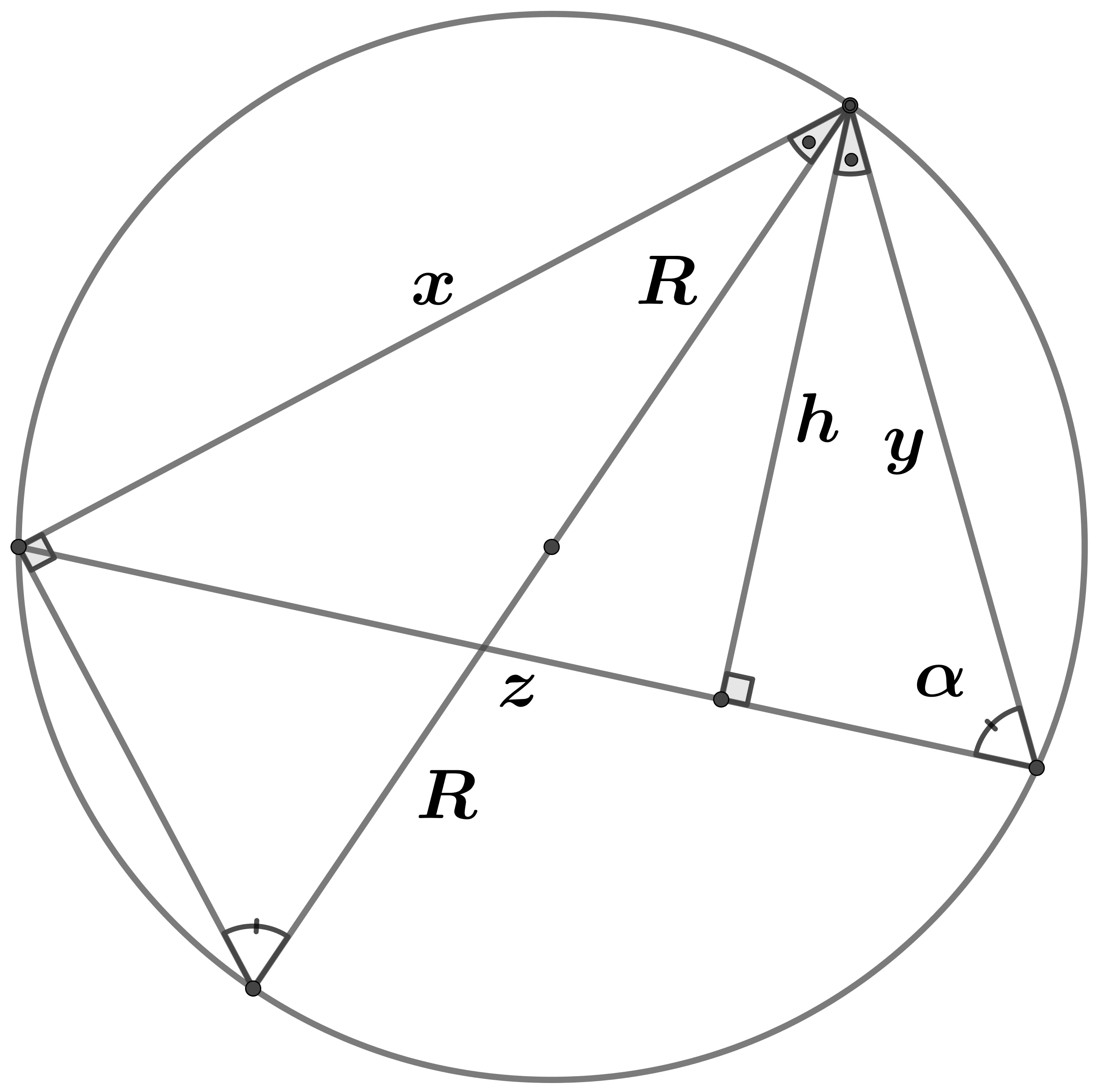}

}\hspace{2mm}%
\parbox{41mm}{\vspace{-1mm}
\begin{equation}
d_z=\frac{y^2-x^2}{2z}
\end{equation}

\noindent This can be seen by subtracting $t$, in eq. (\ref{HauteurDiviseBase}), to $z/2$. Hence
\begin{align}
\nonumber
\delta h&= \frac{d_z}{R}\delta \ell =\frac{y^2-x^2}{2zR}\delta \ell\hspace{-1.8mm}\\
\hspace{-2.4mm}\Rightarrow\hspace{-.5mm}
\delta A&=\frac{y^2-x^2}{4R}\delta \ell
\end{align}

\noindent But in the circumcircle, the angle between the $h$-length height and the $x$-length side is the same as the one between the $y$-length side and the circumradius from de $(x,y)$ vertex, thus $h/x=y/(2R)$. 
}
}

\vspace{1.5mm}
\noindent Isolating $R$ and substituting $2A/z$ to $h$ yields
\begin{equation}
R=\frac{xyz}{4A}
\end{equation}

\noindent Thus 
\begin{equation}
\delta A=A\,\frac{y^2-x^2}{xyz}\delta \ell
\label{DeviationsHeron2}
\end{equation}

\noindent Plugging in equations (\ref{DeviationsCircum1}) and (\ref{DeviationsHeron2}) into equation (\ref{PartialHeron3}), simplifying by $\delta\ell/z$ and re-arranging the terms leads to
\begin{equation}
2 A=x\,\frac{y^2-x^2+z^2}{y^2-x^2}\partial_x A  - y\,\frac{x^2-y^2+z^2}{y^2-x^2}\partial_y A 
\label{PartialHeron4}
\end{equation}

\noindent This equation can be joined to its symmetric couterpart corresponding to the infinitesimal displacement of the $(z,x)$ vertex along the circumsircle, while $\delta y=0$, that is
\begin{equation}
2 A=z\,\frac{x^2-z^2+y^2}{x^2-z^2}\partial_z A  - x\,\frac{z^2-x^2+y^2}{x^2-z^2}\partial_x A 
\label{PartialHeron5}
\end{equation}

\noindent and the equation originating from the infinitesimal scale transformation $R\mapsto R+\delta R$
\begin{equation}
2A=x\,\partial_x A + y\,\partial_y A + z\,\partial_z A
\label{ScaleHeron}
\end{equation}

\noindent in a system of three equations that we can solve in $\partial_x A$ to obtain the partial differential equation
\begin{equation}
\frac{\partial_x A}{A}=\frac{1}{2}\frac{-4x^3+4x(y^2+z^2)}{[-x^4-y^4-z^4+2(y^2 z^2 + x^2 z^2 + x^2 y^2)]}
\label{PartialHeron6}
\end{equation}

\noindent which 
can be integrated out to give the general solution
\begin{align}
&A(x,y,z)\\
&=k(y,z)\sqrt{-x^4-y^4-z^4+2(y^2 z^2 + x^2 z^2 + x^2 y^2)}
\nonumber
\end{align}

\noindent where $k(y,z)$ is an function of $y$ and $z$. Since the function $A$ must be symmetric in $x$, $y$ and $z$, $k(y,z)=c$ where $c$ is a constant. We just have to compute the area of, say, an isocele right triangle of leglength $1$, that is $A(1,1,\sqrt{2})=1/2$, to find out that $c=1/4$. After factorization, the result is equation (\ref{HeronTheorem}), that is, Heron's theorem.

\section*{Conclusion}

\noindent With this paper, we present an alternative way to derive classical theorems in Euclidean geometry. Not all theorems, of course. It does not work for theorems in discrete geometry or involving number theory, for theorems stating that this or that line cuts another at this or that point, is perpendicular or tangent to this or that circle, etc. It has to be a theorem involving an equation that defines a function, which will be seen as a particular solution of a (system of) differential equation(s). The proofs that we propose are not necessarily simpler than others. They do not evade the geometric difficulties at stake. We displace the argument of the proof into the game of infinitesimals, but it remains as geometric. 

The main advantage of this method is that the theorem does not need to be known. We start with a function, any function, and observe the way it behaves, at first order, under small deviations of some quantities. In the best case, it gives us a (system of) differential equation(s) that we can solve and, therefore, discover the theorem. We chose to use it to rediscover about 20 theorems or identitites in the long history of Euclidean geometry, but we hope it may be used to discover new theorems, perhaps in other fields of mathematics. 

\section*{Acknowledgments}

\noindent I thank Rapha{\"e}l Lefevere for a critical reading of the manuscript and a fruitful collaboration. 

The initial idea comes from a previous work in particle physics, looking for natural relations between the mixing angles and the fermion mass ratios \cite{Buysse:2002a}, supervised by Jean-Marc G{\'e}rard, whom I also thank. 

This paper is dedicated to Jacques Weyers, who disappeared last Fall. 

\bibliographystyle{prsty}
\bibliography{BibliMath}

\end{document}